\RequirePackage{fix-cm}
\documentclass[twocolumn]{svjour3}          
\smartqed  
\usepackage{mathptmx}      
\usepackage{algorithm}
\usepackage{algorithmic}
\usepackage{amsfonts}
\usepackage{amsmath}
\usepackage{amssymb}
\usepackage{color}
\usepackage{enumerate}
\usepackage{fullpage}
\usepackage{graphicx}
\usepackage[margin=0.75in]{geometry}
\usepackage{mathtools}
\usepackage{multirow}
\usepackage{stmaryrd}
\usepackage{subfig}

\usepackage[pdftex]{thumbpdf}
\usepackage[pdftex,
pdfstartview=FitH,bookmarks=true,bookmarksnumbered=true,bookmarksopen=true,hypertexnames=false,breaklinks=true,
colorlinks=true,  linkcolor=blue,anchorcolor=blue,
citecolor=blue,filecolor=blue,  menucolor=blue]{hyperref}%
\providecommand{\U}[1]{\protect\rule{.1in}{.1in}}
\usepackage{tikz}
\usetikzlibrary{patterns}

\newcommand{\R}{\mathbb{R}}
\newcommand{\domain}{D}
\newcommand{\boundary}{\partial \domain}

\newcommand{\visc}{\nu}
\newcommand{\effectivevisc}{\visc^*}
\newcommand{\permeability}{\mathbf{K}}
\newcommand{\invpermeability}{\permeability^{-1}}
\newcommand{\unitnormal}{n}
\newcommand{\normalderivative}[1]{\frac{\partial{#1}}{\partial\unitnormal}}

\newcommand{\genericfunction}{\phi}

\newcommand{\velocityspace}{V}
\newcommand{\pressurespace}{Q}

\newcommand{\fe}[1]{{#1}_h}
\newcommand{\dof}[1]{\mathbf{{#1}}}
\newcommand{\rb}[1]{{#1}_r}

\newcommand{\dual}[1]{{#1}^\star}

\newcommand{\dualitybracket}[3]{\langle {#1}, {#2} \rangle_{#3}}

\begin{document}

\title{Application of adaptive ANOVA and reduced basis methods to the stochastic Stokes-Brinkman problem
  \thanks{Bed\v{r}ich Soused\'{\i}k was partially supported by the U.S. National Science Foundation under grant DMS 1913201.}
}


\author{Kevin Williamson \and
  Heyrim Cho \and
  Bed\v{r}ich Soused\'{\i}k
}


 \institute{Kevin Williamson, Corresponding Author \at
  Department of Mathematics and Statistics, University of Maryland, Baltimore County, 1000 Hilltop Circle, Baltimore, MD 21250, USA.
  \email{kwillia1@umbc.edu}           
  \and
  Heyrim Cho \at
  Department of Mathematics, University of California Riverside, Riverside, CA 92521 
  \email{heyrim.cho@ucr.edu}
  \and
  Bed\v{r}ich Soused\'{\i}k \at
  Department of Mathematics and Statistics, University of Maryland, Baltimore County, 1000 Hilltop Circle, Baltimore, MD 21250, USA.
  \email{sousedik@umbc.edu}           
}

\date{Received: date / Accepted: date}

\maketitle

\begin{abstract}
  The Stokes-Brinkman equations model fluid flow in highly heterogeneous porous media. In this paper, we consider the numerical solution of the Stokes-Brinkman equations with stochastic permeabilities, where the permeabilities in subdomains are assumed to be independent and uniformly distributed within a known interval. We employ a truncated anchored ANOVA decomposition alongside stochastic collocation to estimate the moments of the velocity and pressure solutions. Through an adaptive procedure selecting only the most important ANOVA directions, we reduce the number of collocation points needed for accurate estimation of the statistical moments. However, for even modest stochastic dimensions, the number of collocation points remains too large to perform high-fidelity solves at each point. We use reduced basis methods to alleviate the computational burden by approximating the expensive high-fidelity solves with inexpensive approximate solutions on a low-dimensional space. We furthermore develop and analyze rigorous a posteriori error estimates for the reduced basis approximation. We apply these methods to 2D problems considering both isotropic and anisotropic permeabilities.
\end{abstract}

\section{Introduction}

The simulation of flow in porous media has numerous applications, to include reservoir simulation, nuclear waste disposal, and carbon dioxide sequestration. Such simulation is challenging for a variety of reasons. First, the domains tend to be fairly irregular, which complicates the model geometry. Second, the geologic formations consist of many varying materials, with different geologic properties. Third, there are often fractures and vugs within the domain that alter the effective permeabilities. The standard approach to modeling these types of problems is to couple Darcy's law and the Stokes equations and enforce the Beavers-Joseph-Saffman conditions along the interface~\cite{Arbogast2007,Beavers-1967-FPM,URQUIZA2008525}. The free-flow regions (fractures, vugs) are modeled using Stokes flow, whereas the porous region is modeled using Darcy's law~\cite{Darcy-1856-FPM,Whitaker-1986-FPM}. However, the two types of domains are not well-separated in reservoirs, and it may be difficult to determine the appropriate conditions to enforce along the interface. We model flow in porous media using the Stokes-Brinkman equations~\cite{Bear-1972-PM,Brinkman-1948-SBE,Gulbransen-2010-MMF,Laptev-2003-PPM,Popov-2007-MMM}, which combine the Stokes equations and Darcy's law into a single system of equations. The Stokes-Brinkman equations reduce to Stokes or Darcy flow depending upon the coefficients and were suggested as a replacement for the coupled Stokes-Darcy equations in~\cite{Popov-2007-MMM}. By careful selection of coefficients, the equations allow modeling of free-flow and porous domains together, thereby resolving issues along the interface.

In this paper, we consider the case where the exact permeabilities are unknown but are instead specified by a known probability distribution. The first step in many methods to solve such stochastic PDEs is to parametrize the distribution by a finite number of parameters. In some instances, such as the case studied in this paper, such a parametrization arises naturally from a partitioning of the domain into a finite number of subdomains. In other instances, a truncated Karhuen-Lo\`{e}ve expansion~\cite{Ghanem-1991-SFE,Papoulis-1991-PRS} may be used to obtain such a finite parametrization. Once the stochastic space has been parametrized, statistical moments may be approximated through stochastic collocation~\cite{Xiu-SCM-2005} in which deterministic solves are performed at carefully selected points.

The straightforward stochastic collocation method suffers from the curse of dimensionality--complexity grows exponentially with the number of parameters. To reduce the computational complexity of collocation methods, one may utilize Smolyak sparse grid~\cite{Nobile-SSG-2008} or an ANOVA decomposition~\cite{Cao-2009-ANOVA,Fisher-1925-SMR,Foo-2010-MPC}. We consider the latter in this paper. By using a truncated ANOVA decomposition, we may resolve the curse of dimensionality by decomposing the intractable high-dimensional problem into a set of tractable problems of low stochastic dimension which we solve by stochastic collocation. We obtain further improvements by adaptively selecting the ANOVA terms~\cite{Ma-2010-AHS,Yang-2012-AAD}.

Whereas a truncated ANOVA decomposition or Smolyak sparse grid significantly reduces the number of collocation points, the number of remaining collocation points might still be computationally prohibitive if the cost of solving at a single collocation point is high. Such is the case with the discrete Stokes-Brinkman systems, which are highly ill-conditioned due to the heterogeneity in permeabilities. Reduced basis methods~\cite{Boyaval-2010-RBT,Haasdonk-2008-RBM,Hesthaven-2016-CRB,Quarteroni-2016-RBM} alleviate this computational burden by approximating the manifold of solutions by a low-dimensional linear space. An expensive offline step is first performed to generate the reduced basis. However, after the reduced basis is built, all subsequent solves at collocation points are replaced by cheap low-dimensional computations on the reduced basis.

Reduced basis methods in conjunction with ANOVA for parametric partial differential equations were first studied in~\cite{Hesthaven-2016-UAE} where a three step RB-ANOVA-RB method was suggested for solving high-dimensional problems. The first two steps (RB-ANOVA) are similar to the idea presented in this paper in which a reduced basis is built and then used to compute the ANOVA expansion. However, in our work, the construction of the reduced basis and the ANOVA expansion occur simultaneously rather than as separate steps. Furthermore, in~\cite{Hesthaven-2016-UAE}, only a first-level ANOVA expansion is considered for the purpose of identifying parameters which have low sensitivity and can be given fixed values, thereby reducing the overall parametric dimension for the reduced basis solve in the final (RB) step. Instead, we adaptively select the ANOVA terms to allow higher-level terms at minimal cost. Key to the success of our method is the use of anchored ANOVA~\cite{Ma-2010-AHS,Xu-2004-GDR,Yang-2012-AAD,Zhang-2012-EEA,Zhang-2011-APM} which allows us to replace the high-dimensional integrals of the ANOVA expansion with cheaper low-dimensional integrals. In this way, our work is more closely aligned with the work in~\cite{Cho-2017-RBA,Liao-2016-RBA} where reduced basis methods are used in conjunction with anchored ANOVA and adaptive selection of ANOVA terms to solve high-dimensional stochastic partial differential equations.

In this paper, we apply these methods to the stochastic Stokes-Brinkman problem. As a saddle-point problem, the Stokes-Brinkman equations introduce additional complexity to the reduced basis methods because the reduced basis systems must be inf-sup stable. To guarantee inf-sup stability, we follow the approach of~\cite{Rozza-2007-SRB}, which considers the simpler Stokes case. We also devise rigorous a posteriori error estimates based upon the Brezzi stability theory, following the approach of~\cite{Gerner-2012-CRB}. These a posteriori error estimates enable us to be confident in the accuracy of the reduced basis approximations and are useful in building the reduced basis.

This paper is organized as follows. In Section~\ref{sec:stokes-brinkman}, we introduce the Stokes-Brinkman equations, present their discretization using mixed finite elements, and describe the parametrization of the stochastic permeabilities. In Section~\ref{sec:anova}, we discuss the ANOVA decomposition and the adaptive selection of terms. An overview of reduced basis methods is presented in Section~\ref{sec:reduced-basis}, and the rigorous a posteriori error estimates upon which these methods rely are presented in Section~\ref{sec:aee}. Finally, in Section~\ref{sec:numerical}, we present numerical experiments demonstrating the effectiveness of these techniques for stochastic Stokes-Brinkman problems, considering problems with both isotropic and anisotropic permeabilities.

\section{The Stokes-Brinkman Problem and its Discretization} \label{sec:stokes-brinkman}

\subsection{The Stokes-Brinkman Equations}

Let~$\domain \subset \R^d$,~$d = 2, 3$, be a connected, open domain with Lipschitz boundary~$\boundary$. The Stokes-Brinkman equations~\cite{Bear-1972-PM,Brinkman-1948-SBE,Gulbransen-2010-MMF,Laptev-2003-PPM,Popov-2007-MMM} model the flow of a viscous fluid in heterogeneous porous material as

\begin{align}
  -\effectivevisc \Delta u + \visc \invpermeability u + \nabla p &= f \label{eq:sb-momentum} \\
  \nabla \cdot u &= 0, \label{eq:sb-mass}
\end{align}

\noindent where~$\visc > 0$ is the constant viscosity of the fluid,~$\effectivevisc > 0$ is an effective viscosity,~$\permeability$ is a symmetric positive definite permeability tensor,~$u$ is the velocity,~$p$ is the pressure, and~$f$ denotes external forces. In this paper, we restrict attention to the case where~$\permeability$ is a diagonal matrix. We further require that~$\effectivevisc$ and~$\invpermeability$ are bounded above on~$\domain$. Equation~(\ref{eq:sb-momentum}) is derived from the conservation of momentum, and equation~(\ref{eq:sb-mass}) is derived from the conservation of mass. They are accompanied by the following Dirichlet and Neumann boundary conditions

\begin{align*}
  u &= u_D \quad \textrm{ on } \boundary_D \\
  \effectivevisc \normalderivative{u} - p\unitnormal &= u_N \quad \textrm{ on } \boundary_N,
\end{align*}

\noindent where~$\boundary_D$ and~$\boundary_N$ denote the Dirichlet and Neumann boundaries, respectively,~$\unitnormal$ is the unit outward normal, and~$\normalderivative{u}$ is the directional derivative of the velocity in the normal direction.

The Stokes-Brinkman equations may be understood as limiting cases of the Stokes equations and Darcy's law~\cite{Darcy-1856-FPM,Whitaker-1986-FPM}. Indeed, if~$\permeability \gg 0$ and~$\effectivevisc = \visc$, equation~(\ref{eq:sb-momentum}) approximates the Stokes equation~$-\visc \Delta u + \nabla p = f$, and, as~$\effectivevisc \rightarrow 0$, it approximates Darcy's law~$u = -\frac{\permeability}{\visc} (\nabla p - f)$. Thus, the Stokes-Brinkman equations provide a single system of equations to solve in highly heterogeneous porous media. In general, the choice of~$\effectivevisc$ will depend upon the porosity of the material. However, for small permeabilities, the diffusive term~$\effectivevisc \Delta u$ introduces only a small perturbation to Darcy's law~\cite{Laptev-2003-PPM}. In the absence of precise information about the value of~$\effectivevisc$ in regions of small permeability, it is common to choose~$\effectivevisc = \visc$, which is the convention we adopt in this paper. For more details on the choice of~$\effectivevisc$, see~\cite[pp. 26--29]{Laptev-2003-PPM}.

\subsection{Finite Element Discretization}

We seek a weak solution to equations~(\ref{eq:sb-momentum})--(\ref{eq:sb-mass}). Let~$L^2(\domain)$ denote the space of square-integrable functions on~$\domain$,

\begin{align*}
  L^2(\domain) &\equiv \{q \colon \int_\domain q^2 < \infty\},
\end{align*}

\noindent and~$H^1(\domain)$ denote the subspace of~$L^2(\domain)$ with weak derivatives in~$L^2(\domain)$,

\begin{align*}
  H^1(\domain) &\equiv \{q \in L^2(\domain) \colon \frac{\partial q}{\partial x_i} \in L^2(\domain), i = 1,\ldots,d\} .
\end{align*}

\noindent The space~$H^1(\domain)^d$ denotes the space of vector-valued functions whose~$d$ components are each in~$H^1(\domain)$. We define the following spaces for the velocity

\begin{align}
  H^1_E(\domain) &\equiv \{ u \in H^1(\domain)^d \colon u = u_D \textrm{ on } \boundary_D\} \nonumber \\
  H^1_{E_0}(\domain) &\equiv \{ u \in H^1(\domain)^d \colon u = 0 \textrm{ on } \boundary_D\}. \nonumber
\end{align}

\noindent We choose a fixed~$w \in H^1_E(\domain)$ and note that any~$v \in H^1_E(\domain)$ may be written uniquely as~$u + w$ for some~$u \in H^1_{E_0}(\domain)$.

We define the bilinear forms

\begin{align}
  a_S(u, v) &\equiv \int_{\domain} \effectivevisc \nabla u \colon \nabla v, \quad u, v, \in H^1(\domain)^d \label{eq:bilinear-stokes} \\
  a_D(u, v) &\equiv \int_{\domain} u \cdot (\invpermeability v), \quad u, v \in H^1(\domain)^d \label{eq:bilinear-darcy} \\
  a(u, v) &\equiv a_S(u,v) + a_D(u,v), \quad u, v \in H^1(\domain)^d \label{eq:bilinear-a} \\
  b(v, q) &\equiv -\int_{\domain} q \, \nabla \cdot v, \quad v \in H^1(\domain)^d, q \in L^2(\domain)\label{eq:bilinear-b}
\end{align}

\noindent where~$\nabla u \colon \nabla v = \sum_{i=1}^d \nabla u_i \cdot \nabla v_i$. The bilinear forms~$a_S(\cdot,\cdot)$ and~$a_D(\cdot,\cdot)$ denote the Stokes and Darcy parts, respectively, of the bilinear form~$a(\cdot,\cdot)$. We define the linear functionals

\begin{align}
  \ell_1(v) &\equiv \int_{\domain} f \cdot v - a(w,v) + \int_{\boundary_N} u_N \cdot v, \quad v \in H^1_{E_0}(\domain), \label{eq:rhs-l1-continuous} \\
  \ell_2(q) &\equiv -b(w,q), \quad q \in L^2(\domain). \label{eq:rhs-l2-continuous}
\end{align}

\noindent The weak formulation of~(\ref{eq:sb-momentum})--(\ref{eq:sb-mass}) is then to find~$u \in H^1_{E_0}(\domain)$ and~$p \in L^2(\domain)$ such that

\begin{align}
  a(u, v) + b(v, p) &= \ell_1(v), \quad \forall v \in H^1_{E_0}(\domain) \label{eq:sb-weak-momentum} \\
  b(u, q) &= \ell_2(q), \quad \forall q \in L^2(\domain). \label{eq:sb-weak-mass}
\end{align}

\noindent The full velocity solution with proper Dirichlet boundary conditions is then~$u + w$.

In the following, we denote the space~$H^1_{E_0}(\domain)$ by~$\velocityspace$ and the space~$L^2(\domain)$ by~$\pressurespace$. Furthermore, due to the presence of zero Dirichlet boundary conditions, we define the~$\velocityspace$ inner product as~$(u,v)_\velocityspace \equiv \int_\domain \nabla u \colon \nabla v$.

The existence and uniqueness of solutions to~(\ref{eq:sb-weak-momentum})--(\ref{eq:sb-weak-mass}) is a consequence of the Brezzi stability conditions~\cite{Boffi-2013-MFE,Brezzi-1991-MHF}. We state these conditions in Theorem~\ref{thm:stability}, paraphrasing Corollary 4.2.1 of~\cite{Boffi-2013-MFE} for the case of a symmetric~$a(\cdot,\cdot)$. This theorem additionally provides stability estimates which we will use to derive our a posteriori error estimates in Section~\ref{sec:aee}.

\begin{theorem} \label{thm:stability}
  Let~$\velocityspace$ and~$\pressurespace$ be Hilbert Spaces. Let~$a \colon \velocityspace \times \velocityspace \rightarrow \R$ be a continuous symmetric bilinear form that satisfies the coercivity condition: there exists an~$\alpha > 0$ such that~$a(u,u) \geq \alpha \|u\|_\velocityspace$ for all~$u \in \velocityspace$. Furthermore, let~$b \colon \velocityspace \times \pressurespace \rightarrow \R$ be a continuous bilinear form that satisfies the inf-sup condition: there exists~$\beta > 0$ such that

  \begin{equation}
    \inf_{q \in \pressurespace} \sup_{v \in \velocityspace} \frac{b(v,q)}{\|q\|_\pressurespace \|v\|_\velocityspace} \geq \beta. \label{eq:inf-sup-condition}
  \end{equation}

  \noindent Then, for any~$\ell_1 \in \dual{\velocityspace}$ and~$\ell_2 \in \dual{\pressurespace}$, the saddle-point system~(\ref{eq:sb-weak-momentum})--(\ref{eq:sb-weak-mass}) has a unique solution~$(u,p) \in \velocityspace \times \pressurespace$ which satisfies the stability bounds

  \begin{eqnarray}
    \|u\|_\velocityspace &\leq \frac{1}{\alpha}\|\ell_1\|_{\dual{\velocityspace}} + \frac{2}{\beta}\sqrt{\frac{\gamma}{\alpha}} \|\ell_2\|_{\dual{\pressurespace}} \label{eq:stability-bounds-u} \\
    \|p\|_\pressurespace &\leq \frac{2}{\beta}\sqrt{\frac{\gamma}{\alpha}} \|\ell_1\|_{\dual{\velocityspace}} + \frac{\gamma}{\beta^2} \|\ell_2\|_{\dual{\pressurespace}} \label{eq:stability-bounds-p}
  \end{eqnarray}

  \noindent where~$\alpha$ and~$\beta$ are the above coercivity and inf-sup constants and~$\gamma > 0$ is the continuous constant satisfying~$a(u,v) \leq \gamma \|u\|_\velocityspace \|v\|_\velocityspace$ for all~$u,v \in \velocityspace$.
\end{theorem}

Coercivity and continuity of~(\ref{eq:bilinear-stokes}) and continuity of~(\ref{eq:bilinear-b}) are immediately seen to be satisfied as these are the forms used in the Stokes equations. The fact that~$\invpermeability$ is symmetric positive definite implies that~$a(u,u) = a_S(u,u) + a_D(u,u) > a_S(u,u)$ and coercivity of~(\ref{eq:bilinear-a}) immediately follows. Continuity of~(\ref{eq:bilinear-a}) follows as the bilinear forms~$a_S(\cdot,\cdot)$ and~$a_D(\cdot,\cdot)$ are both bounded due to boundedness of the viscosities and inverse permeabilities. The inf-sup condition is more delicate, but the spaces~$\velocityspace$ and~$\pressurespace$ above were chosen to satisfy this property~\cite[Chapter~3]{Elman-2014-FEF}.

\begin{remark}
  In the case where no Neumann boundary conditions are specified, the inf-sup condition~(\ref{eq:inf-sup-condition}) fails as~$b(q,v) = 0$ for all~$v \in \velocityspace$ where~$q$ is constant. In this case, we may define~$\pressurespace$ to be the quotient space~$L^2(\domain) / \R$ in which two functions~$q_1,q_2 \in L^2(\domain)$ are identified if they differ by a constant function. By equipping this space with the norm~$\|q\|_\pressurespace = \|q - 1/|\domain| \int_{\domain} q\|_{L^2(\domain)}$, the inf-sup condition is satisfied. This demonstrates that the pressure solution in this case is unique up to a constant. In this paper, we will always specify Neumann conditions on a portion of the boundary.
\end{remark}

We solve~(\ref{eq:sb-weak-momentum})--(\ref{eq:sb-weak-mass}) using the mixed finite element method~\cite{Boffi-2013-MFE}. Given conforming finite element spaces~$\fe{\velocityspace} \subset \velocityspace$ and~$\fe{\pressurespace} \subset \pressurespace$, we seek~$(\fe{u},\fe{p}) \in \fe{\velocityspace} \times \fe{\pressurespace}$ such that

\begin{align}
  a(\fe{u},\fe{v}) + b(\fe{v},\fe{p}) &= \ell_1(\fe{v}), \quad \forall \fe{v} \in \fe{\velocityspace}, \nonumber \\
  b(\fe{u},\fe{q}) &= \ell_2(\fe{q}), \quad \forall \fe{q} \in \fe{\pressurespace}. \nonumber
\end{align}

\noindent By selecting finite element bases for~$\fe{\velocityspace}$ and~$\fe{\pressurespace}$, we obtain the discrete saddle-point problem

\begin{equation}
  \left[
    \begin{array}{cc}
      \dof{A} & \dof{B}^T \\
      \dof{B} & \dof{0}
    \end{array}
    \right]
  \left[
    \begin{array}{c}
      \dof{u} \\
      \dof{p}
    \end{array}
    \right]
    =
    \left[
      \begin{array}{c}
        \dof{f} \\
        \dof{g}
      \end{array}
      \right], \label{saddle-point}
\end{equation}

\noindent where~$\dof{A}$ and~$\dof{B}$ are discrete analogs of~(\ref{eq:bilinear-a}) and~(\ref{eq:bilinear-b}), respectively,~$\dof{u}$ and~$\dof{p}$ are the vectors of degrees of freedom for the velocity and pressure, respectively, and~$\dof{f}$ and~$\dof{g}$ are the discretizations of the functionals~(\ref{eq:rhs-l1-continuous}) and~(\ref{eq:rhs-l2-continuous}), respectively.

The coercivity and continuity properties are immediately inherited by the finite element discretization. However, care must be taken to ensure that inf-sup stability is maintained. In this paper, we obtain inf-sup stability by using the~$Q_2-P_{-1}$ approximation in which the domain is discretized into shape-regular quadrilaterals on which the velocity space is piecewise continuous biquadratic and the pressure space is piecewise discontinuous linear~\cite[Chapter~3]{Elman-2014-FEF}.

Throughout this paper, we will use the notation~$\dof{M}_H$ to denote the mass matrix for a finite-dimensional Hilbert space~$H$.

\subsection{The Stochastic Stokes-Brinkman Problem}

Let~$(\Omega, \Sigma, \mathcal{P})$ be a complete probability space with sample space~$\Omega$,~$\sigma$-algebra~$\Sigma$, and probability measure~$\mathcal{P}$. We assume stochasticity in the permeability tensor~$\permeability(\omega)$ for~$\omega \in \Omega$. For simplicity, we focus only on stochastic permeability, although stochastic viscosity, forcing term, or boundary conditions may be treated in a similar manner.

We assume that~$\omega$ may be well-approximated by a finite number of random variables~$\xi = (\xi_1,\ldots,\xi_M)^T$ with~$\xi_m \in \Gamma_m \equiv [a_m,b_m]$ and~$\xi \in \Gamma \equiv \prod_{m=1}^M \Gamma_m$. This may arise, for example, from a truncated Karhuen-Lo\`{e}ve expansion~\cite{Ghanem-1991-SFE,Papoulis-1991-PRS} or from a partitioning of the domain into subdomains. When applying the reduced basis methods discussed in Section~\ref{sec:reduced-basis}, it will be convenient if the parametrization admits an affine decomposition for the bilinear form~(\ref{eq:bilinear-a}) and the linear functional~(\ref{eq:rhs-l1-continuous}) such that

\begin{align*}
  a(u,v;\xi) &= \sum_{i=1}^{n_A} \theta_i^A(\xi) a_i(u,v) \\
  \ell_1(v;\xi) &= \sum_{i=1}^{n_f} \theta_i^f(\xi) f_i(v),
\end{align*}

\noindent where~$\{a_i\}_{i=1}^{n_A}$ are parameter-independent bilinear forms, $\{f_i\}_{i=1}^{n_f}$ are parameter-independent linear functionals, and~$\{\theta_i^A\}_{i=1}^{n_A}$ and~$\{\theta_i^f\}_{i=1}^{n_f}$ are functions mapping the parameter~$\xi$ to the coefficients of the affine decomposition. We note that the bilinear form~(\ref{eq:bilinear-b}) and the linear functional~(\ref{eq:rhs-l2-continuous}) are independent of the permeability and are, therefore, parameter-independent.

Assuming such an affine decomposition, the saddle-point problem~(\ref{saddle-point}) is now parametrized as

\begin{equation}
  \left[
    \begin{array}{cc}
      \dof{A}(\xi) & \dof{B}^T \\
      \dof{B} & \dof{0}
    \end{array}
    \right]
  \left[
    \begin{array}{c}
      \dof{u}(\xi) \\
      \dof{p}(\xi)
    \end{array}
    \right]
    =
    \left[
      \begin{array}{c}
        \dof{f}(\xi) \\
        \dof{g}
      \end{array}
      \right]. \label{eq:saddle-point-parametrized}
\end{equation}

\noindent The matrix~$\dof{A}(\xi)$ exhibits an affine decomposition

\begin{equation}
  \dof{A}(\xi) = \sum_{i=1}^{n_A} \theta_i^A(\xi) \dof{A}_i, \label{eq:affine-decomposition-A}
\end{equation}

\noindent where~$\{\dof{A}_i\}_{i=1}^{n_A}$ are discretizations of the parameter-independent bilinear forms~$\{a_i\}_{i=1}^{n_A}$. Similarly,~$\dof{f}(\xi)$ exhibits an affine decomposition

\begin{equation}
  \dof{f}(\xi) = \sum_{i=1}^{n_f} \theta_i^f(\xi) \dof{f}_i, \label{eq:affine-decomposition-f}
\end{equation}

\noindent where~$\{\dof{f}_i\}_{i=1}^{n_f}$ are discretizations of the parameter-independent linear functionals~$\{f_i\}_{i=1}^{n_f}$. The affine decompositions~(\ref{eq:affine-decomposition-A})--(\ref{eq:affine-decomposition-f}) greatly simplify the use of reduced basis methods, as is discussed in Section~\ref{sec:reduced-basis}.

\section{Adaptive ANOVA and Stochastic Collocation} \label{sec:anova}

In this section, we introduce the ANOVA decomposition~\cite{Cao-2009-ANOVA,Fisher-1925-SMR,Foo-2010-MPC} which is a useful tool for analyzing a multivariate function. In particular, we utilize a truncated ANOVA decomposition to approximate the computationally expensive high-dimensional stochastic problem with a set of cheaper low-dimensionsal problems. In addition, we adaptively select the most effective ANOVA terms, following the ideas presented in~\cite{Ma-2010-AHS,Yang-2012-AAD}.

\subsection{ANOVA Decomposition}

Let~$\genericfunction(x;\xi)$ denote a parametrized function for~$x \in \domain$ and~$\xi \in \Gamma$. In our application,~$\genericfunction$ may be either the velocity~$u$ or pressure~$p$. Let~$\mathcal{I}$ denote the index set~$\{1,\ldots,M\}$. We will refer to a subset~$T \subseteq \mathcal{I}$ as a \emph{direction}. ANOVA decomposes the function~$\genericfunction$ into contributions from individual directions~$T$. More formally, the decomposition has the form

\begin{equation}
  \genericfunction(x;\xi) = \sum_{T \subseteq \mathcal{I}} \genericfunction_T(x;\xi_T) \label{eq:anova},
\end{equation}

\noindent where~$\xi_T$ is a restriction of~$\xi$ to the coefficients in~$T$ and

\begin{equation}
  \genericfunction_T(x; \xi_T) = \int_{\Gamma_{T'}} \genericfunction(x;\xi) d\xi_{T'} - \sum_{S \subset T} \genericfunction_S(x; \xi_S), \label{eq:anova-term}
\end{equation}

\noindent with~$T' \subseteq \mathcal{I}$ denoting the complement of~$T$ in~$\mathcal{I}$ and~$\Gamma_S = \prod_{m \in S} \Gamma_m$ for~$S \subseteq \mathcal{I}$. We define the order of a direction~$T$ to be its cardinality~$|T|$. Thus, the ANOVA decomposition attempts to decompose~$\genericfunction$ into the contributions of each individual direction by removing the contributions from lower-order subdirections. The ANOVA decomposition of~$\genericfunction$ is built by first computing the order-zero term

\begin{equation*}
  \genericfunction_\emptyset (x) = \int_{\Gamma} \genericfunction(x; \xi) d\xi,
\end{equation*}

\noindent then successively computing higher-order terms by first marginalizing out the unused variables with~$\int_{\Gamma_{T'}} \genericfunction(x;\xi) d\xi_{T'}$ and then removing the computed effects of previous subdirections with~$\sum_{S \subset T} \genericfunction_S(x; \xi_S)$. While the cost of the full ANOVA decomposition is prohibitive in high-dimensional spaces, it may be used to obtain a useful approximation by truncating the expansion by keeping only the low-order terms. Thus, a high-dimensional multivariate function may be approximated through the sum of low-dimensional functions.

The ANOVA decomposition defined by~(\ref{eq:anova})--(\ref{eq:anova-term}) involves the computation of high-dimensional integrals which require expensive Monte Carlo computations. An alternative is anchored ANOVA~\cite{Ma-2010-AHS,Xu-2004-GDR,Yang-2012-AAD,Zhang-2012-EEA,Zhang-2011-APM} in which an anchor point~$c \in \Gamma$ is carefully chosen and the measure in~(\ref{eq:anova-term}) is replaced by the Dirac measure~$\delta(\xi-c)$. The order-zero term may then be computed by evaluating~$\genericfunction$ at the anchor point~$c$:

\begin{equation*}
  \genericfunction_\emptyset (x) = \genericfunction(x; c).
\end{equation*}

\noindent Given the order-zero term, we may compute the first-order terms and then the second-order terms as

\begin{align*}
  \genericfunction_{\{i\}} (x; \xi_{\{i\}}) &= \genericfunction(x; c, \xi, {\{i\}}) - \genericfunction_\emptyset(x) \\
  \genericfunction_{\{i,j\}} (x; \xi_{\{i,j\}}) &= \genericfunction(x; c, \xi, {\{i,j\}}) - \genericfunction_{\{i\}}(x; \xi_{\{i\}}) \\
  \, &\quad - \genericfunction_{\{j\}}(x; \xi_{\{j\}}) - \genericfunction_\emptyset(x).
\end{align*}

\noindent The notation~$\genericfunction(x;c,\xi, T)$ is to be understood as evaluating~$\genericfunction$ at the parameter~$\widehat{\xi}$ which takes on the values

\begin{equation*}
  \widehat{\xi}_i = \left\{
  \begin{array}{ll}
    c_i & i \notin T \\
    \xi_i & i \in T.
  \end{array} \right.
\end{equation*}

\noindent Third and higher-order terms may be computed in a similar fashion.

It was shown in~\cite{Ma-2010-AHS} that the mean of~$\xi$ often serves as a good choice of anchor point, and we shall use this choice.

\subsection{Stochastic Collocation}

We may utilize the ANOVA decomposition to compute moments of~$\genericfunction$. The mean of~$\genericfunction$ may simply be obtained by summing the means of the individual terms:

\begin{equation}
  \mathbb{E}[\genericfunction(x;\xi)] = \sum_{T \subseteq \mathcal{I}} \mathbb{E}[\genericfunction_T(x;\xi_T)] . \label{eq:anova-mean}
\end{equation}

\noindent A property of the standard ANOVA decomposition is the orthogonality of its terms. Due to this, higher-order moments such as the variance of~$\genericfunction$ may also be obtained by summing the moments of the individual terms. However, the orthogonality property only holds if the same measure is used in both the computation of the ANOVA terms and in the statistical moments. This fails in the case of anchored ANOVA. We remark that, for a good choice of an anchor point, the orthogonality may be approximately preserved. For general choices of anchor points, a method involving the covariances of individual terms is needed for good accuracy~\cite{Tang-2015-SAA}. That is, we compute the variance by summing the covariances of all pairs of directions:

\begin{equation*}
  \mathbb{E}[(\genericfunction(x;\xi) - \mu)^2] = \sum_{T,S \subseteq \mathcal{I}} \mathbb{E}[(\genericfunction_T(x;\xi_T) - \mu_T)(\genericfunction_S(x;\xi_S) - \mu_S)]
\end{equation*}

\noindent where~$\mu$ denotes the mean~(\ref{eq:anova-mean}) and~$\mu_T$ denotes the mean of~$\genericfunction_T$.

To estimate the moments of individual ANOVA terms, we use stochastic collocation. In stochastic collocation, we interpolate the stochastic solution by a set of Lagrange polynomials at collocation points

\begin{displaymath}
  \genericfunction(x;\xi) \approx \sum_{\xi^{(k)} \in \Theta} \genericfunction_c(x; \xi^{(k)}) L_{\xi^{(k)}}(\xi)
\end{displaymath}

\noindent where~$\Theta \subset \Gamma$ is the set of collocation points and~$\{L_{\xi^{(k)}}\}$ are the Lagrange polynomials. The coefficients~$\genericfunction_c(x; \xi^{(k)})$ are obtained by function evaluations at specific realizations~$\xi^{(k)}$ of~$\xi$. Thus, stochastic collocation methods reduce solution of the stochastic problem to a set of deterministic solves at specific sample points.

We select collocation points for each direction~$T$ as follows. For each index~$i \in T$, a set~$\Theta_i$ of~$p_i$ points and associated weights are selected as nodes according to a quadrature rule on~$\Gamma_i$. Then the full set of quadrature points is obtained by using a tensor product~$\Theta_T = \bigotimes_{i \in T} \Theta_i$. For a quadrature point~$\widehat{\xi} = \bigotimes_{i \in T} \widehat{\xi}_i \in \Theta_T$, the corresponding weight is the product~$w(\widehat{\xi}) = \prod_{i \in T} w_i(\widehat{\xi}_i)$ where~$w_i(\widehat{\xi}_i)$ denotes the weight of~$\widehat{\xi}_i$ in the quadrature rule for index~$i$. The collocation means for direction~$T$ and for the full function may then be computed as

\begin{align}
  \mathbb{E}_{\mathrm{sc}}[\genericfunction_T(x;\xi_T)] &= \sum_{\widehat{\xi} \in \Theta_T} \genericfunction_T(x;\widehat{\xi}) w(\widehat{\xi}) \nonumber \\
  \mathbb{E}_{\mathrm{sc}}[\genericfunction(x;\xi)] &= \sum_{T \subseteq \mathcal{I}} \mathbb{E}[\genericfunction_T(x;\xi_T)] . \label{eq:collocation-mean}
\end{align}

For a direction S, let~$\widetilde{\genericfunction}_S(x;\xi_S) = \genericfunction_S(x;\xi_S) - \mathbb{E}_{\mathrm{sc}}[\genericfunction_S(x;\xi_S)]$. To compute the covariance of~$S$ and~$T$, we need quadrature points for~$S \cup T$, although we need not form these quadrature points explicitly. We partition~$S \cup T$ into three disjoint sets:~$S \cap T$,~$S \setminus T$, and~$T \setminus S$. For~$\widehat{\xi}_{S \cap T} \in \Theta_{S \cap T}$, define

\begin{displaymath}
  Z_{S \setminus T}(\widehat{\xi}_{S \cap T}) = \sum_{\substack{\widehat{\xi}_S \in \Theta_S\\\widehat{\xi}_S\mid_{S \cap T} = \widehat{\xi}_{S \cap T}}} w_S(\widehat{\xi}_S) \widetilde{\genericfunction}_S(x;\widehat{\xi}_S).
\end{displaymath}

\noindent We may then compute the collocation covariances and variances as

\begin{align}
  \mathbb{E}_{\mathrm{sc}}[\widetilde{\genericfunction}_S(x;\xi_S) \widetilde{\genericfunction}_T(x;\xi_T)] &= \sum_{\widehat{\xi} \in \Theta_{S \cap T}} \frac{1}{w(\widehat{\xi})} Z_{S \setminus T}(\widehat{\xi}) Z_{T \setminus S}(\widehat{\xi}) \nonumber \\
  \mathbb{E}_{\mathrm{sc}}[(\genericfunction(x;\xi) - \mu)^2] &= \sum_{S,T \subseteq \mathcal{I}} \mathbb{E}_{\mathrm{sc}}[\widetilde{\genericfunction}_S(x;\xi_S) \widetilde{\genericfunction}_T(x;\xi_T)]. \label{eq:collocation-variance}
\end{align}

\noindent For truncated ANOVA, we replace~(\ref{eq:collocation-mean})--(\ref{eq:collocation-variance}) with sums over the used directions.

If a truncated ANOVA decomposition is used to compute~(\ref{eq:collocation-mean})--(\ref{eq:collocation-variance}), immediate savings over the full tensor product collocation are apparent. Assuming each index has polynomial order~$p$, the full tensor product collocation requires function evaluations at~$p^M$ collocation points, which quickly becomes prohibitive even for moderate~$M$. However, if an ANOVA decomposition truncated at level~$\ell$ is used, then the total number of collocation points is reduced to~$\sum_{l=0}^{\ell} \binom{M}{l} p^l$, which is substantially less than~$p^M$.

\subsection{Adaptive ANOVA}

Instead of simply truncating~(\ref{eq:anova}) at a specified level, we could attempt to adaptively select which directions contribute the most significantly. Such an adaptive scheme was developed in~\cite{Ma-2010-AHS,Yang-2012-AAD}. For each direction~$T$, we can score its contribution to the mean using

\begin{equation}
  \eta^{\mathrm{A}}(T) = \frac{\|\mathbb{E}_{\mathrm{sc}}[u_T]\|_{\velocityspace} + \|\mathbb{E}_{\mathrm{sc}}[p_T]\|_{\pressurespace}}{\|\sum_{|S| < |T|} \mathbb{E}_{\mathrm{sc}}[u_S]\|_{\velocityspace} + \|\sum_{|S| < |T|} \mathbb{E}_{\mathrm{sc}}[p_S]\|_{\pressurespace}}. \label{eq:anova-indicator}
\end{equation}

\noindent The form of this indicator was used in~\cite{Liao-2016-RBA} when applying adaptive ANOVA techniques to the Stokes equations due to the mixed formulation of the model problem. An alternative indicator which uses the variance rather than the mean may be used instead. In this paper, we use the mean.

We consider a direction~$T$ \emph{active} if it is included in the ANOVA decomposition. We consider any active direction~$T$ whose indicator~(\ref{eq:anova-indicator}) exceeds a given tolerance~$\epsilon^{\mathrm{A}}$ to be \emph{effective}. Let~$\mathcal{J}_l$ denote the set of active directions at level~$l$, and let~$\widetilde{\mathcal{J}}_l \subseteq \mathcal{J}_l$ denote the set of effective directions at level~$l$:

\begin{align}
  \widetilde{\mathcal{J}}_l &\equiv \{ T \in \mathcal{J}_l \colon \eta^{\mathrm{A}}(T) > \epsilon^{\mathrm{A}}\} . \label{eq:effective-directions}
\end{align}

\noindent Then we choose the next-level active directions~$\mathcal{J}_{l+1}$ as those level~$l+1$ directions~$T$ satisfying~$S \in \widetilde{\mathcal{J}}_l$ for every level~$l$ subdirection~$S$ of~$T$:

\begin{align}
  \mathcal{J}_{l+1} &\equiv \{ T \colon |T| = l+1 \text{ and } S \in \widetilde{\mathcal{J}}_{l} \, \forall S \subset T, |S| = l\} . \label{eq:next-level-directions}
\end{align}

\noindent Thus, we build level~$l+1$ active directions by considering only those directions which can be built from effective level~$l$ directions. This heuristic is based upon the idea that, if a direction~$T$ is important, then its subdirections at the previous level will likely also be important.

\section{Reduced Basis Methods} \label{sec:reduced-basis}

Whereas the anchored ANOVA method reduces the number of collocation points needed to accurately estimate the moments, solving~(\ref{eq:saddle-point-parametrized}) at each collocation point may still be prohibitive. In the context of reduced basis methods, we refer to solving~(\ref{eq:saddle-point-parametrized}) as a \emph{high-fidelity} solve. Reduced basis methods~\cite{Boyaval-2010-RBT,Haasdonk-2008-RBM,Hesthaven-2016-CRB,Quarteroni-2016-RBM} may be used in conjunction with ANOVA~\cite{Cho-2017-RBA,Liao-2016-RBA} to reduce the number of high-fidelity solves. Suppose that the solutions at the collocation points may be well-approximated by a low-dimensional space. Significant computational savings may then be obtained by reducing the full computation to this low-dimensional space. We seek spaces~$\rb{\velocityspace} \subset \fe{\velocityspace}$ and~$\rb{\pressurespace} \subset \fe{\pressurespace}$ of dimensions much smaller than those of~$\fe{\velocityspace}$ and~$\fe{\pressurespace}$ but for which the solution~$(\rb{u}(\xi),\rb{p}(\xi)) \in \rb{\velocityspace} \times \rb{\pressurespace}$ to the variational problem

\begin{align}
  a(\rb{u}(\xi),\rb{v};\xi) + b(\rb{v},\rb{p}(\xi)) &= \ell_1(\rb{v};\xi), \quad \forall \rb{v} \in \rb{\velocityspace}, \label{eq:rb-weak-momentum} \\
  b(\rb{u}(\xi),\rb{q}) &= \ell_2(\rb{q}), \quad \forall \rb{q} \in \rb{\pressurespace} \label{eq:rb-weak-mass}
\end{align}

\noindent accurately approximates the solution~$(\fe{u}(\xi),\fe{p}(\xi)) \in \fe{\velocityspace} \times \fe{\pressurespace}$ to the high-fidelity problem.

Let~$\rb{\dof{\velocityspace}}$ denote a matrix whose columns form a~$\velocityspace$-orthogonal basis for the space~$\rb{\velocityspace}$, and, similarly, let~$\rb{\dof{\pressurespace}}$ denote a matrix whose columns form a~$\pressurespace$-orthogonal basis for the space~$\rb{\pressurespace}$. Then~(\ref{eq:rb-weak-momentum})--(\ref{eq:rb-weak-mass}) is equivalent to solving the following saddle-point problem as in~(\ref{eq:saddle-point-parametrized})

\begin{equation}
  \left[
    \begin{array}{cc}
      \rb{\dof{\velocityspace}}^T \dof{A}(\xi) \rb{\dof{\velocityspace}} & \rb{\dof{\velocityspace}}^T \dof{B}^T \rb{\dof{\pressurespace}} \\
      \rb{\dof{\pressurespace}}^T \dof{B} \rb{\dof{\velocityspace}} & \dof{0}
    \end{array}
    \right]
  \left[
    \begin{array}{c}
      \rb{\dof{u}}(\xi) \\
      \rb{\dof{p}}(\xi)
    \end{array}
    \right]
  =
  \left[
    \begin{array}{c}
      \rb{\dof{\velocityspace}}^T \dof{f}(\xi) \\
      \rb{\dof{\pressurespace}}^T \dof{g}
    \end{array}
    \right] . \label{eq:reduced-basis-saddle-point}
\end{equation}

\noindent The approximate solutions in~$\fe{\velocityspace}$ and~$\fe{\pressurespace}$ may then be obtained as~$\rb{\dof{\velocityspace}}\rb{\dof{u}}(\xi)$ and~$\rb{\dof{\pressurespace}}\rb{\dof{p}}(\xi)$.

Care must be taken to ensure that the reduced basis problem~(\ref{eq:rb-weak-momentum})--(\ref{eq:rb-weak-mass}) is inf-sup stable. Whereas the coercivity and continuity of~$a$ on the reduced basis space follow directly from the coercivity and continuity of~$a$ on the high-fidelity space, the inf-sup condition~(\ref{eq:inf-sup-condition}) is not immediately satisfied. The condition may be satisfied by properly enriching the velocity space~$\rb{\velocityspace}$ relative to the pressure space. Let~$T \colon \pressurespace \rightarrow \velocityspace$ denote the \emph{supremizer} operator such that, for a given~$q \in \pressurespace$,~$Tq$ satisfies

\begin{align}
  (Tq,v)_\velocityspace &= b(v,q) \quad \forall v \in \velocityspace. \label{eq:supremizer}
\end{align}

\noindent Then the inf-sup condition is satisfied if, for each~$q \in \rb{\pressurespace}$,~$Tq \in \rb{\velocityspace}$~\cite{Rozza-2007-SRB}. Let~$\{\rb{u}^1,\ldots,\rb{u}^s\}$ be a set of~$s$ high-fidelity velocity solutions from which we seek to build the reduced basis and~$\{\rb{p}^1,\ldots,\rb{p}^s\}$ be the corresponding set of~$s$ high-fidelity pressure solutions. We then obtain an inf-sup stable approximation by defining

\begin{align}
  \rb{\velocityspace} &\equiv \mathrm{span}\left\{\rb{u}^1,\ldots,\rb{u}^s,T\rb{p}^1,\ldots,T\rb{p}^s\right\} \nonumber \\
  \rb{\pressurespace} &\equiv \mathrm{span}\left\{\rb{p}^1,\ldots,\rb{p}^s\right\} . \nonumber
\end{align}

\noindent Note that, with this procedure, the dimension of~$\rb{\velocityspace}$ is always twice that of~$\rb{\pressurespace}$.

The efficiency of the reduced basis method stems from the dimension of~(\ref{eq:reduced-basis-saddle-point}) being significantly smaller than the dimension of~(\ref{eq:saddle-point-parametrized}). The accuracy relies upon the construction of rigorous a posteriori error estimates which bound the errors of the reduced basis solutions in the~$\velocityspace$ and~$\pressurespace$ norms. These will be the subject of Section~\ref{sec:aee}. For now, we simply assume such error estimates are available.

We desire the entire computation on the reduced space to be independent of the size of the high-fidelity spaces. For solving~(\ref{eq:reduced-basis-saddle-point}), this means that the assembly of the saddle-point system must be done independently of the size of the high-fidelity problem. For this, we exploit the affine decomposition~(\ref{eq:affine-decomposition-A})--(\ref{eq:affine-decomposition-f}) to obtain the reduced basis affine decompositions

\begin{align}
  \rb{\dof{\velocityspace}}^T \dof{A}(\xi) \rb{\dof{\velocityspace}} &= \sum_{i=1}^{n_A} \theta_i^A(\xi) \rb{\dof{\velocityspace}}^T \dof{A}_i \rb{\dof{\velocityspace}} \nonumber \\
  \rb{\dof{\velocityspace}}^T \dof{f}(\xi) &= \sum_{i=1}^{n_f} \theta_i^f(\xi) \rb{\dof{\velocityspace}}^T \dof{f}_i . \nonumber
\end{align}

\noindent Thus, after constructing the reduced basis, we compute and store the parameter-independent quantities $\rb{\dof{\velocityspace}}^T \dof{A}_i \rb{\dof{\velocityspace}}$, $\rb{\dof{\velocityspace}}^T \dof{f}_i$, $\rb{\dof{\pressurespace}}^T \dof{B} \rb{\dof{\velocityspace}}$, and $\rb{\dof{\pressurespace}}^T \dof{g}$.

\subsection{Constructing the Reduced Basis}

We form the reduced basis using a variation of the standard greedy algorithm~\cite[Chapter 7]{Quarteroni-2016-RBM} which is summarized in Algorithm~\ref{alg:greedy-alg}. In the standard greedy algorithm, the reduced basis is formed in a potentially expensive offline step. A finite training sample set~$\Xi \subset \Gamma$ is selected from the space of parameters such that~$\Xi$ is a good representation of the parameter space. To fit within our full reduced basis ANOVA algorithm (Algorithm~\ref{alg:rb-anova}), we allow the case where the greedy training algorithm can extend an existing reduced basis. If no reduced basis currently exists, a sample point~$\xi \in \Xi$ is selected at which a high-fidelity solve is performed to initialize the reduced basis. For each iteration of training, reduced basis solves are performed over all points in~$\Xi$ and the a posteriori error estimates are computed for each reduced basis solution. If all error estimates are below a prescribed tolerance~$\epsilon^{\mathrm{RB}}$, then training stops and the current reduced basis is output. Otherwise, the parameter which attains the largest error estimate is selected and a high-fidelity solve is performed at this parameter. The high-fidelity solution is then added to the reduced basis. This process repeats until all error estimates are below the desired tolerance. To ensure inf-sup stability, in addition to adding the high-fidelity velocity solution to the velocity reduced basis, we must also add the application of the supremizer~(\ref{eq:supremizer}) to the corresponding pressure solution. A Gram-Schmidt orthonormalization procedure should be used to ensure that the reduced basis matrices remain orthonormal with respect to the appropriate inner product.

\begin{algorithm}[h!]
  \caption{Greedy algorithm for generating a reduced basis using a training set~$\Xi$ and a tolerance~$\epsilon^{\mathrm{RB}}$.}
  \begin{algorithmic}[1]
    \IF{$\rb{\dof{\velocityspace}}$ and $\rb{\dof{\pressurespace}}$ are uninitialized}
    \STATE choose~$\widehat{\xi} \in \Xi$
    \STATE compute high-fidelity solution~(\ref{eq:saddle-point-parametrized})~$(\dof{u}(\widehat{\xi}),\dof{p}(\widehat{\xi}))$
    \STATE compute supremizer~(\ref{eq:supremizer}) $\dof{v}(\widehat{\xi}) = \dof{M}_\velocityspace^{-1} \dof{B}^T \dof{p}(\widehat{\xi})$
    \STATE $\rb{\dof{\velocityspace}} \leftarrow \left[\dof{u}(\widehat{\xi}) \mid \dof{v}(\widehat{\xi}) \right]$
    \STATE $\rb{\dof{\pressurespace}} \leftarrow \left[\dof{p}(\widehat{\xi}) \right]$
    \ENDIF
    \LOOP
    \FOR{$\xi \in \Xi$}
    \STATE compute RB solution~(\ref{eq:reduced-basis-saddle-point})~$(\rb{\dof{u}}(\xi),\rb{\dof{p}}(\xi))$
    \STATE compute a posteriori error estimate~(\ref{eq:aee})~$\rb{\Delta}(\xi)$
    \ENDFOR
    \STATE let~$\Delta = \max_{\xi \in \Xi} \rb{\Delta}(\xi)$
    \IF{$\Delta < \epsilon^{\mathrm{RB}}$}
    \RETURN reduced basis and offline quantities
    \ELSE
    \STATE let~$\widehat{\xi} = \mathop{\mathrm{argmax}}_{\xi \in \Xi} \rb{\Delta}(\xi)$
    \STATE compute high-fidelity solutions~(\ref{eq:saddle-point-parametrized})~$(\dof{u}(\widehat{\xi}),\dof{p}(\widehat{\xi}))$
    \STATE compute supremizer~(\ref{eq:supremizer}) $\dof{v}(\widehat{\xi}) = \dof{M}_\velocityspace^{-1} \dof{B}^T \dof{p}(\widehat{\xi})$
    \STATE $\rb{\dof{\velocityspace}} \leftarrow \left[\rb{\dof{\velocityspace}} \mid \dof{u}(\widehat{\xi}) \mid \dof{v}(\widehat{\xi}) \right]$
    \STATE $\rb{\dof{\pressurespace}} \leftarrow \left[\rb{\dof{\pressurespace}} \mid \dof{p}(\widehat{\xi}) \right]$
    \ENDIF
    \ENDLOOP
  \end{algorithmic}
  \label{alg:greedy-alg}
\end{algorithm}

Our variation of the standard greedy algorithm follows that of~\cite{Elman-2013-RBC} for reduced basis collocation in which the formation of the collocation points and the generation of the reduced basis are not separated. That is, as we are building the collocation points at each ANOVA level, we use those collocation points as our training set~$\Xi$. In~\cite{Elman-2013-RBC}, the authors proposed a single sweep through these points by augmenting the reduced basis whenever a single parameter's error estimate exceeds the tolerance. In this paper, however, we use the standard approach of computing error estimates over all training points and augmenting the reduced basis with a high-fidelity solve at the point with the largest error estimate. We remark that this approach is easily parallelized as the reduced basis solves and error estimate computations are embarrassingly parallel. The approach in~\cite{Elman-2013-RBC}, however, is better suited for serial computation. The full reduced basis ANOVA algorithm is summarized in Algorithm~\ref{alg:rb-anova}.

\begin{algorithm}[h!]
  \caption{Reduced Basis Adaptive ANOVA with initial level~$\ell$, maximum level~$L$, reduced basis tolerance~$\epsilon^{\mathrm{RB}}$, and ANOVA tolerance~$\epsilon^{\mathrm{A}}$.}
  \begin{algorithmic}[1]
    \STATE compute anchor point~$c$ as mean of distribution
    \STATE solve~$\dof{u}(c), \dof{p}(c)$ of high-fidelity problem~(\ref{eq:saddle-point-parametrized}) at anchor point~$c$
    \STATE compute supremizer~(\ref{eq:supremizer}) $\dof{v}(c) = \dof{M}_\velocityspace^{-1} \dof{B}^T \dof{p}(c)$
    \STATE $\rb{\dof{\velocityspace}} \leftarrow \left[\dof{u}(c) \mid \dof{v}(c) \right]$
    \STATE $\rb{\dof{\pressurespace}} \leftarrow \left[\dof{p}(c) \right]$
    \STATE initialize set of active directions~$\mathcal{J}_{\ell}$ as all directions of order up to~$\ell$
    \STATE initialize~$\Xi_0 \leftarrow \{c\}$
    \FOR{$l = \ell$ \TO $L$}
    \STATE generate all collocation points~$\Xi_l$ for the active directions~$\mathcal{J}_{\ell}$
    \STATE form reduced basis training set~$\Xi \leftarrow \Xi_l \setminus \Xi_0$
    \STATE update reduced basis matrices~$\rb{\dof{\velocityspace}}$ and~$\rb{\dof{\pressurespace}}$ using training set~$\Xi$ and tolerance~$\epsilon^{\mathrm{RB}}$ (see Algorithm~\ref{alg:greedy-alg})
    \STATE $\Xi_0 \leftarrow \Xi_l \cup \Xi_0$
    \STATE perform reduced basis solves on~$\Xi_0$ and compute ANOVA indicators~(\ref{eq:anova-indicator})
    \IF{$l < L$}
    \STATE compute effective directions~$\widetilde{\mathcal{J}}_l$~(\ref{eq:effective-directions})
    \STATE compute next-level active directions~$\mathcal{J}_{l+1}$~(\ref{eq:next-level-directions})
    \IF{$|\mathcal{J}_{l+1}| = 0$}
    \RETURN reduced basis and ANOVA solutions
    \ENDIF
    \ENDIF
    \ENDFOR
    \RETURN reduced basis and ANOVA solutions
  \end{algorithmic}
  \label{alg:rb-anova}
\end{algorithm}

\section{A Posteriori Error Estimates} \label{sec:aee}

Essential for the successful application of reduced basis methods is the computation of rigorous a posteriori error estimates. We follow the approach of~\cite{Gerner-2012-CRB} for saddle-point problems, which develops residual-based error estimates based upon the Brezzi stability theory. Let~$\rb{u}(\xi) \in \rb{\velocityspace}$ and~$\rb{p}(\xi) \in \rb{\pressurespace}$ denote the reduced basis velocity and pressure solutions, respectively, for parameter~$\xi$, and let~$\fe{u}(\xi) \in \fe{\velocityspace}$ and~$\fe{p}(\xi) \in \fe{\pressurespace}$ denote the high-fidelity velocity and pressure solutions, respectively. We are interested in bounding the velocity and pressure errors

\begin{align*}
  \rb{e}^u(\xi) &\equiv \fe{u}(\xi) - \rb{u}(\xi) \in \fe{\velocityspace} \\
  \rb{e}^p(\xi) &\equiv \fe{p}(\xi) - \rb{p}(\xi) \in \fe{\pressurespace}.
\end{align*}

\noindent We note that these errors are the solution to the Stokes-Brinkman problem

\begin{align}
  a(\rb{e}^u(\xi),v;\xi) + b(v,\rb{e}^p(\xi)) &= r_1(v;\xi), \quad \forall v \in \fe{\velocityspace} \label{eq:residual-equation1} \\
  b(\rb{e}^u(\xi),q) &= r_2(q;\xi), \quad \forall q \in \fe{\pressurespace}  \label{eq:residual-equation2}
\end{align}

\noindent where~$r_1(\cdot;\xi) \in \dual{\fe{\velocityspace}}$ and~$r_2(\cdot;\xi) \in \dual{\fe{\pressurespace}}$ are the reduced basis residuals defined by

\begin{align}
  r_1(v;\xi) &\equiv \ell_1(v;\xi) - a(\rb{u}(\xi),v;\xi) - b(v,\rb{p}(\xi)) \label{eq:rb-velocity-residual} \\
  r_2(q;\xi) &\equiv \ell_2(q) - b(\rb{u}(\xi),q). \label{eq:rb-pressure-residual}
\end{align}

We are interested in tight upper bounds~$\rb{\Delta}^u(\xi)$ and~$\rb{\Delta}^p(\xi)$ which satisfy

\begin{align}
  \|\rb{e}^u(\xi)\|_\velocityspace &\leq \rb{\Delta}^u(\xi) \nonumber \\
  \|\rb{e}^p(\xi)\|_\pressurespace &\leq \rb{\Delta}^p(\xi). \nonumber
\end{align}

\noindent We may then form an a posteriori error estimate over the whole solution as

\begin{align}
  \rb{\Delta}(\xi) &= \sqrt{\rb{\Delta}^u(\xi)^2 + \rb{\Delta}^p(\xi)^2}. \label{eq:aee}
\end{align}

\noindent Applying the stability bounds~(\ref{eq:stability-bounds-u})--(\ref{eq:stability-bounds-p}) from Theorem~\ref{thm:stability} to~(\ref{eq:residual-equation1})--(\ref{eq:residual-equation2}), we obtain the following upper bounds for the errors:

\begin{align}
  \|\rb{e}^u(\xi)\|_\velocityspace &\leq \frac{1}{\alpha(\xi)}\|r_1(\cdot;\xi)\|_{\dual{\velocityspace}} + \frac{2}{\beta}\sqrt{\frac{\gamma(\xi)}{\alpha(\xi)}}\|r_2(\cdot;\xi)\|_{\dual{\pressurespace}} \label{eq:stability-estimate1} \\
  \|\rb{e}^p(\xi)\|_\pressurespace &\leq \frac{2}{\beta}\sqrt{\frac{\gamma(\xi)}{\alpha(\xi)}}\|r_1(\cdot;\xi)\|_{\dual{\velocityspace}} + \frac{\gamma(\xi)}{\beta^2} \|r_2(\cdot;\xi)\|_{\dual{\pressurespace}}, \label{eq:stability-estimate2}
\end{align}

\noindent where~$\alpha(\xi)$ is the coercivity constant

\begin{equation*}
  \alpha(\xi) \equiv \inf_{v \in \fe{\velocityspace}} \frac{a(v,v;\xi)}{\|v\|^2_\velocityspace},
\end{equation*}

\noindent $\gamma(\xi)$ is the continuity constant

\begin{equation*}
  \gamma(\xi) \equiv \sup_{u,v \in \fe{\velocityspace}} \frac{a(u,v;\xi)}{\|u\|_\velocityspace \|v\|_\velocityspace},
\end{equation*}

\noindent and~$\beta$ is the inf-sup constant

\begin{equation*}
  \beta \equiv \inf_{q \in \fe{\pressurespace}} \sup_{v \in \fe{\velocityspace}} \frac{b(v,q)}{\|q\|_\pressurespace \|v\|_\velocityspace}.
\end{equation*}

Note that the inf-sup constant~$\beta$ is defined independently of the parameter~$\xi$. Thus, we only need to solve a single eigenvalue problem to determine its value over all parameters. The coercivity~$\alpha(\xi)$ and continuity~$\gamma(\xi)$ constants are parameter-dependent so that computing their exact values would require solving separate eigenvalue problems for each parameter. To reduce the computational cost, we instead compute tight bounds in an efficient manner that is independent of the dimension of~$\fe{\velocityspace}$. Details of this procedure are described in Section~\ref{sec:scm}.

Using a lower bound~$\alpha^{\mathrm{LB}}(\xi) \leq \alpha(\xi)$ for the coercivity constant and an upper bound~$\gamma^{\mathrm{UB}}(\xi) \geq \gamma(\xi)$, for the continuity constant, we may define our error estimates as

\begin{align}
  \rb{\Delta}^u(\xi) &\equiv \frac{1}{\alpha^{\mathrm{LB}}(\xi)}\|r_1(\cdot;\xi)\|_{\dual{\velocityspace}} + \frac{2}{\beta}\sqrt{\frac{\gamma^{\mathrm{UB}}(\xi)}{\alpha^{\mathrm{LB}}(\xi)}}\|r_2(\cdot;\xi)\|_{\dual{\pressurespace}} \label{eq:rb-velocity-bound} \\
  \rb{\Delta}^p(\xi) &\equiv \frac{2}{\beta}\sqrt{\frac{\gamma^{\mathrm{UB}}(\xi)}{\alpha^{\mathrm{LB}}(\xi)}}\|r_1(\cdot;\xi)\|_{\dual{\velocityspace}} + \frac{\gamma^{\mathrm{UB}}(\xi)}{\beta^2} \|r_2(\cdot;\xi)\|_{\dual{\pressurespace}}. \label{eq:rb-pressure-bound}
\end{align}

\noindent Since these are upper bounds for the stability estimates~(\ref{eq:stability-estimate1})--(\ref{eq:stability-estimate2}), we observe that these are true upper bounds for the reduced basis errors.

\subsection{Efficient Computation of the Dual Norm}

Computing the error estimates~(\ref{eq:rb-velocity-bound})--(\ref{eq:rb-pressure-bound}) requires computation of the dual norms of the residuals. Dual norm computations are best performed through use of the Riesz representation. Let~$R_H \colon H \rightarrow \dual{H}$ denote the Riesz isomorphism from a Hilbert space~$H$ to its dual defined by~$\dualitybracket{R_H u}{v}{H} = (u,v)_H$ for all~$u,v \in H$. Here,~$\dualitybracket{\cdot}{\cdot}{H}$ denotes the duality pairing and~$(\cdot,\cdot)_H$ denotes the inner product. Then, by the Riesz representation theorem, for any~$f \in \dual{H}$, we have~$\|f\|_{\dual{H}} = \|R_H^{-1}f\|_H$. For finite-dimensional spaces, it is apparent that the operator~$R_H$ is represented by the mass matrix for the~$H$-norm. Thus, if~$\dof{f}$ denotes the vector of coefficients for~$f \in \dual{H}$, then~$\|f\|_{\dual{H}}^2 = \dof{f}^T \dof{M}_H^{-1} \dof{f}$.

In order for the error estimates to be computed efficiently during the online phase, we desire the cost to be independent of the dimension of the finite element spaces~$\fe{\velocityspace}$ and~$\fe{\pressurespace}$. To this effect, we utilize the affine decompositions~(\ref{eq:affine-decomposition-A})--(\ref{eq:affine-decomposition-f}). The dual norm of the residual~(\ref{eq:rb-velocity-residual}) is then computed as

\begin{align}
  \|r_1(\cdot;\xi)\|_{\dual{\velocityspace}}^2 &= \sum_{i=1}^{n_f} \sum_{j=1}^{n_f} \theta_i^f(\xi) \theta_j^f(\xi) \dof{C}_{ij} \nonumber \\
  \, &\quad - 2 \left( \sum_{i=1}^{n_f} \sum_{j=1}^{n_A} \theta_i^f(\xi) \theta_j^A(\xi) \dof{D}_{ij} \right) \rb{\dof{u}} \nonumber \\
  \, &\quad - 2 \left( \sum_{i=1}^{n_f} \theta_i^f(\xi) \dof{E}_i \right) \rb{\dof{p}} \nonumber \\
  \, &\quad + \rb{\dof{u}}^T \left( \sum_{i=1}^{n_A} \sum_{j=1}^{n_A} \theta_i^A(\xi) \theta_j^A(\xi) \dof{F}_{ij} \right) \rb{\dof{u}} \nonumber \\
  \, &\quad + 2 \rb{\dof{u}}^T \left( \sum_{i=1}^{n_A} \theta_i^A(\xi) \dof{G_i} \right) \rb{\dof{p}} \nonumber \\
  \, &\quad + \rb{\dof{p}}^T \dof{H} \rb{\dof{p}}, \label{eq:offline-dual1}
\end{align}

\noindent where

\begin{align*}
  \dof{C}_{ij} &= \dof{f}_i^T \dof{M}_\velocityspace^{-1} \dof{f}_j \\
  \dof{D}_{ij} &= \dof{f}_i^T \dof{M}_\velocityspace^{-1} \dof{A}_j \rb{\dof{\velocityspace}} \\
  \dof{E}_i &= \dof{f}_i^T \dof{M}_\velocityspace^{-1} \dof{B}^T \rb{\dof{\pressurespace}} \\
  \dof{F}_{ij} &= \rb{\dof{\velocityspace}}^T \dof{A}_i^T \dof{M}_\velocityspace^{-1} \dof{A}_j \rb{\dof{\velocityspace}} \\
  \dof{G}_i &= \rb{\dof{\velocityspace}}^T \dof{A}_i \dof{M}_\velocityspace^{-1} \dof{B}^T \rb{\dof{\pressurespace}} \\
  \dof{H} &= \rb{\dof{\pressurespace}}^T \dof{B} \dof{M}_\velocityspace^{-1} \dof{B}^T \rb{\dof{\pressurespace}},
\end{align*}

\noindent and the dual norm of the residual~(\ref{eq:rb-pressure-residual}) is computed as

\begin{align}
  \|r_2(\cdot;\xi)\|_{\dual{\pressurespace}}^2 &= \dof{R} - 2 \dof{S} \rb{\dof{u}} + \rb{\dof{u}}^T \dof{T} \rb{\dof{u}} \nonumber
\end{align}

\noindent where

\begin{align*}
  \dof{R} &= \dof{g}^T \dof{M}_\pressurespace^{-1} \dof{g} \\
  \dof{S} &= \dof{g}^T \dof{M}_\pressurespace^{-1} \dof{B} \rb{\dof{\velocityspace}} \\
  \dof{T} &= \rb{\dof{\velocityspace}}^T \dof{B}^T \dof{M}_\pressurespace^{-1} \dof{B} \rb{\dof{\velocityspace}}.
\end{align*}

\noindent The terms~$\dof{C}_{ij}, \dof{D}_{ij}, \dof{E}_i, \dof{F}_{ij}, \dof{G}_i, \dof{H}, \dof{R}, \dof{S}, \dof{T}$ are all parameter-independent and can be computed and stored during the offline phase. Thus, the dual norms for the residuals may be computed efficiently during the online phase in a cost that is independent of the dimensions of~$\fe{\velocityspace}$ and~$\fe{\pressurespace}$.

\begin{remark}
  For moderately large reduced bases, the offline quantities may consume a significant amount of memory. In particular, there are~$n_A^2$ matrices of the form~$\dof{F}_{ij}$. If the size of the velocity reduced basis is~$N$, then each of these matrices will contain~$N^2$ entries, for a storage requirement of~$n_A^2 N^2$ double-precision floating-point numbers. Nonetheless, the computation is still independent of the size of the high-fidelity problem. This, however, should be considered when evaluating the efficiency of the reduced basis method.
\end{remark}

\subsection{Efficient Computation of Stability Bounds Using SCM} \label{sec:scm}

As previously remarked, we seek efficiently computable bounds on the coercivity and continuity constants with cost that is independent of the high-fidelity problem. A popular method for obtaining such bounds is the Successive Constraint Method (SCM), first proposed in~\cite{Huynh-2007-SCM} and subsequently refined in~\cite{Chen-2009-SCM}. We recall here the method as presented in~\cite{Chen-2009-SCM}.

For simplicity, we focus on the coercivity constant, since the method may easily be extended to the continuity case. Recall that the coercivity constant is the largest~$\alpha(\xi) > 0$ such that~$a(u,u;\xi) \geq \alpha(\xi) \|u\|_\velocityspace$ for all~$u \in \velocityspace$. Using the affine decomposition of~$a$, we may write

\begin{align}
  \alpha(\xi) &= \inf_{u \in \velocityspace} \sum_{i=1}^{n_A} \theta_i^A(\xi) \frac{a_i(u,u)}{\|u\|_\velocityspace} . \label{eq:scm-1}
\end{align}

\noindent Defining the set

\begin{align*}
  \mathcal{Y} &\equiv \{y = (y_1,\ldots,y_{n_A})^T \in \R^{n_A} \colon \nonumber \\
  \, &\quad \exists u \in \velocityspace \ni y_i = \frac{a_i(u,u)}{\|u\|_\velocityspace}, i = 1, \ldots, n_A \},
\end{align*}

\noindent we may express~(\ref{eq:scm-1}) as the minimization problem

\begin{align*}
  \alpha(\xi) &= \inf_{y \in \mathcal{Y}} \sum_{i=1}^{n_A} \theta_i^A(\xi) y_i .
\end{align*}

\noindent We may obtain a lower bound for~$\alpha(\xi)$ by replacing~$\mathcal{Y}$ with a set~$\mathcal{Y}_{\mathrm{LB}}(\xi)$ such that~$\mathcal{Y} \subset \mathcal{Y}_{\mathrm{LB}}(\xi)$. By choosing~$\mathcal{Y}_{\mathrm{LB}}(\xi)$ to be a convex polyhedron, we obtain the following linear programming problem

\begin{align*}
  \alpha^{\mathrm{LB}}(\xi) &= \min_{y \in \mathcal{Y}_{\mathrm{LB}}(\xi)} \sum_{i=1}^{n_A} \theta_i^A(\xi) y_i .
\end{align*}

To define~$\mathcal{Y}_{\mathrm{LB}}(\xi)$, we first note that the box~$\mathcal{B}$ defined as

\begin{align*}
  \mathcal{B} &\equiv \{ y \in \R^{n_A} \colon \nonumber \\
  \, &\quad \inf_{u \in \velocityspace} \frac{a_i(u,u)}{\|u\|_\velocityspace} \leq y_i \leq \sup_{u \in \velocityspace} \frac{a_i(u,u)}{\|u\|_\velocityspace}, i=1,\ldots,n_A \}
\end{align*}

\noindent clearly satisfies~$\mathcal{Y} \subset \mathcal{B}$. Forming~$\mathcal{B}$ involves the solution of~$2 n_A$ eigenvalue problems. Suppose we are given two sets of parameters--$\Xi_E$ for which exact coercivity constants have been computed and~$\Xi_P$ for which lower bounds have been computed. Then a finer~$\mathcal{Y}_{\mathrm{LB}}(\xi)$ may be obtained as follows. First, choose~$M_E, M_P \geq 0$ and locate the~$M_E$ points~$\{\xi_E^1,\ldots,\xi_E^{M_E}\}$ in~$\Xi_E$ that are closest to~$\xi$ and the~$M_P$ points~$\{\xi_P^1,\ldots,\xi_P^{M_P}\}$ in~$\Xi_P$ that are closest to~$\xi$. We then include the following~$M_E + M_P$ inequality constraints in addition to those defined by the box~$\mathcal{B}$:

\begin{align}
  \sum_{i=1}^{n_A} \theta_i^A(\xi_E^j) y_i & \geq \alpha(\xi_E^j), j = 1, \ldots, M_E \nonumber \\
  \sum_{i=1}^{n_A} \theta_i^A(\xi_P^j) y_i & \geq \alpha^{\mathrm{LB}}(\xi_P^j), j = 1, \ldots, M_P . \nonumber
\end{align}

\noindent It is clear that the set~$\mathcal{Y}_{\mathrm{LB}}(\xi)$ is a superset of~$\mathcal{Y}$. We thus reduce the problem of computing~$\alpha^{\mathrm{LB}}(\xi)$ to solving a linear program in~$n_A$ variables with~$2n_A + M_E + M_P$ linear inequality constraints, which is independent of the dimension of~$\fe{\velocityspace}$, as desired.

What remains is to demonstrate how~$\Xi_E$ is selected and the lower bounds for the points in~$\Xi_P$ are obtained. These may be built through an offline training procedure. This procedure will make use of the following upper bound for~$\alpha(\xi)$. Given any~$\xi'$ for which~$\alpha(\xi')$ has been computed, we may use the corresponding eigenvector~$u'$ and compute~$y' = (y_1,\ldots,y_{n_A})^T \in \mathcal{Y}$ where~$y_i = a_i(u',u')/\|u'\|_\velocityspace$. Collecting all such~$y'$ into a set~$\mathcal{Y}_{\mathrm{UB}}$, we obtain an upper bound as

\begin{align}
  \alpha^{\mathrm{UB}}(\xi) &= \min_{y \in \mathcal{Y}_{\mathrm{UB}}} \sum_{i=1}^{n_A} \theta_i^A(\xi) y_i . \nonumber
\end{align}

\noindent Since~$\mathcal{Y}_{\mathrm{UB}}$ is small and finite, we may compute this upper bound efficiently by simply enumerating over all elements of~$\mathcal{Y}_{\mathrm{UB}}$. Furthermore, since clearly~$\mathcal{Y}_{\mathrm{UB}} \subset \mathcal{Y}$, we have a true upper bound. Given upper and lower bounds for the coercivity constant, we define the indicator

\begin{align}
  \eta^{\mathrm{SCM}}(\xi) &\equiv 1 - \frac{\alpha^{\mathrm{LB}}(\xi)}{\alpha^{\mathrm{UB}}(\xi)} \label{eq:scm-indicator}
\end{align}

\noindent which measures the relative gap between the upper and lower bounds.

Having determined an indicator for the quality of the SCM approximation, we may now describe the offline training algorithm. We begin with a training set~$\Xi$ which should be representative of the parameter set. Choose~$\xi' \in \Xi$, initialize~$\Xi_E = \{\xi'\}$ and~$\Xi_P = \Xi \setminus \{\xi'\}$, and set~$\alpha^{\mathrm{LB}}(\xi) = 0$ for all~$\xi \in \Xi$. Then loop as follows. Solve for the exact coercivity constant for~$\xi'$ and compute lower bounds using SCM and indicators for all~$\xi \in \Xi_P$. If the largest indicator is below a desired tolerance~$\epsilon^{\mathrm{SCM}}$, terminate. Otherwise, select a new point~$\xi' \in \Xi_P$ with the largest indicator, compute and store~$\alpha(\xi')$, move~$\xi'$ from~$\Xi_P$ to~$\Xi_E$, and update the lower bounds for points in~$\Xi_P$ using SCM. Repeat until termination. This procedure is summarized in Algorithm~\ref{alg:scm-training}.

\begin{algorithm}
  \caption{Greedy training algorithm for SCM given training set~$\Xi$.}
  \begin{algorithmic}[1]
    \STATE $\Xi_E \leftarrow \emptyset$,~$\Xi_P \leftarrow \Xi$
    \STATE choose~$\xi' \in \Xi$
    \STATE $\alpha^{\mathrm{LB}}(\xi) \leftarrow 0$ for all~$\xi \in \Xi$
    \LOOP
    \STATE compute and store~$\alpha(\xi')$
    \STATE $\Xi_E \leftarrow \Xi_E \cup \{\xi'\}$,~$\Xi_P \leftarrow \Xi_P \setminus \{\xi'\}$
    \IF{$\Xi_P = \emptyset$}
    \RETURN SCM data
    \ENDIF
    \FOR{$\xi \in \Xi_P$}
    \STATE update~$\alpha^{\mathrm{LB}}(\xi)$ using SCM
    \STATE compute indicator~$\eta^{\mathrm{SCM}}(\xi)$~(\ref{eq:scm-indicator})
    \ENDFOR
    \STATE $\eta \leftarrow \max_{\xi \in \Xi_P} \eta^{\mathrm{SCM}}(\xi)$
    \IF{$\eta < \epsilon^{\mathrm{SCM}}$}
    \RETURN SCM data
    \ENDIF
    \STATE $\xi' \leftarrow \mathop{\mathrm{argmax}}_{\xi \in \Xi_P} \eta^{\mathrm{SCM}}(\xi)$
    \ENDLOOP
  \end{algorithmic}
  \label{alg:scm-training}
\end{algorithm}

\section{Numerical Experiments} \label{sec:numerical}

We compare the effectiveness of reduced basis ANOVA on a set of three model problems. These model problems are introduced in Section~\ref{sec:model-problems} and consider both isotropic and anisotropic permeabilities. In Section~\ref{sec:numerical-scm}, we analyze the performance of the SCM method on the three problems and provide details on the necessary eigenvalue computations. Finally, in Section~\ref{sec:rb-numerical}, we study the performance of reduced basis ANOVA.

\subsection{Model Problems} \label{sec:model-problems}

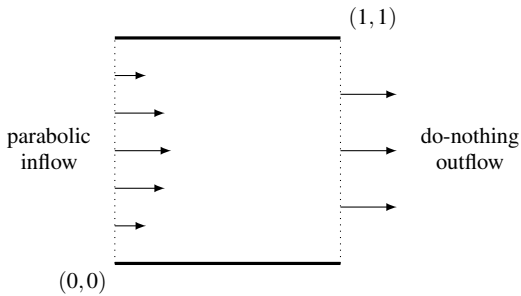
\begin{figure}[h!]
  \centering
  \begin{tikzpicture}[scale=3,>=latex]
    \draw[very thick] (0,0) -- (1,0);
    \draw[very thick] (0,1) -- (1,1);
    \draw[dotted] (0,0) -- (0,1);
    \draw[dotted] (1,0) -- (1,1);

    \draw[->] (0,0.166667) -- (0.138889,0.166667);
    \draw[->] (0,0.333333) -- (0.222222,0.333333);
    \draw[->] (0,0.500000) -- (0.250000,0.500000);
    \draw[->] (0,0.666667) -- (0.222222,0.666667);
    \draw[->] (0,0.833333) -- (0.138889,0.833333);
    \node[anchor=east] at (0,0.5) {\begin{tabular}{c}
      parabolic \\
      inflow
    \end{tabular} };

    \draw[->] (1,0.75) -- (1.25,0.75);
    \draw[->] (1,0.5) -- (1.25,0.5);
    \draw[->] (1,0.25) -- (1.25,0.25);
    \node[anchor=west] at (1.25,0.5) {\begin{tabular}{c}
      do-nothing \\
      outflow
    \end{tabular} };

    \node[anchor=north east] at (0,0) {$(0,0)$};
    \node[anchor=south west] at (1,1) {$(1,1)$};
  \end{tikzpicture}
  \caption{Domain in 2D used for numerical experiments}
  \label{fig:main-domain}
\end{figure}

We consider three model problems. All three make use of the simple 2D square domain pictured in Figure~\ref{fig:main-domain}. The top and bottom boundaries are designated as walls through which no flow occurs. The left boundary is a parabolic inflow boundary. The right boundary is a do-nothing out-flow boundary. All three problems partition this domain into~$n \times n$ uniform subdomains in which permeabilities are constant. All three problems use constant viscosity~$\visc = \effectivevisc = 10^{-3}$ and a zero forcing term.

\begin{figure}[h!]
  \centering
  \includegraphics[width=8.4cm]{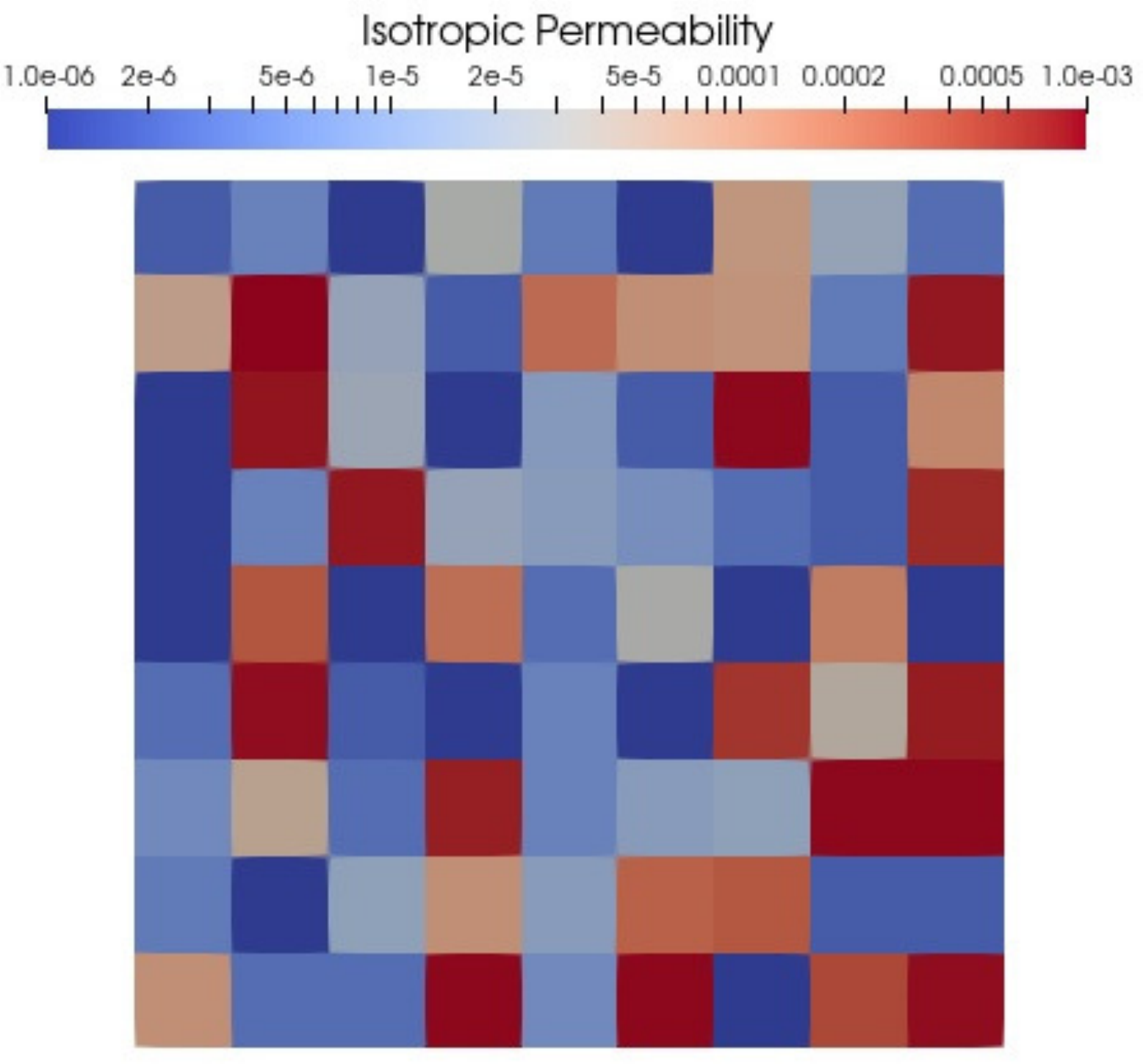}
  \caption{Mean permeabilities for the \emph{iso} problem on a log scale}
  \label{fig:isotropic-permeabilities}
\end{figure}

For the first problem, labeled \emph{iso} in the following discussion, we consider isotropic flow in which the permeability tensor in each subdomain is of the form~$k \dof{I}$ for a scalar~$k > 0$. In this case, we partition the domain into~$9 \times 9$ subdomains. Each subdomain is then partitioned into~$12 \times 12$ regular quadrilateral elements to form a~$108 \times 108$ grid. The permeability coefficient of each subdomain is uniformly distributed in an interval~$[a,b]$ chosen as follows. First, for each subdomain, a value~$c$ is sampled from a~$\mathrm{Beta}(0.5,0.5)$ distribution and then mapped to the interval~$[-6,-3]$. A value~$r$ is uniformly sampled in the interval~$[0.05,0.15]$, and the chosen interval is~$[(1-r)\cdot 10^c, (1+r)\cdot 10^c]$. The parameters of the beta distribution were chosen to skew the permeabilities to slightly favor a mix of high permeability regions (near~$10^{-3}$) and low permeability regions (near~$10^{-6}$). The mean permeabilities for each subdomain are displayed in Figure~\ref{fig:isotropic-permeabilities} on a log scale.

The resulting stochastic space is then of dimension~$N = 81$. The parametrization coincides with the permeabilities in each subdomain, i.e.,~$\xi \in \R^{81}$ with~$\xi_i$ equal to the permeability in subdomain~$i$. The affine decomposition~(\ref{eq:affine-decomposition-A}) consists of~$n_A = N + 1 = 82$ parameter independent matrices as follows:

\begin{equation}
  A_i = \left\{
    \begin{array}{ll}
      A_S & i = 1 \\
      A_D^{i-1} & i > 1 ,
    \end{array} \right. \label{eq:affine-decomposition-A-matrices-iso}
\end{equation}

\noindent where~$A_S$ is the discretization of the Stokes bilinear form~(\ref{eq:bilinear-stokes}), and~$A_D^k$ is the discretization of the Darcy bilinear form~(\ref{eq:bilinear-darcy}) with support only on subdomain~$k$. The mappings~$\{\theta_i^A(\xi)\}_{i=1}^{n_A}$ then take the form

\begin{equation}
  \theta_i^A(\xi) = \left\{
    \begin{array}{ll}
      1 & i = 1 \\
      1/\xi_{i-1} & i > 1 .
    \end{array} \right. \label{eq:affine-decomposition-A-mappings-iso}
\end{equation}

\noindent To apply the Dirichlet boundary conditions, we use the finite element function which interpolates the Dirichlet conditions on the Dirichlet boundary and has zero value on all remaining nodes. Let~$\fe{w}$ denote this finite element function. Then, since we define zero forcing term and do-nothing Neumann conditions, the discretizations of~(\ref{eq:rhs-l1-continuous})--(\ref{eq:rhs-l2-continuous}) are simply

\begin{align}
  \ell_1(\fe{v};\xi) &= -a(\fe{w},\fe{v};\xi) \quad \forall \fe{v} \in \fe{\velocityspace} \label{eq:rhs-l1-discrete-isotropic}\\
  \ell_2(\fe{q}) &= -b(\fe{w},\fe{q}) \quad \forall \fe{q} \in \fe{\pressurespace} \label{eq:rhs-l2-discrete-isotropic}.
\end{align}

\noindent Note that, since~$b$ is parameter-independent and the Dirichlet boundary conditions are parameter-independent, the discretization~(\ref{eq:rhs-l2-discrete-isotropic}) is parameter-independent. However, due to the dependence of~$a$ on the parameter, the discretization~(\ref{eq:rhs-l1-discrete-isotropic}) is parameter-dependent. It admits an affine decomposition of the form~(\ref{eq:affine-decomposition-A-matrices-iso})--(\ref{eq:affine-decomposition-A-mappings-iso}). We note, however, that our choice of~$\fe{w}$ suggests that we need only~$n_f = 10$ parameter-independent vectors, corresponding to the vector arising from the Stokes bilinear form and the 9 Darcy bilinear forms on the subdomains which border the in-flow boundary.

For the second and third problems, we consider anisotropic flow in which the permeability tensor in each subdomain is of the form~$\mathrm{diag}(k_x,k_y)$ for positive scalars~$k_x \neq k_y$. The second problem, labeled \emph{aniso1}, considers the case where~$k_x < k_y$, that is, where vertical flow is favored. The third problem, labeled \emph{aniso2}, considers the case where~$k_x > k_y$, that is, where horizontal flow is favored. Both problems are partitioned into~$6 \times 6$ subdomains, each subdomain partitioned into~$18 \times 18$ elements. This partitioning is to make the size of the problems similar to that of the isotropic problem. Here, the total number of elements is~$108 \times 108$, matching the isotropic problem exactly; and the number of parameters is~$72$, as there are two parameters per subdomain.

\begin{figure*}[h!]
  \centering
  \includegraphics[width=17.4cm]{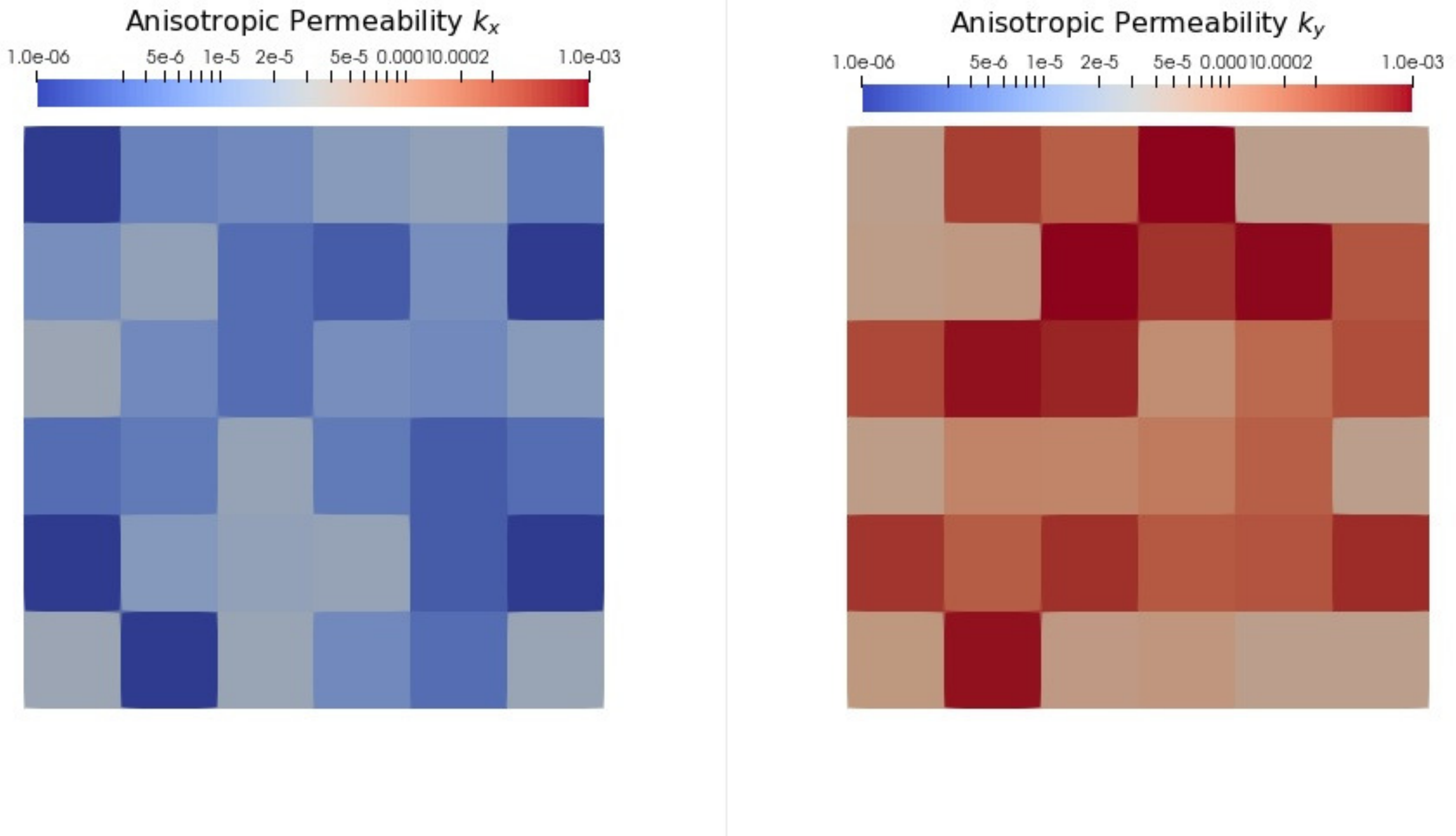}
  \caption{Mean permeabilities for the \emph{aniso1} problem on a log scale}
  \label{fig:aniso1-permeabilities}
\end{figure*}

\begin{figure*}[h!]
  \centering
  \includegraphics[width=17.4cm]{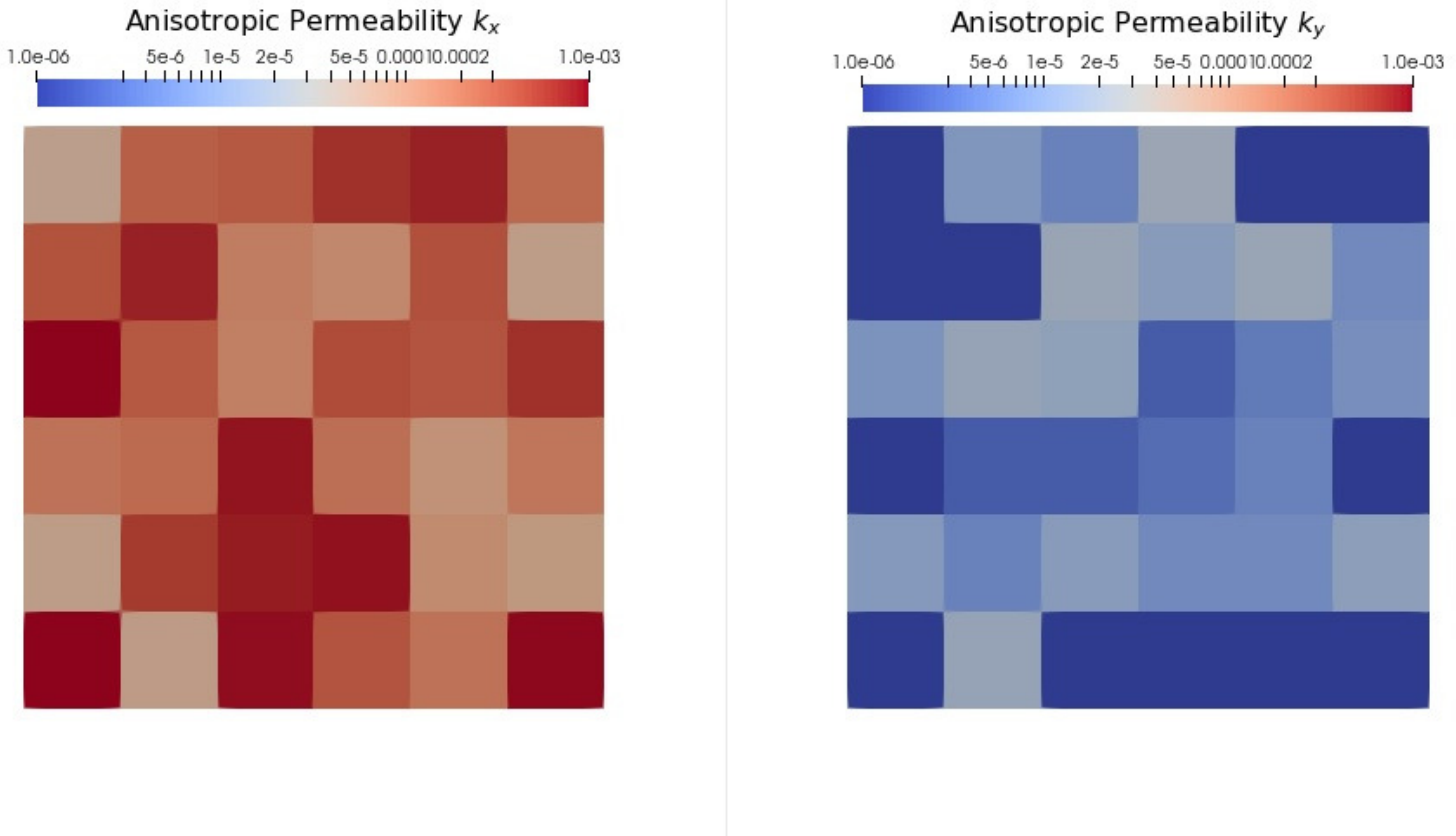}
  \caption{Mean permeabilities for the \emph{aniso2} problem on a log scale}
  \label{fig:aniso2-permeabilities}
\end{figure*}

The distributions on permeabilities are chosen in a manner similar to the isotropic problem. The main difference is that the smaller permeability ($k_x$ in the case of \emph{aniso1},~$k_y$ in the case of \emph{aniso2}) has the beta random variable mapped to~$[-6,-4.75]$, and the larger permeability ($k_y$ in the case of \emph{aniso1},~$k_x$ in the case of \emph{aniso2}) has the beta random variable mapped to~$[-4.25,-3]$. The mean permeabilities for each subdomain and direction are displayed in Figures~\ref{fig:aniso1-permeabilities} and~\ref{fig:aniso2-permeabilities} on a log scale.

The anisotropic problems admit a similar affine decomposition to the isotropic problem. However, now there are two Darcy parameter-independent matrices per subdomain, each discretizations of~(\ref{eq:bilinear-darcy}) with support only from basis functions on the given subdomain and direction. Thus,~$n_A = 73$. Furthermore, we have~$n_f = 7$, as our choice of~$\fe{w}$ only admits non-zero parameter-independent vectors arising from the Stokes bilinear form and then one for each of the~$6$ subdomains bordering the in-flow boundary (as the Dirichlet in-flow conditions are zero in the~$y$-direction).

Table~\ref{tab:finite-element-sizes} compares the sizes of the finite element discretizations of each problem. As can be observed, all three problems have the same dimension for the high-fidelity space, suggesting they should all be of similar difficulty. Condition number estimates, computed using the \texttt{condest} function in \textsc{Matlab}, for the assembled matrices when choosing the mean of each distribution as the parameter are presented in Table~\ref{tab:condition-numbers}. Here, we observe that all three problems are fairly ill-conditioned. We note that the sizes of the high-fidelity problems were chosen to be small enough that we could perform the high-fidelity solves efficiently using a direct solver. The log permeabilities were also chosen to be within the interval~$[-6,-3]$ so as to keep the condition numbers moderate. These condition numbers have an effect on the sharpness of the error bounds, as will be discussed in Section~\ref{sec:rb-numerical}. In short, larger condition numbers imply that the error bounds may not be as sharp, leading to potentially larger reduced bases than necessary.

\begin{table*}[h!]
  \centering
  \begin{tabular}{ccccccc}
    Problem & subdomains & elements & total elements & velocity dof & pressure dof & total dof \\ \hline
    \emph{iso} & $9 \times 9$ & $12 \times 12$ & 11664 & 92880 & 34992 & 127872 \\
    \emph{aniso1} & $6 \times 6$ & $18 \times 18$ & 11664 & 92880 & 34992 & 127872 \\
    \emph{aniso2} & $6 \times 6$ & $18 \times 18$ & 11664 & 92880 & 34992 & 127872
  \end{tabular}
  \caption{Comparison of sizes of the finite element discretization for each model problem}
  \label{tab:finite-element-sizes}
\end{table*}

\begin{table}[h!]
  \centering
  \begin{tabular}{cc}
    Problem & condition number \\ \hline
    \emph{iso} & $2.6894 \times 10^4$ \\
    \emph{aniso1} & $5.0727 \times 10^4$ \\
    \emph{aniso2} & $4.0431 \times 10^3$
  \end{tabular}
  \caption{Comparison of condition number estimates for the assembled matrices with the mean parameter}
  \label{tab:condition-numbers}
\end{table}

Table~\ref{tab:parametrization-sizes} compares the sizes of the stochastic dimension and number of parameter-independent components in the affine decompositions. The isotropic problem is larger in this case due to the use of more subdomains.

\begin{table}[h!]
  \centering
  \begin{tabular}{cccc}
    Problem & stochastic dimension & $n_A$ & $n_f$ \\ \hline
    \emph{iso} & 81 & 82 & 10 \\
    \emph{aniso1} & 72 & 73 & 7 \\
    \emph{aniso2} & 72 & 73 & 7
  \end{tabular}
  \caption{Comparison of parametrization sizes for each model problem}
  \label{tab:parametrization-sizes}
\end{table}

\subsection{SCM} \label{sec:numerical-scm}

For the eigenvalue problems solved during the SCM training (see Algorithm~\ref{alg:scm-training}), we used the LOBPCG method~\cite{Knyazev-LOBPCG-2001} as implemented in the BLOPEX \textsc{Matlab} package~\cite{Knyazev-BLOPEX}. The inf-sup constant, being parameter-independent, required only a one-time eigenvalue computation. The specific eigenvalue problem solved was to find the smallest eigenvalue~$\lambda$ of

\begin{align}
  \dof{B} \dof{M}_\velocityspace^{-1} \dof{B}^T x &= \lambda \dof{M}_\pressurespace x \label{eq:infsup-eigen}
\end{align}

\noindent and choosing~$\beta = \sqrt{\lambda}$. No preconditioner was used. The coercivity and continuity constants are parameter-dependent and require solving for the smallest and largest, respectively, eigenvalues~$\lambda(\xi)$ of the generalized eigenvalue problem

\begin{align}
  \dof{A}(\xi) x &= \lambda(\xi) \dof{M}_\velocityspace x . \label{eq:coercivity-eigen}
\end{align}

\noindent For the coercivity constant, we solved for the smallest eigenvalue~$\lambda(\xi)$ of~(\ref{eq:coercivity-eigen}) using an incomplete Cholesky factorization of~$\dof{A}(\xi)$ as a preconditioner and set~$\alpha(\xi) = \lambda(\xi)$. Computing the continuity constant requires solving for the largest eigenvalue of~(\ref{eq:coercivity-eigen}). For this, we transformed the problem into solving for the smallest eigenvalue~$\mu(\xi)$ of the generalized eigenvalue problem

\begin{align}
  \dof{M}_\velocityspace &= \mu(\xi) \dof{A}(\xi) x . \label{eq:continuity-eigen}
\end{align}

\noindent The continuity constant is then chosen as~$\gamma(\xi) = 1/\mu(\xi)$. All LOBPCG computations were performed with a maximum of 1000 iterations and a tolerance of $10^{-10}$.

We remark that inverting~$\dof{M}_\velocityspace$ is performed during the computation of~(\ref{eq:offline-dual1}), computation of~(\ref{eq:infsup-eigen}), and application of the supremizer~(\ref{eq:supremizer}). For small enough problems, we may compute and store the sparse Cholesky factorization of~$\dof{M}_\velocityspace$. In these cases, we may use these factors to perform an exact solve as a preconditioner for~(\ref{eq:continuity-eigen}). If a preconditioned iterative method is required for inverting~$\dof{M}_\velocityspace$, then we may use the same preconditioner for~(\ref{eq:continuity-eigen}). In this paper, we used the Cholesky factors for the preconditioner.

Each application of SCM requires solving a linear programming problem. For this, we used the \texttt{linprog} function in \textsc{Matlab} with the dual-simplex algorithm.

\begin{figure}[h!]
  \centering
  \includegraphics[width=8.4cm]{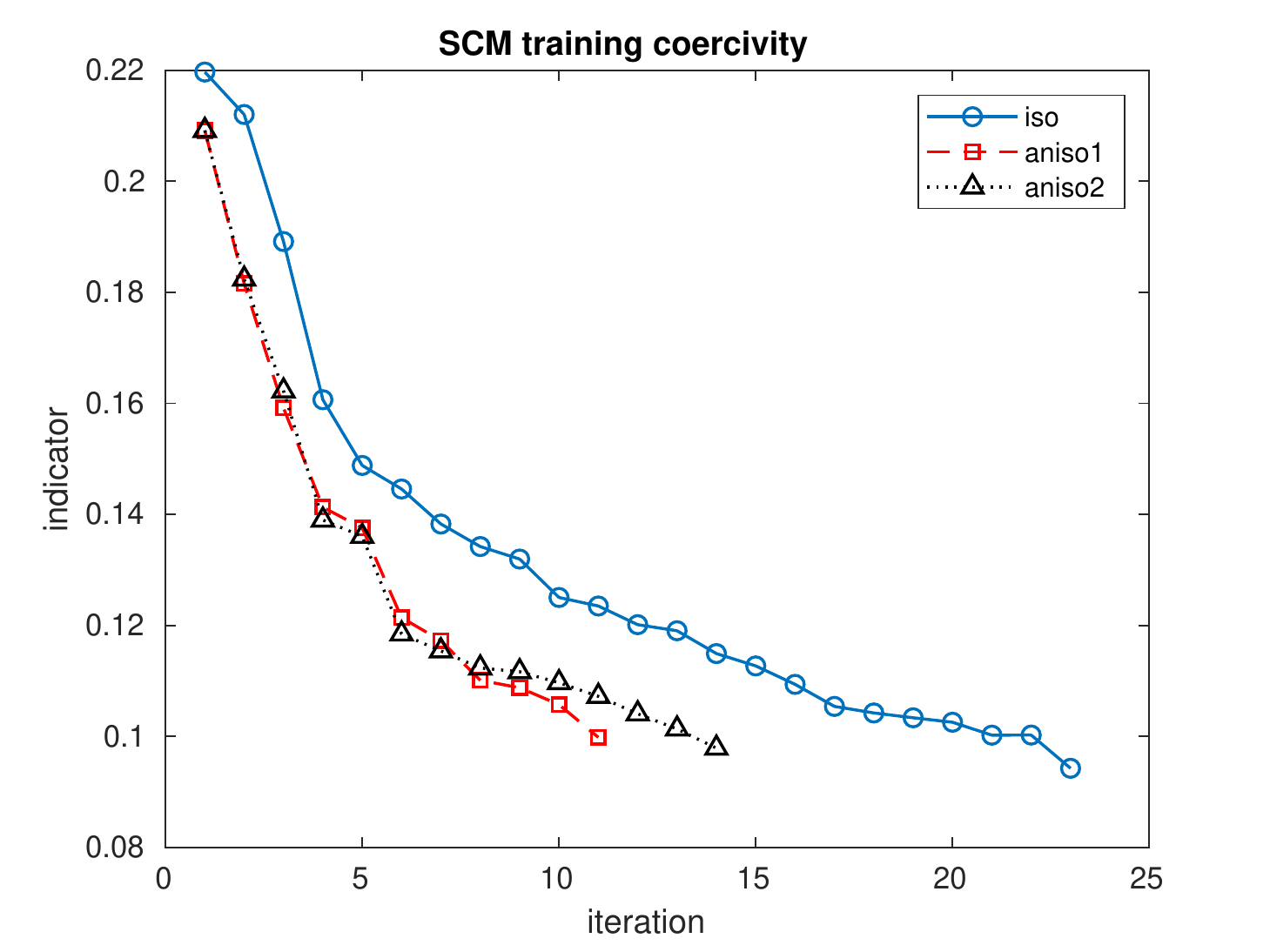}
  \caption{Largest indicators for coercivity during each iteration of SCM training using a tolerance of~$\epsilon^{\mathrm{SCM}} = 0.1$ with a training size of 50000 points and parameters~$M_E = M_P = 100$}
  \label{fig:scm-iterations-coercivity}
\end{figure}

\begin{figure}[h!]
  \centering
  \includegraphics[width=8.4cm]{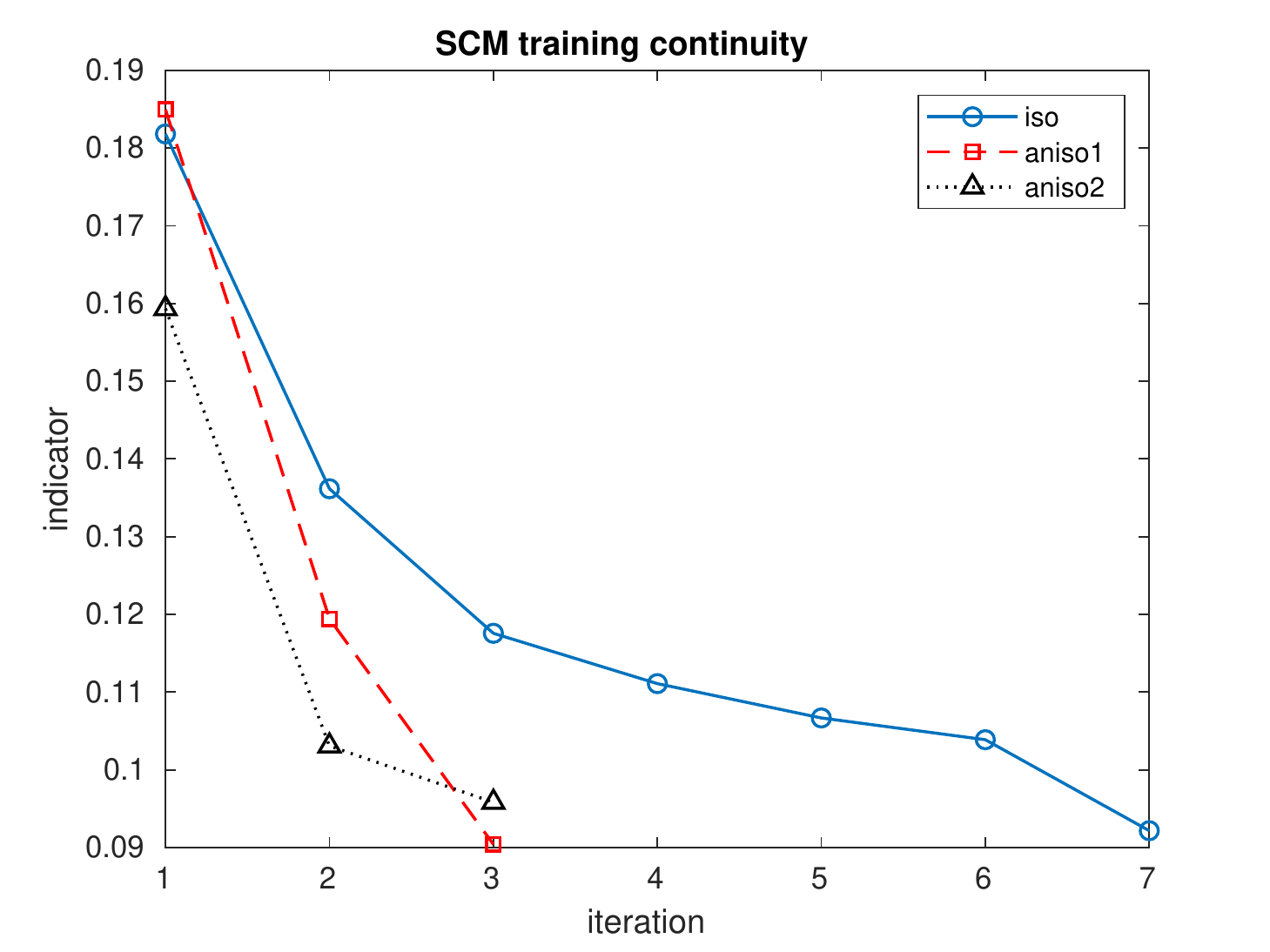}
  \caption{Largest indicators for continuity during each iteration of SCM training using a tolerance of~$\epsilon^{\mathrm{SCM}} = 0.1$ with a training size of 50000 points and parameters~$M_E = M_P = 100$}
  \label{fig:scm-iterations-continuity}
\end{figure}

For the SCM training of each model problem, we generated a Halton set of size 50000. We chose~$M_E = M_P = 100$ and the tolerance~$\epsilon^{\mathrm{SCM}} = 0.1$. Table~\ref{tab:scm-iterations} summarizes the number of iterations required for each model problem to attain the prescribed tolerance. As can be observed, the training terminates fairly quickly in all three cases, with the coercivity constant being slightly harder than continuity and the isotropic case being slightly harder than the anisotropic cases. The isotropic case being more difficult is likely due to the larger number of parameters. Figures~\ref{fig:scm-iterations-coercivity} and~\ref{fig:scm-iterations-continuity} display the largest indicator~(\ref{eq:scm-indicator}) over each iteration of training. The indicators for the \emph{iso} problem are plotted with the solid blue line with circles; the indicators for the \emph{aniso1} problem are plotted with the dashed red line with squares; the indicators for the \emph{aniso2} problem are plotted with the dotted black line with triangles.

\begin{table}[h!]
  \centering
  \begin{tabular}{ccc}
    Problem & coercivity & continuity \\ \hline
    \emph{iso} & 23 & 7 \\
    \emph{aniso1} & 11 & 3 \\
    \emph{aniso2} & 12 & 3
  \end{tabular}
  \caption{Number of SCM training iterations to obtain a tolerance of~$\epsilon^{\mathrm{SCM}} = 0.1$ with a training size of 50000 points and parameters~$M_E = M_P = 100$}
  \label{tab:scm-iterations}
\end{table}

\subsection{Reduced Basis ANOVA} \label{sec:rb-numerical}

To study the effectiveness of reduced basis ANOVA on the three model problems, we compute reduced basis ANOVA approximations using the following tensor product of reduced basis and ANOVA tolerances. For the reduced basis tolerances, we consider the set~$\{1,0.1,0.01\}$. For the ANOVA tolerances, we consider the set~$\{10^{-4},10^{-5},10^{-6}\}$. Thus, we performed a total of 9 reduced basis ANOVA computations for each problem. In all problems, we used all ANOVA directions up to level 1 and adaptively selected directions for higher levels. In all problems, the adaptive selection of ANOVA directions terminated at level 2. We used Gauss-Legendre quadrature with polynomial order 5 for generating collocation points. We compare each reduced basis ANOVA computation with the results from a Monte Carlo simulation performed using~$10^6$ points generated from a Halton set. We assume that the moments computed using the Monte Carlo are accurate enough to be used as the true moments so that we may use it to approximate the error in the moments from the reduced basis ANOVA computations.

Table~\ref{tab:L2-errors-iso} summarizes the~$L^2$ norms for the errors in mean and variance for the~$\emph{iso}$ problem, where the error here is taken as the difference between the reduced basis ANOVA approximation and the Monte Carlo simulation. These values are normalized by the norm of the Monte Carlo approximation so as to provide relative errors. The table shows the errors in the velocity, pressure, and combined mean and variances. It is clear that, as the reduced basis tolerance is reduced, better approximations to the mean and variance are obtained. A reduced basis tolerance of~$10^{-2}$ results in excellent approximations even for an ANOVA tolerance of~$10^{-4}$. However, there appear to be little gains from reducing the ANOVA tolerance. Accuracy improves by increasing the tolerance from~$10^{-4}$ to~$10^{-5}$; however, increasing the tolerance from~$10^{-5}$ to~$10^{-6}$ yields no improvement. Tables~\ref{tab:L2-errors-aniso1} and~\ref{tab:L2-errors-aniso2} present the same information for the \emph{aniso1} and \emph{aniso2} problems, respectively. In these cases, we observe the same trend of increased accuracy when reducing the reduced basis tolerance and stagnating accuracy when reducing the ANOVA tolerance.

\begin{table*}[h!]
  \centering
  \resizebox{\textwidth}{!}{%
  \begin{tabular}{lr|ccc|ccc|ccc|c}
    \cline{3-11}
        \, & \, & \multicolumn{3}{c|}{$\epsilon^{\mathrm{A}}$} & \multicolumn{3}{c|}{$\epsilon^{\mathrm{A}}$} & \multicolumn{3}{c|}{$\epsilon^{\mathrm{A}}$} & \, \\
        \, &  \, &  $10^{-4}$ &    $10^{-5}$ &    $10^{-6}$ &    $10^{-4}$ &    $10^{-5}$ &    $10^{-6}$ &    $10^{-4}$ &    $10^{-5}$ &    $10^{-6}$ & \, \\
        \hline
        \multicolumn{1}{ |l }{\multirow{3}{*}{$\epsilon^{\mathrm{RB}}$}} &    1 & $4.89 \times 10^{-4}$ & $4.50 \times 10^{-4}$ & $4.50 \times 10^{-4}$ & $3.27 \times 10^{-4}$ & $2.42 \times 10^{-4}$ & $2.42 \times 10^{-4}$ & $3.28 \times 10^{-4}$ & $2.43 \times 10^{-4}$ & $2.43 \times 10^{-4}$ & \multicolumn{1}{c|}{\multirow{3}{*}{mean}} \\
        \multicolumn{1}{ |l}{} &  0.1 & $3.23 \times 10^{-4}$ & $2.10 \times 10^{-4}$ & $2.14 \times 10^{-4}$ & $3.78 \times 10^{-5}$ & $3.65 \times 10^{-5}$ & $3.66 \times 10^{-5}$ & $3.88 \times 10^{-5}$ & $3.70 \times 10^{-5}$ & $3.71 \times 10^{-5}$ & \multicolumn{1}{c|}{} \\
        \multicolumn{1}{ |l}{} & 0.01 & $1.11 \times 10^{-4}$ & $9.82 \times 10^{-5}$ & $9.97 \times 10^{-5}$ & $6.63 \times 10^{-6}$ & $3.61 \times 10^{-6}$ & $4.34 \times 10^{-6}$ & $7.29 \times 10^{-6}$ & $4.51 \times 10^{-6}$ & $5.13 \times 10^{-6}$ & \multicolumn{1}{c|}{} \\
        \hline
        \multicolumn{1}{ |l }{\multirow{3}{*}{$\epsilon^{\mathrm{RB}}$}}&    1 & $3.41 \times 10^{-1}$ & $3.22 \times 10^{-1}$ & $3.22 \times 10^{-1}$ & $2.41 \times 10^{-2}$ & $1.97 \times 10^{-2}$ & $1.97 \times 10^{-2}$ & $2.41 \times 10^{-2}$ & $1.98 \times 10^{-2}$ & $1.98 \times 10^{-2}$ & \multicolumn{1}{c|}{\multirow{3}{*}{variance}} \\
        \multicolumn{1}{ |l}{} &  0.1 & $2.03 \times 10^{-1}$ & $9.42 \times 10^{-2}$ & $9.40 \times 10^{-2}$ & $2.71 \times 10^{-3}$ & $2.31 \times 10^{-3}$ & $2.31 \times 10^{-3}$ & $2.73 \times 10^{-3}$ & $2.31 \times 10^{-3}$ & $2.32 \times 10^{-3}$ & \multicolumn{1}{c|}{} \\
        \multicolumn{1}{ |l}{} & 0.01 & $5.49 \times 10^{-3}$ & $5.07 \times 10^{-3}$ & $4.93 \times 10^{-3}$ & $7.11 \times 10^{-4}$ & $9.57 \times 10^{-5}$ & $1.01 \times 10^{-4}$ & $7.11 \times 10^{-4}$ & $9.61 \times 10^{-5}$ & $1.01 \times 10^{-4}$ & \multicolumn{1}{c|}{} \\
        \hline
        \, & \, & \multicolumn{3}{c|}{velocity} & \multicolumn{3}{c|}{pressure} & \multicolumn{3}{c|}{combined} & \, \\
    \cline{3-11}
  \end{tabular}
  }
  \caption{$L^2$-norm of moment errors for reduced basis ANOVA estimates of the \emph{iso} problem}
  \label{tab:L2-errors-iso}
\end{table*}

\begin{table*}[h!]
  \centering
  \resizebox{\textwidth}{!}{%
  \begin{tabular}{lr|ccc|ccc|ccc|c}
    \cline{3-11}
    \, & \, & \multicolumn{3}{c|}{$\epsilon^{\mathrm{A}}$} & \multicolumn{3}{c|}{$\epsilon^{\mathrm{A}}$} & \multicolumn{3}{c|}{$\epsilon^{\mathrm{A}}$} & \, \\
    \, & \,      &    $10^{-4}$ &    $10^{-5}$ &    $10^{-6}$ &    $10^{-4}$ &    $10^{-5}$ &    $10^{-6}$ &    $10^{-4}$ &    $10^{-5}$ &    $10^{-6}$ & \, \\
    \hline
    \multicolumn{1}{ |l }{\multirow{3}{*}{$\epsilon^{\mathrm{RB}}$}}    &    1 & $4.49 \times 10^{-4}$ & $3.95 \times 10^{-4}$ & $3.95 \times 10^{-4}$ & $7.70 \times 10^{-5}$ & $4.65 \times 10^{-5}$ & $4.66 \times 10^{-5}$ & $7.74 \times 10^{-5}$ & $4.70 \times 10^{-5}$ & $4.71 \times 10^{-5}$ & \multicolumn{1}{c|}{\multirow{3}{*}{mean}} \\
    \multicolumn{1}{ |l}{}    &  0.1 & $2.38 \times 10^{-4}$ & $2.45 \times 10^{-4}$ & $2.46 \times 10^{-4}$ & $9.48 \times 10^{-6}$ & $1.16 \times 10^{-5}$ & $1.25 \times 10^{-5}$ & $1.03 \times 10^{-5}$ & $1.24 \times 10^{-5}$ & $1.32 \times 10^{-5}$ & \multicolumn{1}{c|}{} \\
    \multicolumn{1}{ |l}{}    & 0.01 & $9.15 \times 10^{-5}$ & $7.80 \times 10^{-5}$ & $7.76 \times 10^{-5}$ & $3.26 \times 10^{-6}$ & $2.48 \times 10^{-6}$ & $2.46 \times 10^{-6}$ & $3.62 \times 10^{-6}$ & $2.82 \times 10^{-6}$ & $2.80 \times 10^{-6}$ & \multicolumn{1}{c|}{} \\
    \hline
    \multicolumn{1}{ |l }{\multirow{3}{*}{$\epsilon^{\mathrm{RB}}$}}&    1 & $9.35 \times 10^{-2}$ & $6.82 \times 10^{-2}$ & $6.82 \times 10^{-2}$ & $9.51 \times 10^{-3}$ & $8.28 \times 10^{-3}$ & $8.22 \times 10^{-3}$ & $9.53 \times 10^{-3}$ & $8.29 \times 10^{-3}$ & $8.23 \times 10^{-3}$ & \multicolumn{1}{c|}{\multirow{3}{*}{variance}} \\
    \multicolumn{1}{ |l}{}    &  0.1 & $1.06 \times 10^{-2}$ & $1.03 \times 10^{-2}$ & $1.04 \times 10^{-2}$ & $1.27 \times 10^{-3}$ & $6.07 \times 10^{-4}$ & $6.04 \times 10^{-4}$ & $1.27 \times 10^{-3}$ & $6.10 \times 10^{-4}$ & $6.07 \times 10^{-4}$ & \multicolumn{1}{c|}{} \\
    \multicolumn{1}{ |l}{}    & 0.01 & $2.74 \times 10^{-3}$ & $2.02 \times 10^{-3}$ & $2.00 \times 10^{-3}$ & $9.99 \times 10^{-4}$ & $7.72 \times 10^{-5}$ & $7.72 \times 10^{-5}$ & $9.99 \times 10^{-4}$ & $7.81 \times 10^{-5}$ & $7.81 \times 10^{-5}$ & \multicolumn{1}{c|}{} \\
    \hline
    \, & \, & \multicolumn{3}{c|}{velocity} & \multicolumn{3}{c|}{pressure} & \multicolumn{3}{c|}{combined} & \, \\
    \cline{3-11}
  \end{tabular}
  }
  \caption{$L^2$-norm of moment errors for reduced basis ANOVA estimates of the \emph{aniso1} problem}
  \label{tab:L2-errors-aniso1}
\end{table*}

\begin{table*}[h!]
  \centering
  \resizebox{\textwidth}{!}{%
  \begin{tabular}{lr|ccc|ccc|ccc|c}
    \cline{3-11}
    \, & \, & \multicolumn{3}{c|}{$\epsilon^{\mathrm{A}}$} & \multicolumn{3}{c|}{$\epsilon^{\mathrm{A}}$} & \multicolumn{3}{c|}{$\epsilon^{\mathrm{A}}$} & \, \\
    \, & \, & $10^{-4}$ & $10^{-5}$ & $10^{-6}$ & $10^{-4}$ & $10^{-5}$ & $10^{-6}$ & $10^{-4}$ & $10^{-5}$ & $10^{-6}$ & \, \\
    \cline{1-12}
    \multicolumn{1}{ |l }{\multirow{3}{*}{$\epsilon^{\mathrm{RB}}$}} &    1 & $1.18 \times 10^{-4}$ & $1.06 \times 10^{-4}$ & $1.06 \times 10^{-4}$ & $1.76 \times 10^{-4}$ & $1.57 \times 10^{-4}$ & $1.57 \times 10^{-4}$ & $1.70 \times 10^{-4}$ & $1.51 \times 10^{-4}$ & $1.51 \times 10^{-4}$ & \multicolumn{1}{c|}{\multirow{3}{*}{mean}} \\
    \multicolumn{1}{ |l}{} &  0.1 & $3.79 \times 10^{-5}$ & $4.24 \times 10^{-5}$ & $4.24 \times 10^{-5}$ & $2.47 \times 10^{-5}$ & $2.42 \times 10^{-5}$ & $2.42 \times 10^{-5}$ & $2.67 \times 10^{-5}$ & $2.71 \times 10^{-5}$ & $2.71 \times 10^{-5}$ & \multicolumn{1}{c|}{} \\
    \multicolumn{1}{ |l}{} & 0.01 & $2.51 \times 10^{-5}$ & $1.48 \times 10^{-5}$ & $1.49 \times 10^{-5}$ & $1.73 \times 10^{-5}$ & $9.29 \times 10^{-6}$ & $7.56 \times 10^{-6}$ & $1.85 \times 10^{-5}$ & $1.01 \times 10^{-5}$ & $8.82 \times 10^{-6}$ & \multicolumn{1}{c|}{} \\
    \cline{1-12}
    \multicolumn{1}{ |l }{\multirow{3}{*}{$\epsilon^{\mathrm{RB}}$}}&    1 & $6.19 \times 10^{-2}$ & $4.94 \times 10^{-2}$ & $4.94 \times 10^{-2}$ & $3.33 \times 10^{-2}$ & $2.33 \times 10^{-2}$ & $2.34 \times 10^{-2}$ & $3.34 \times 10^{-2}$ & $2.35 \times 10^{-2}$ & $2.35 \times 10^{-2}$ & \multicolumn{1}{c|}{\multirow{3}{*}{variance}} \\
    \multicolumn{1}{ |l}{} &  0.1 & $3.82 \times 10^{-3}$ & $4.72 \times 10^{-3}$ & $4.73 \times 10^{-3}$ & $1.62 \times 10^{-3}$ & $1.81 \times 10^{-3}$ & $1.83 \times 10^{-3}$ & $1.63 \times 10^{-3}$ & $1.83 \times 10^{-3}$ & $1.84 \times 10^{-3}$ & \multicolumn{1}{c|}{} \\
    \multicolumn{1}{ |l}{} & 0.01 & $2.24 \times 10^{-3}$ & $1.35 \times 10^{-3}$ & $1.58 \times 10^{-3}$ & $8.36 \times 10^{-4}$ & $3.58 \times 10^{-4}$ & $4.03 \times 10^{-4}$ & $8.44 \times 10^{-4}$ & $3.65 \times 10^{-4}$ & $4.12 \times 10^{-4}$ & \multicolumn{1}{c|}{} \\
    \cline{1-12}
    \, & \, & \multicolumn{3}{c|}{velocity} & \multicolumn{3}{c|}{pressure} & \multicolumn{3}{c|}{combined} & \, \\
    \cline{3-11}
  \end{tabular}
  }
  \caption{$L^2$-norm of moment errors for reduced basis ANOVA estimates of the \emph{aniso2} problem}
  \label{tab:L2-errors-aniso2}
\end{table*}

\begin{figure*}[p!]
  \centering
  \begin{tabular}{ccc}
    \includegraphics[width=5cm]{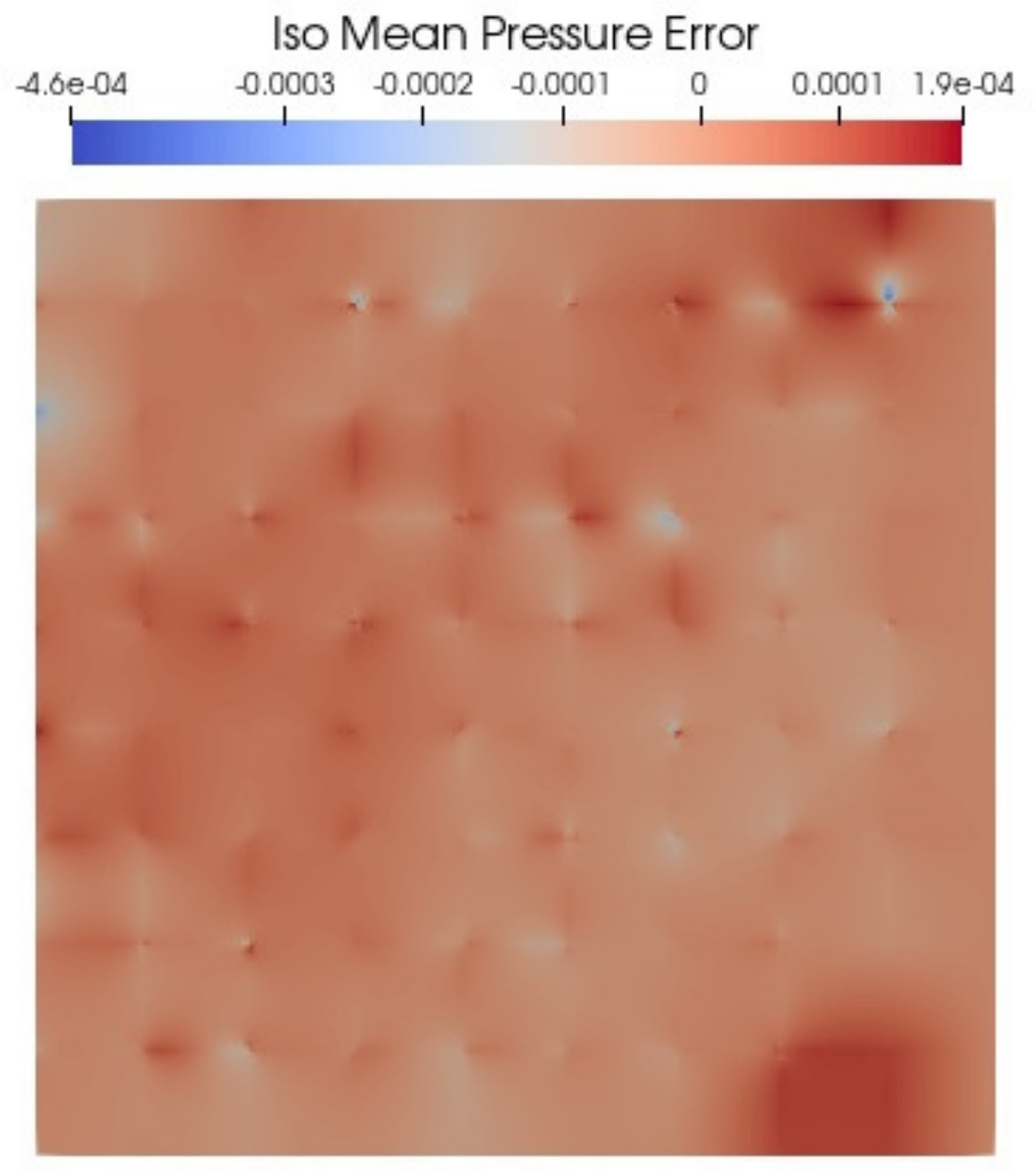} &
    \includegraphics[width=5cm]{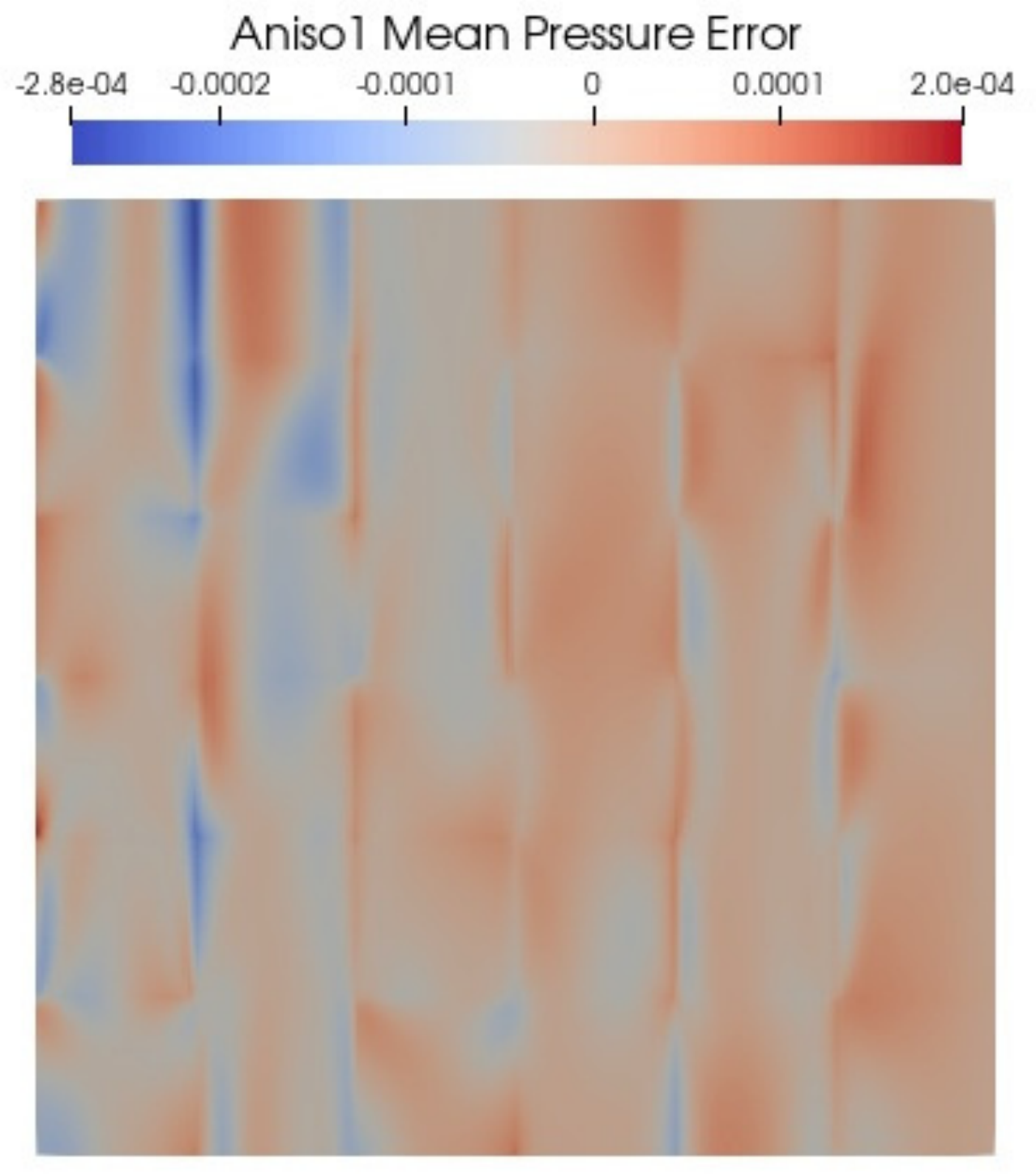} &
    \includegraphics[width=5cm]{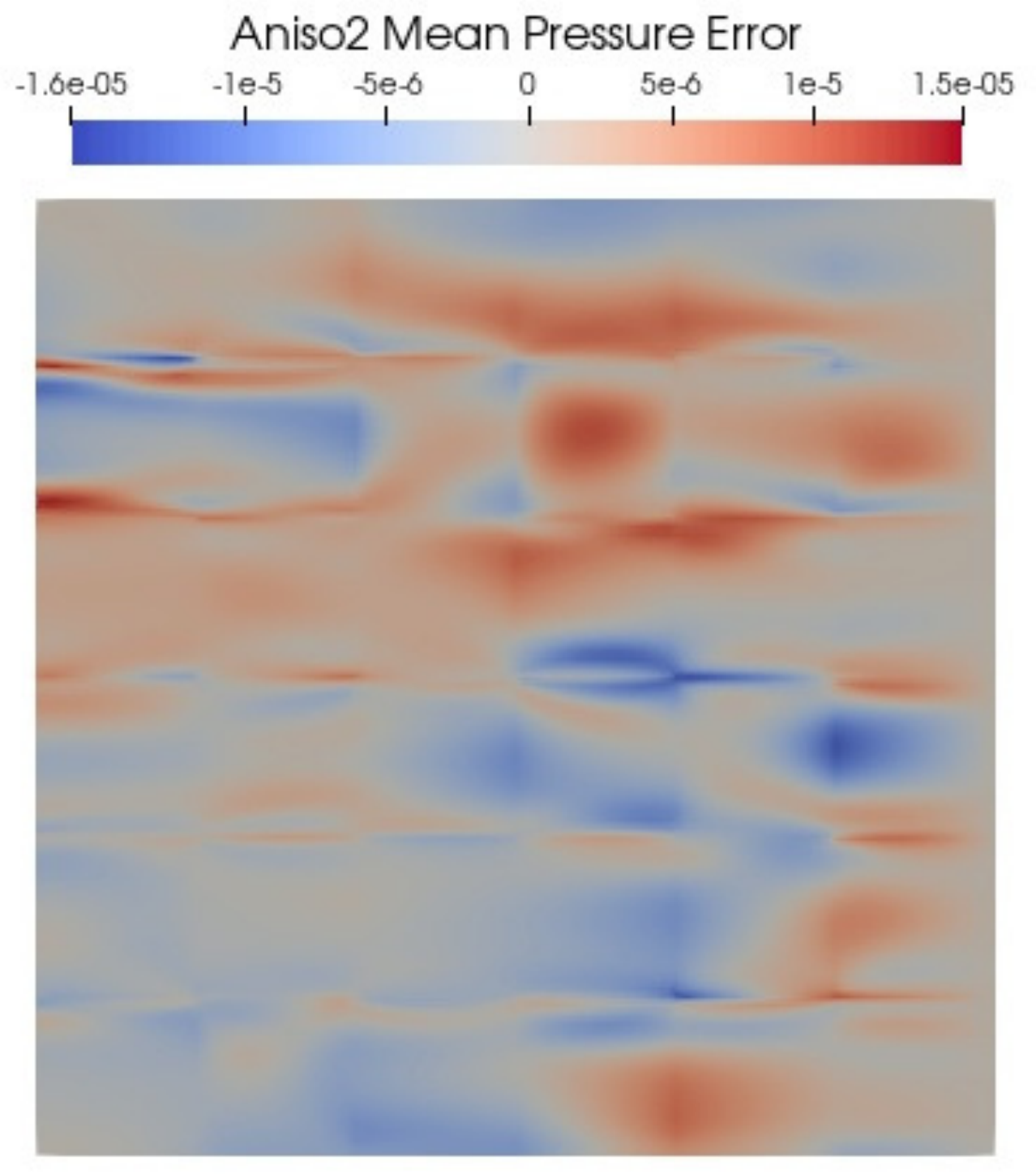} \\
    \includegraphics[width=5cm]{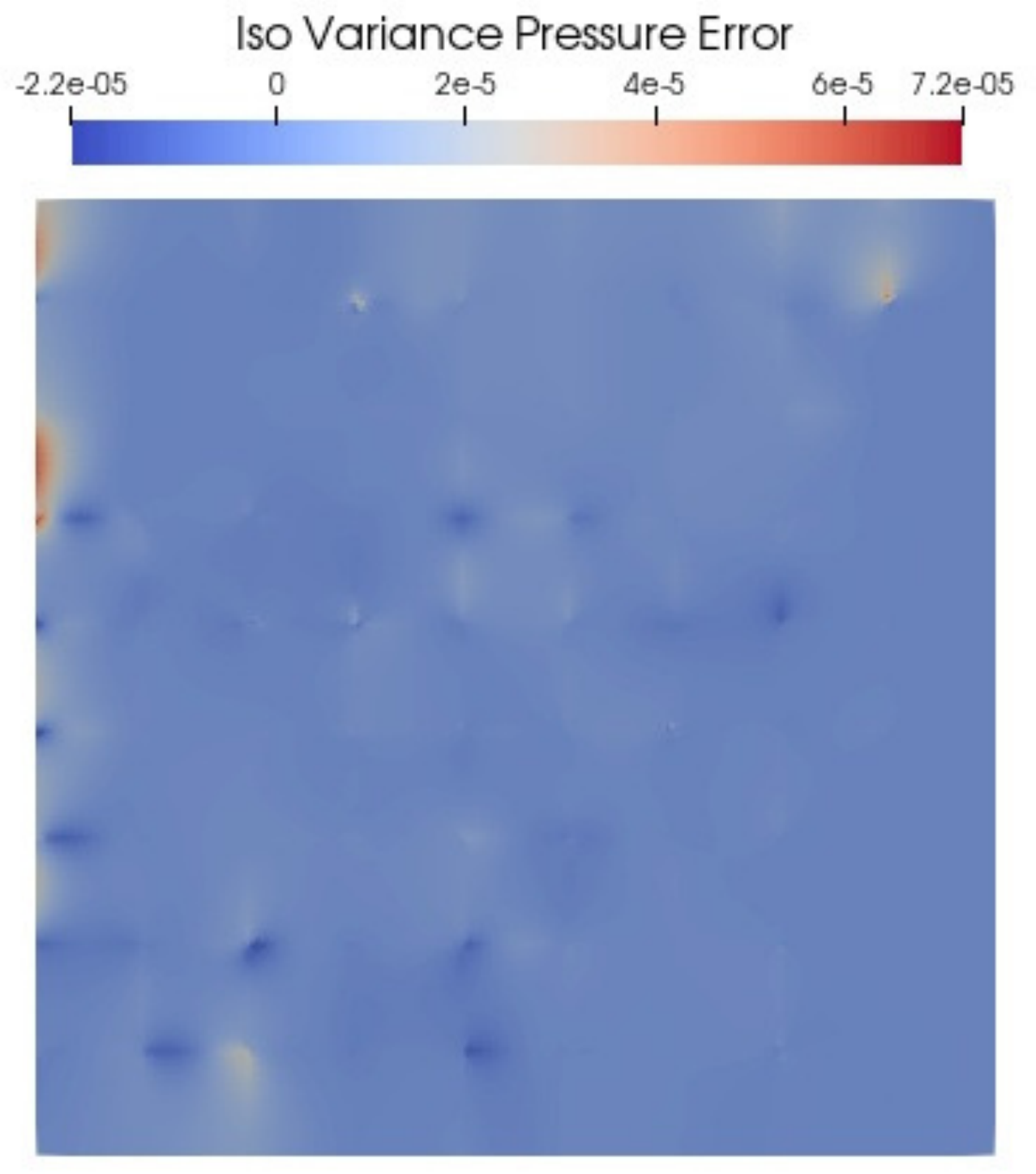} &
    \includegraphics[width=5cm]{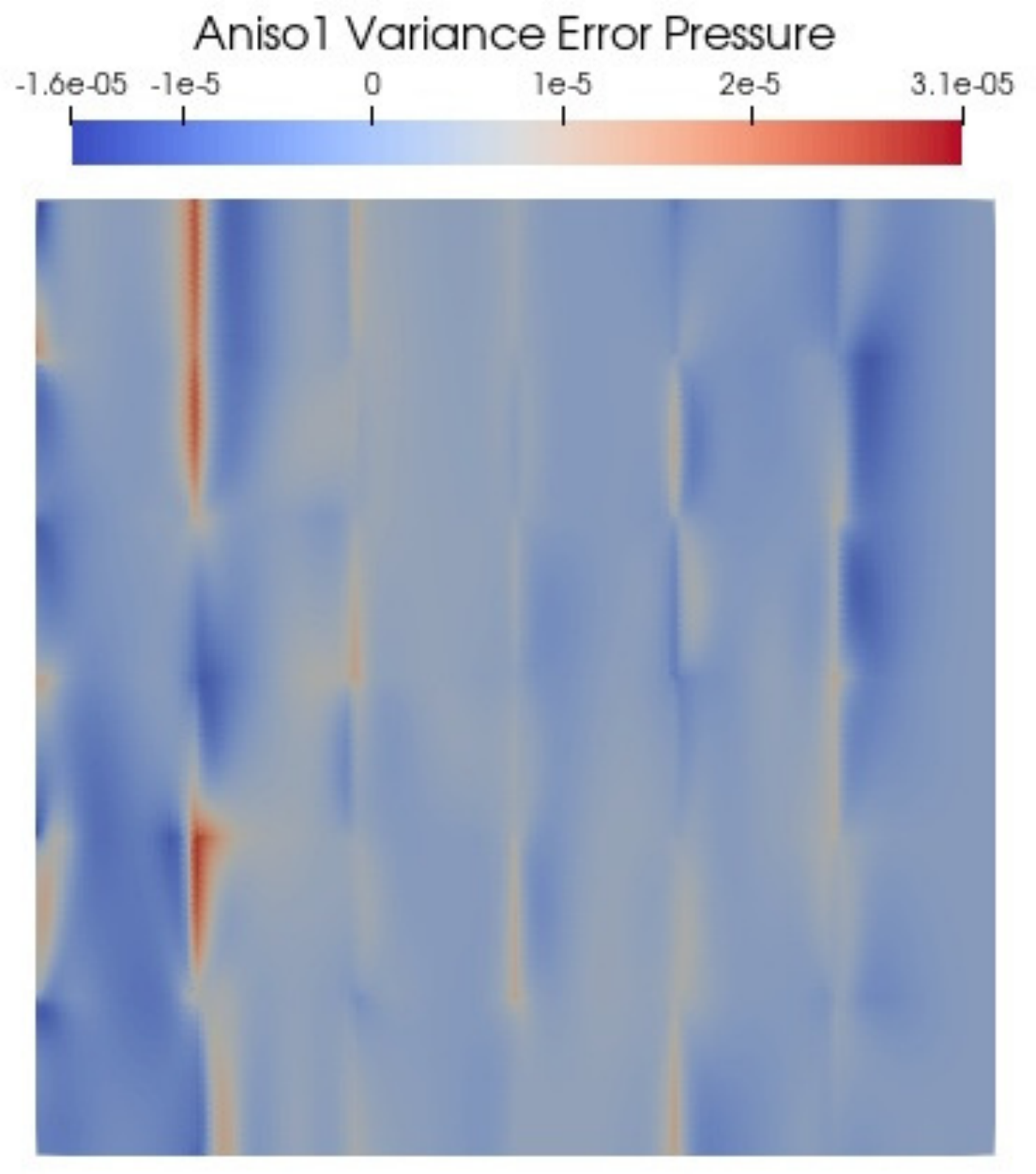} &
    \includegraphics[width=5cm]{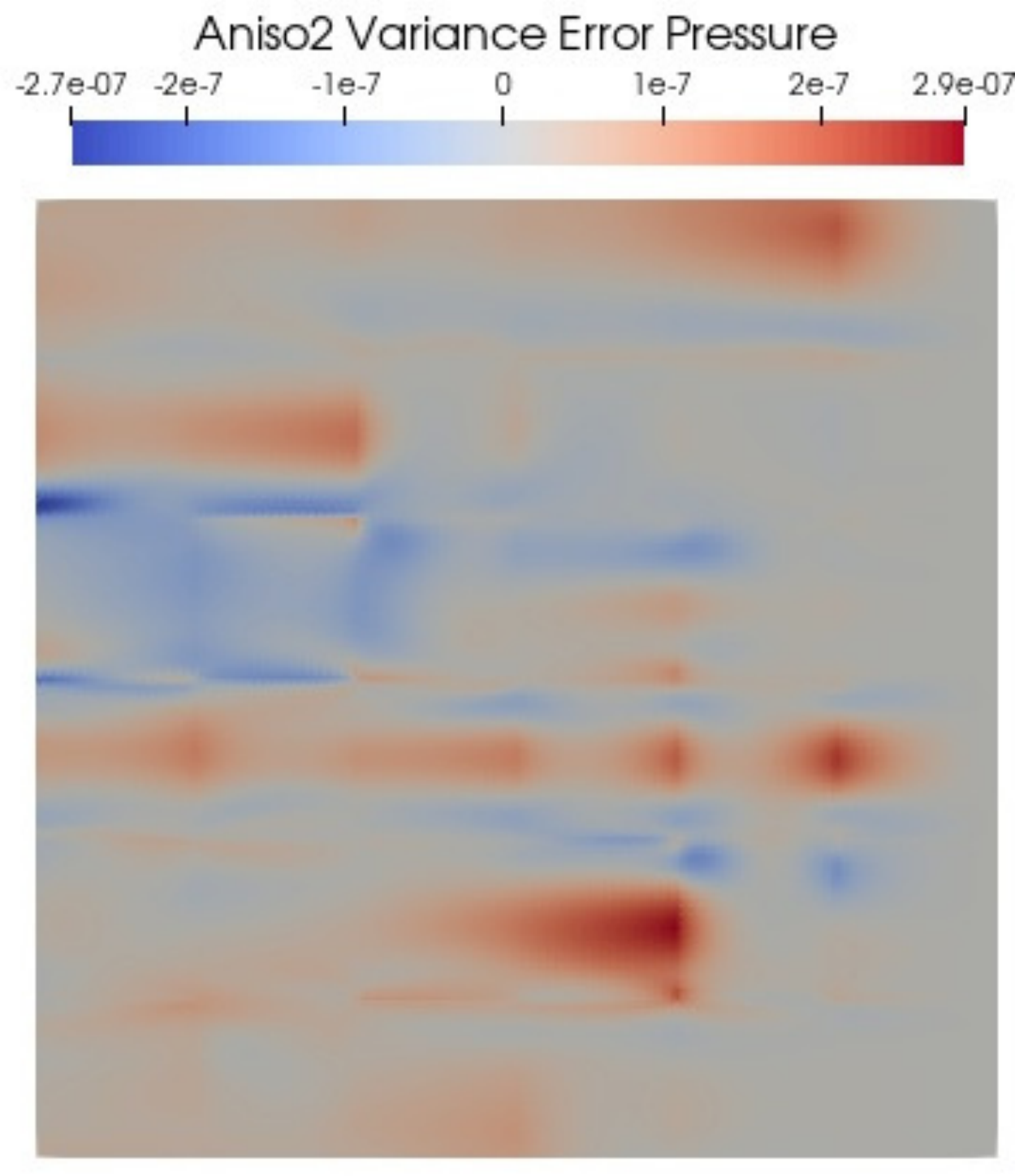}
  \end{tabular}
  \caption{Moment errors in pressure for each of the three problems}
  \label{fig:pressure-mean-and-variance-error}
\end{figure*}

\begin{figure*}[p!]
  \centering
  \begin{tabular}{ccc}
    \includegraphics[width=5cm]{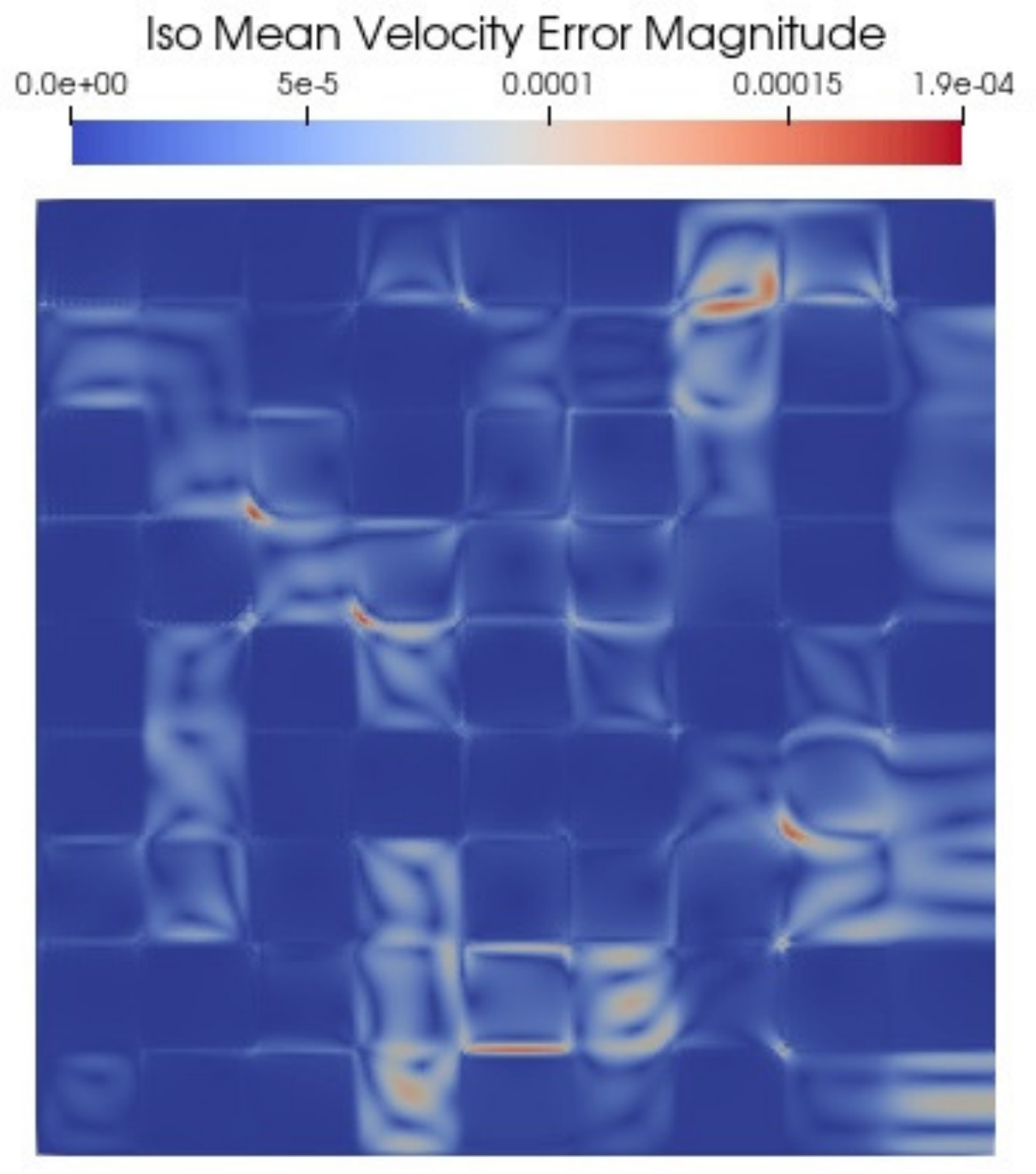} &
    \includegraphics[width=5cm]{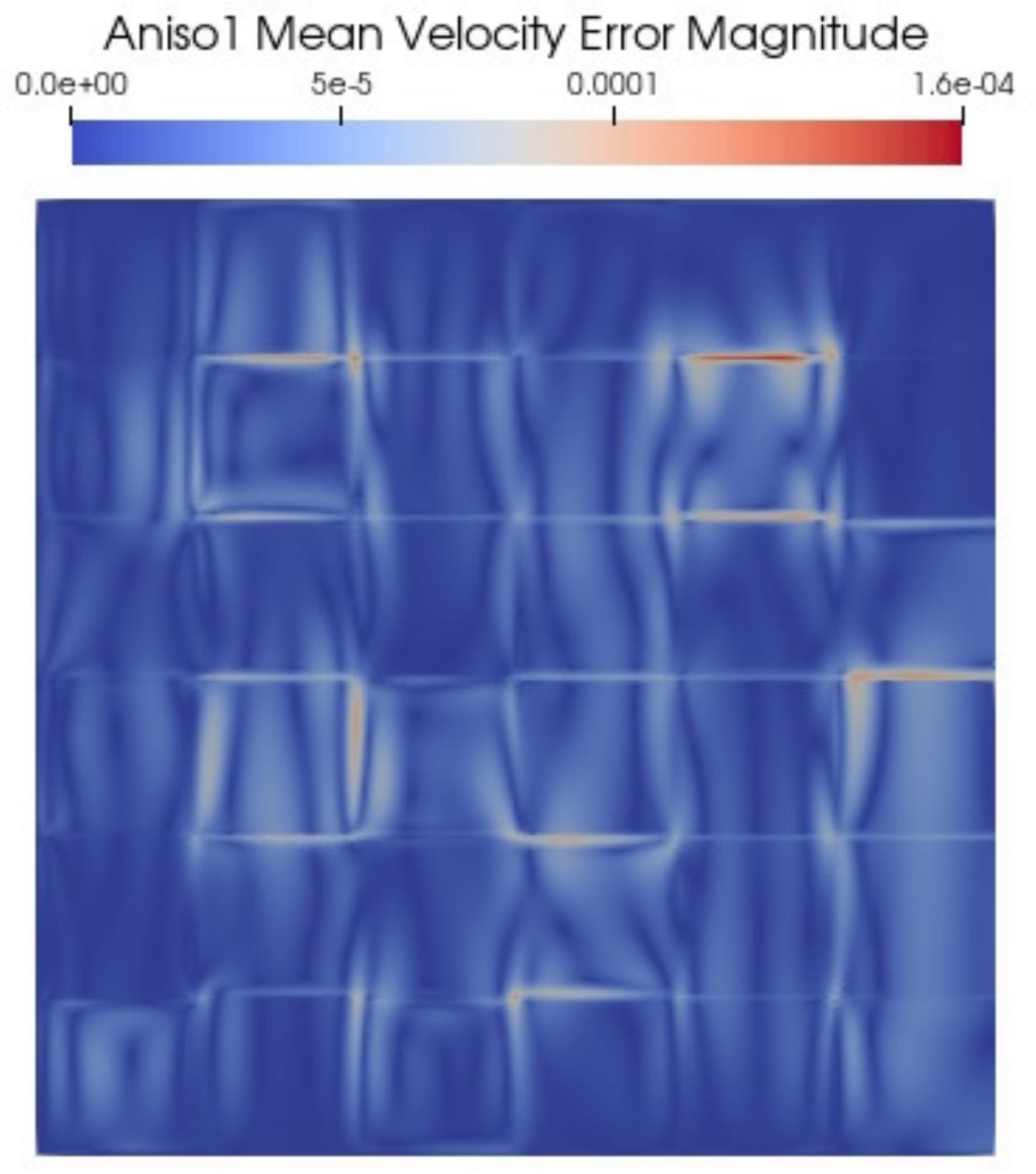} &
    \includegraphics[width=5cm]{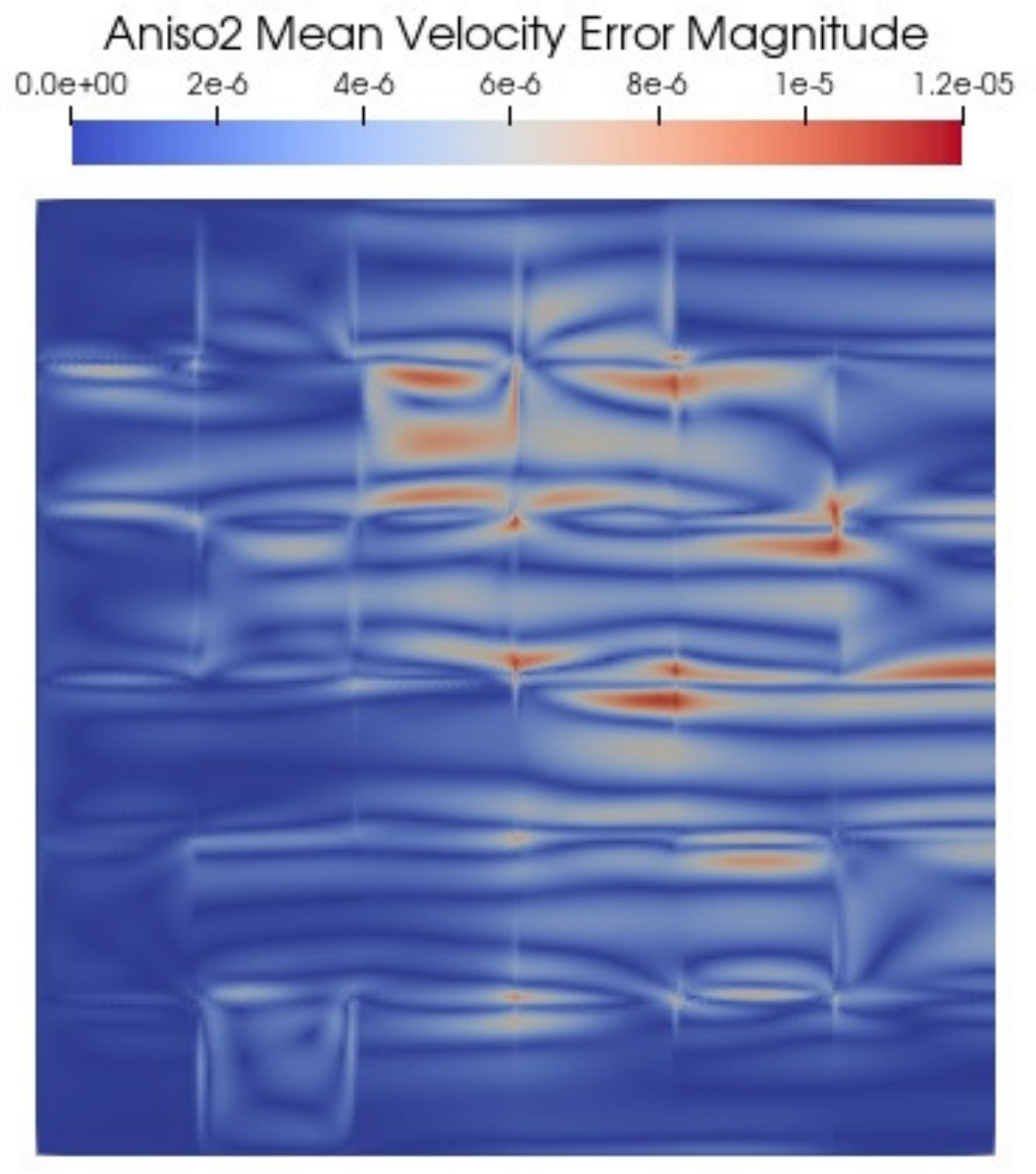} \\
    \includegraphics[width=5cm]{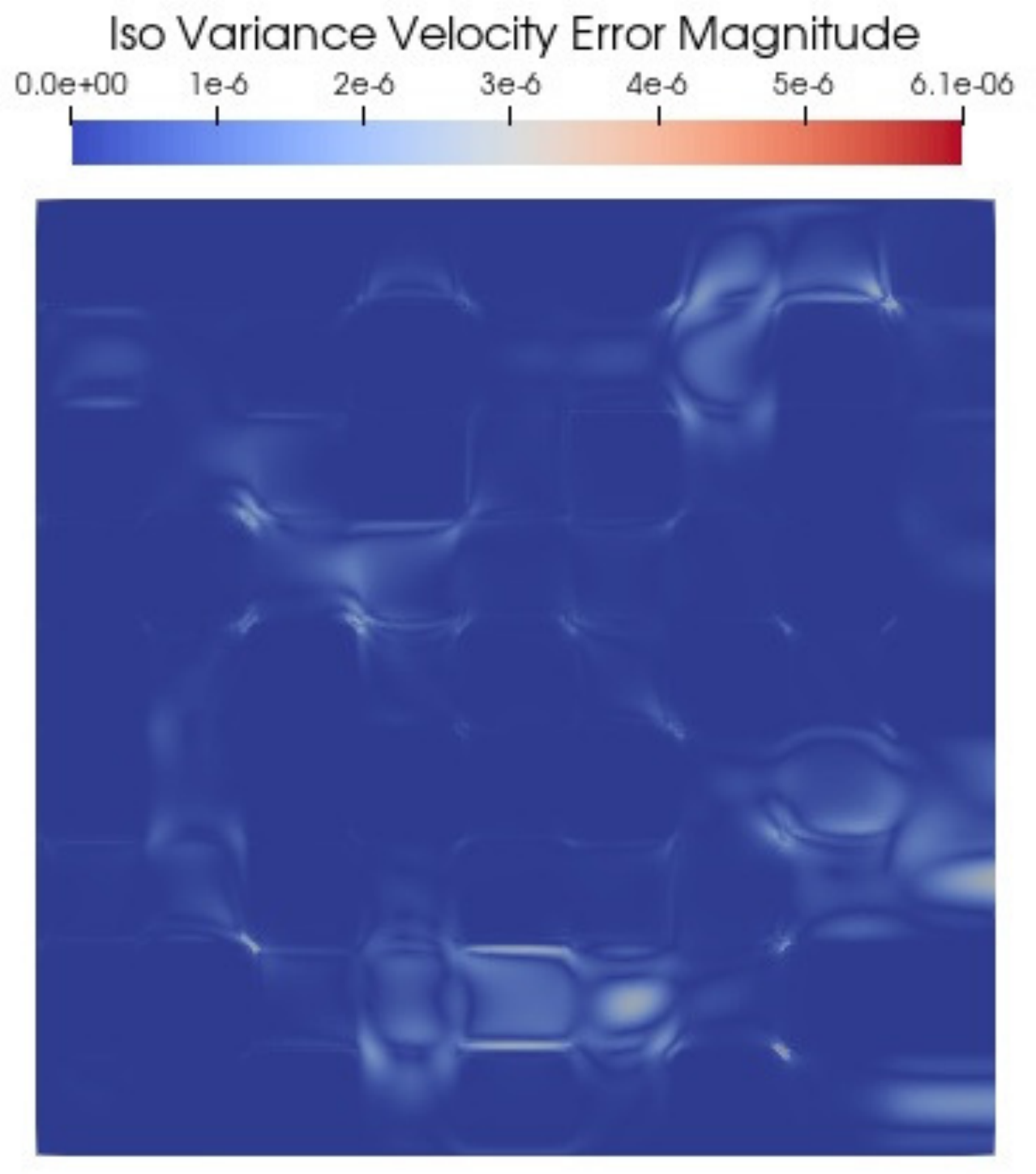} &
    \includegraphics[width=5cm]{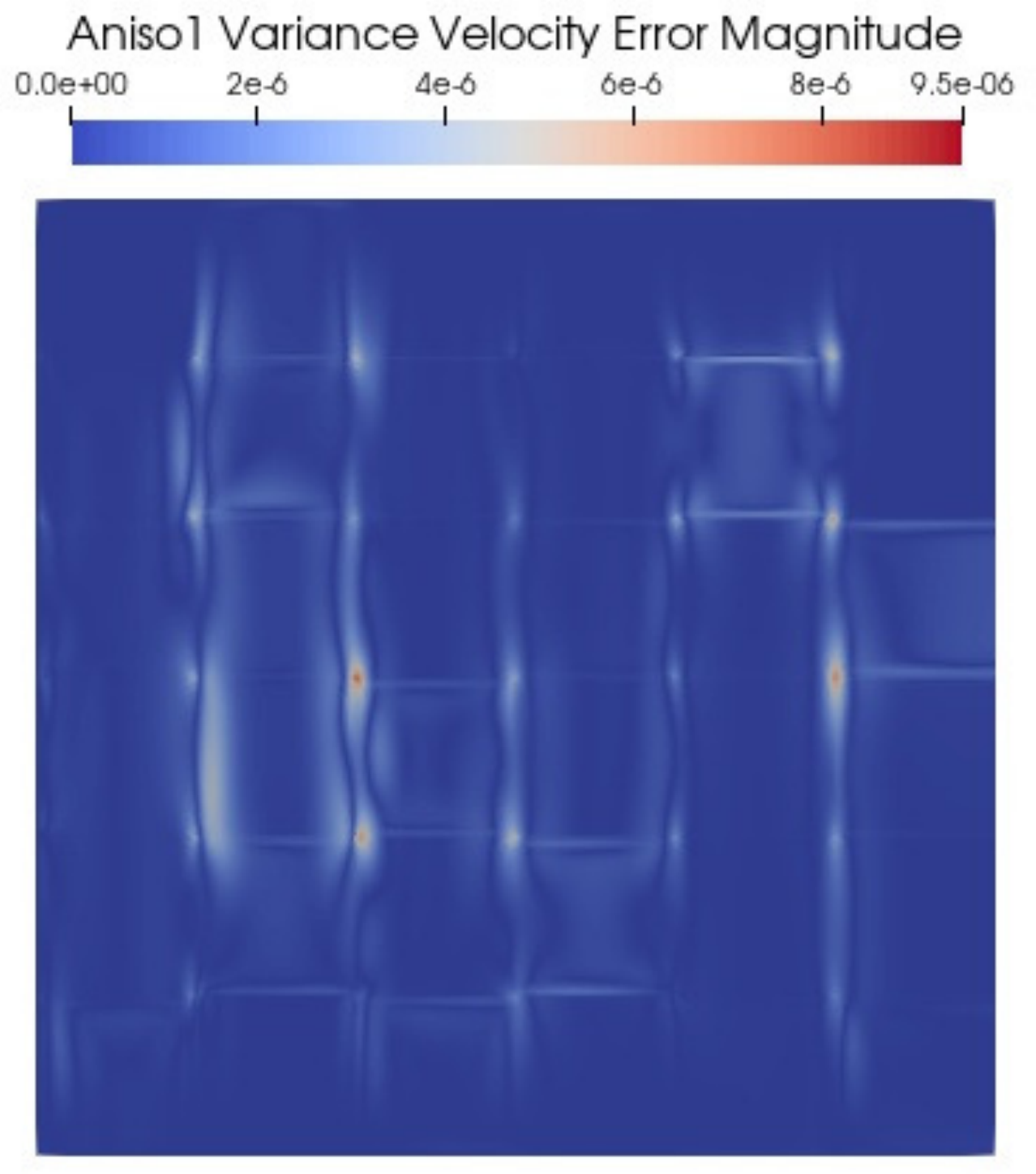} &
    \includegraphics[width=5cm]{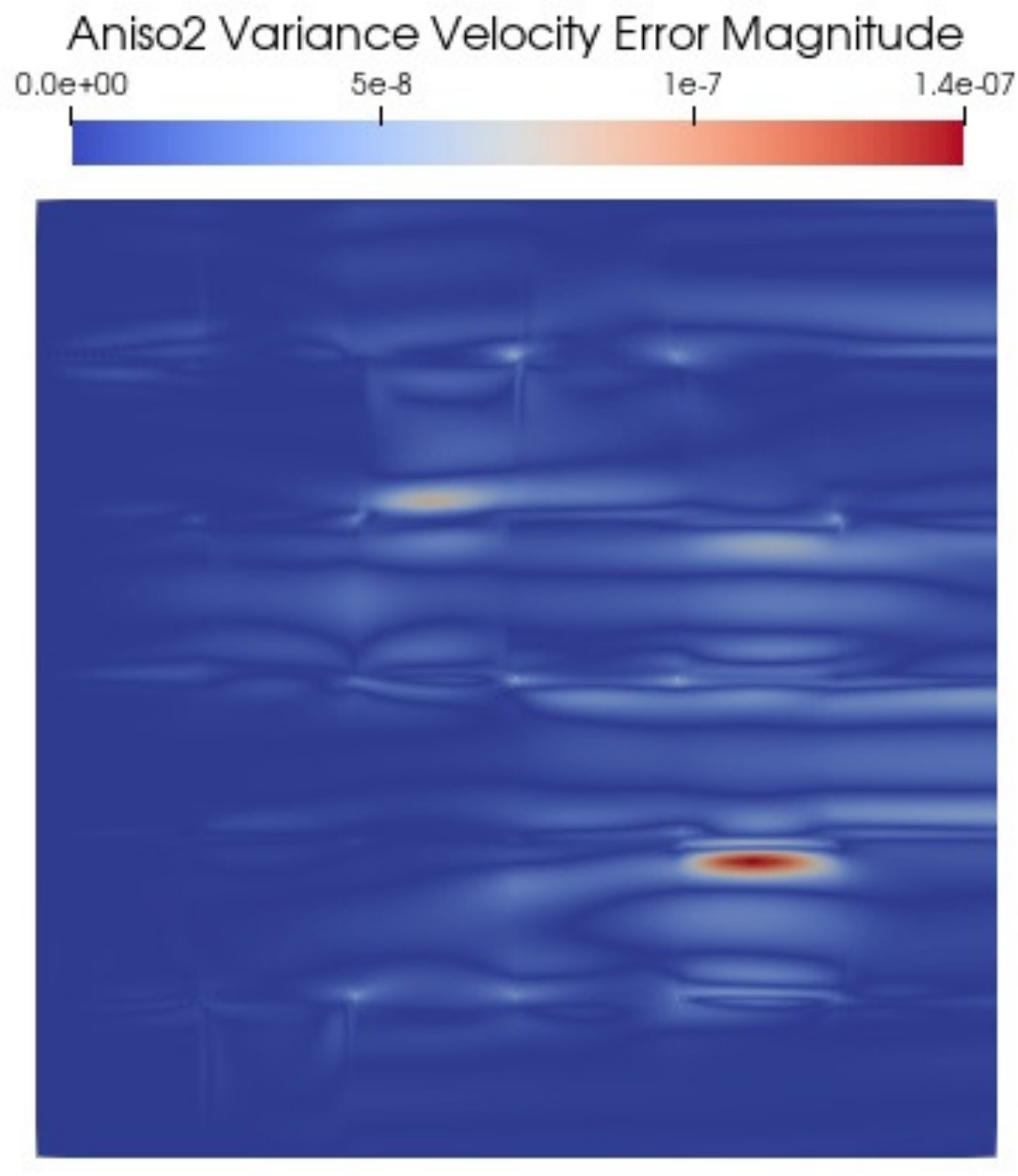}
  \end{tabular}
  \caption{Moment errors in velocity magnitude for each of the three problems}
  \label{fig:velocity-magnitude-mean-and-variance-error}
\end{figure*}

\begin{figure*}[p!]
  \centering
  \begin{tabular}{ccc}
    \includegraphics[width=5cm]{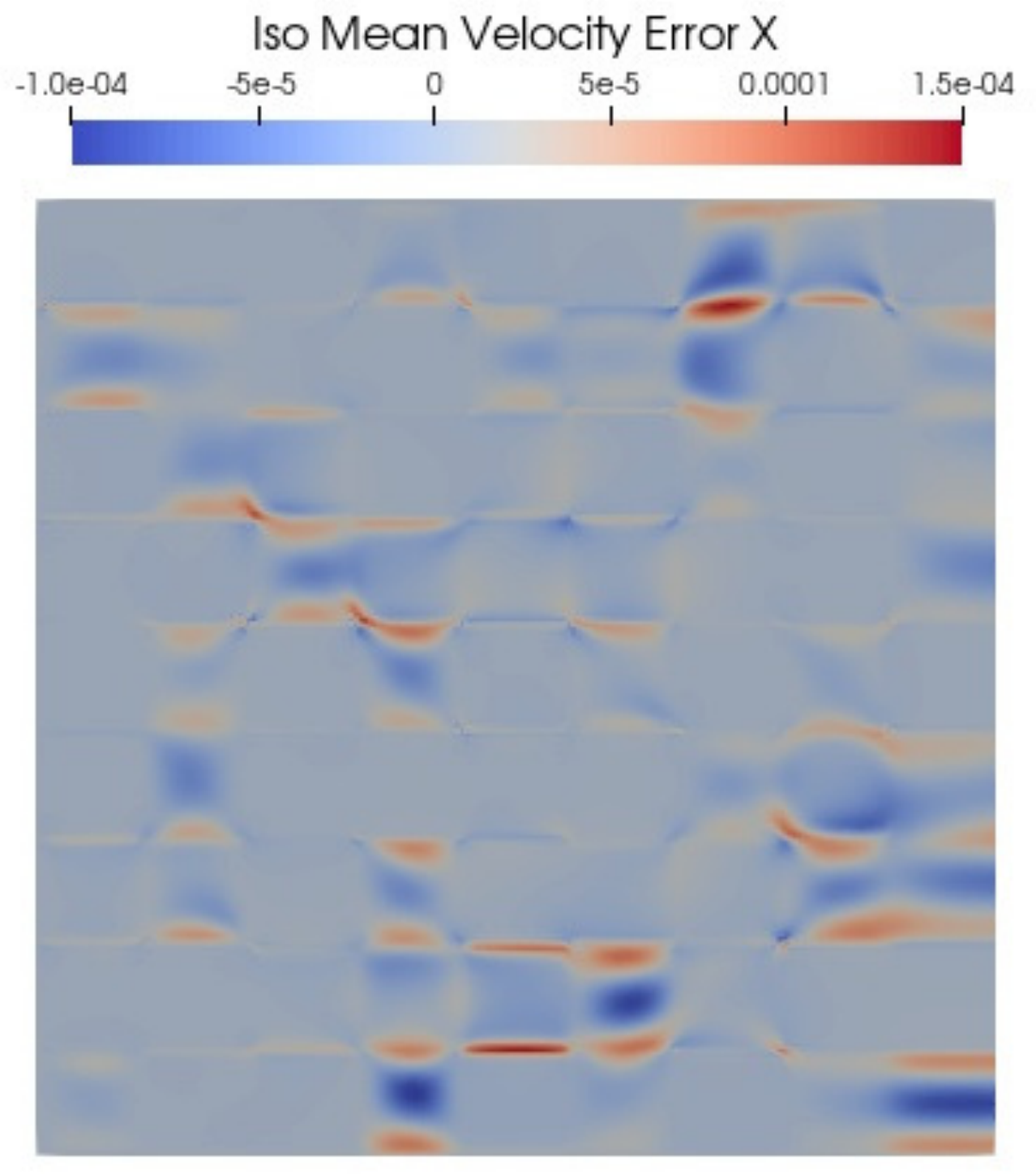} &
    \includegraphics[width=5cm]{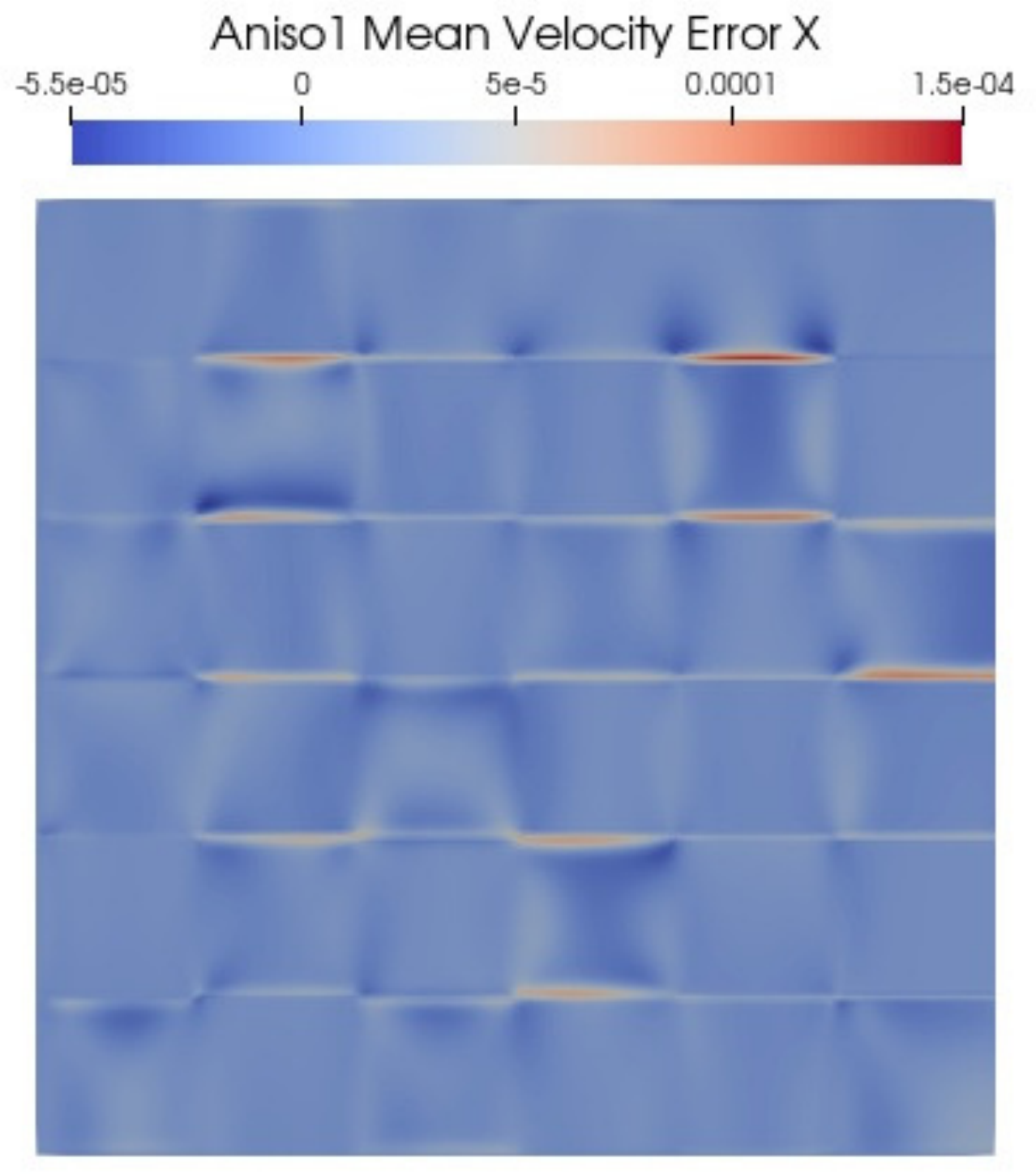} &
    \includegraphics[width=5cm]{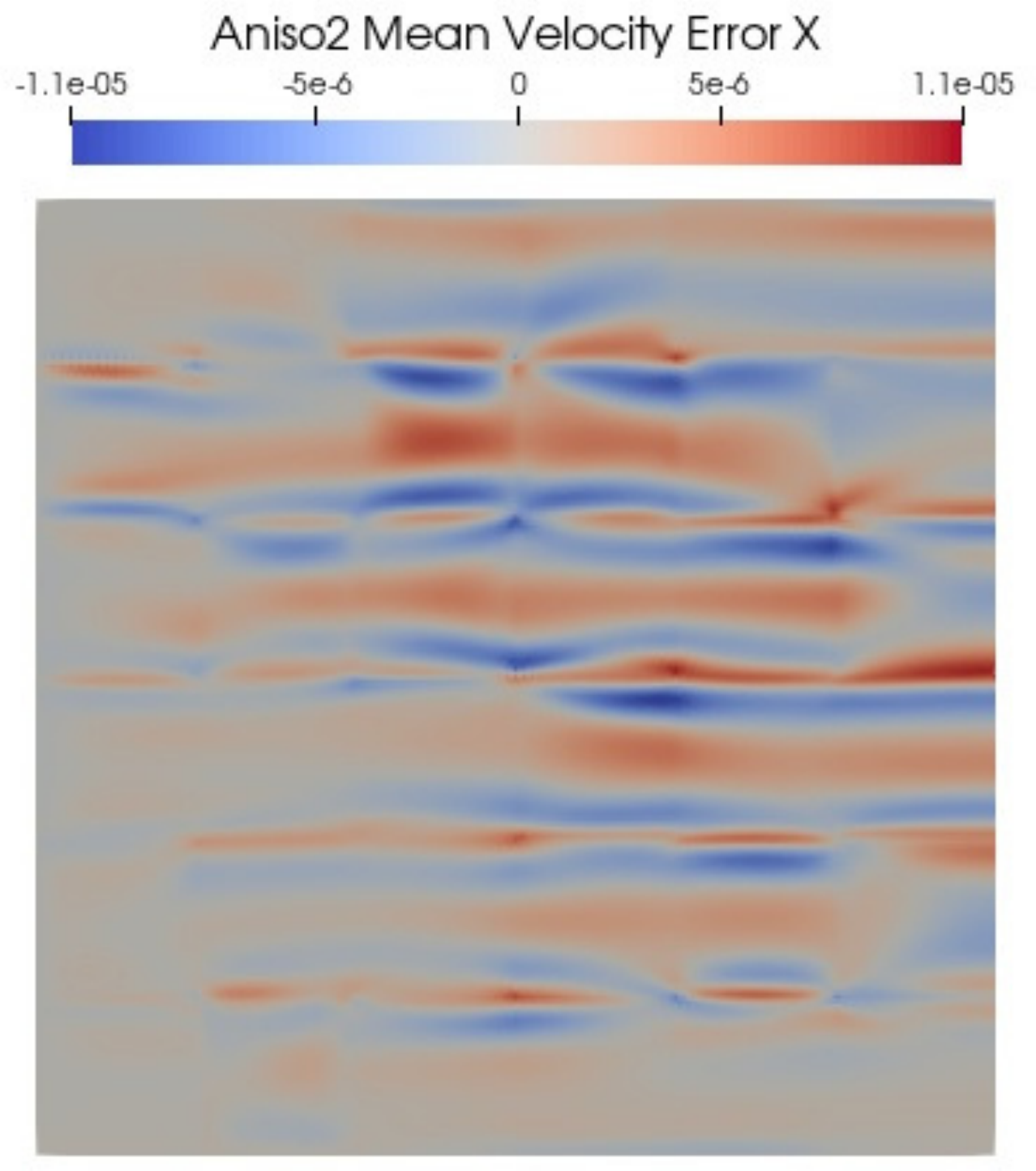} \\
    \includegraphics[width=5cm]{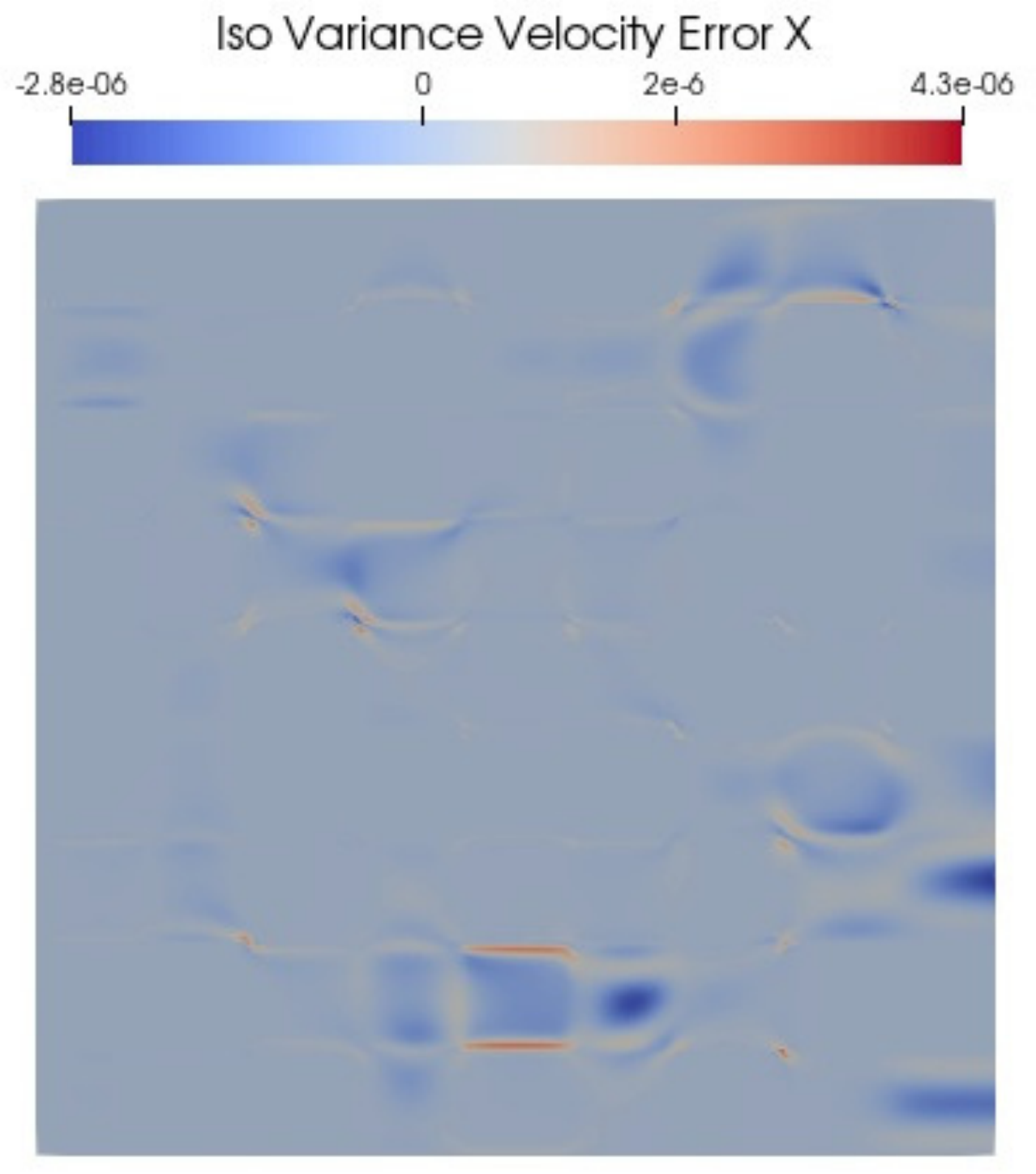} &
    \includegraphics[width=5cm]{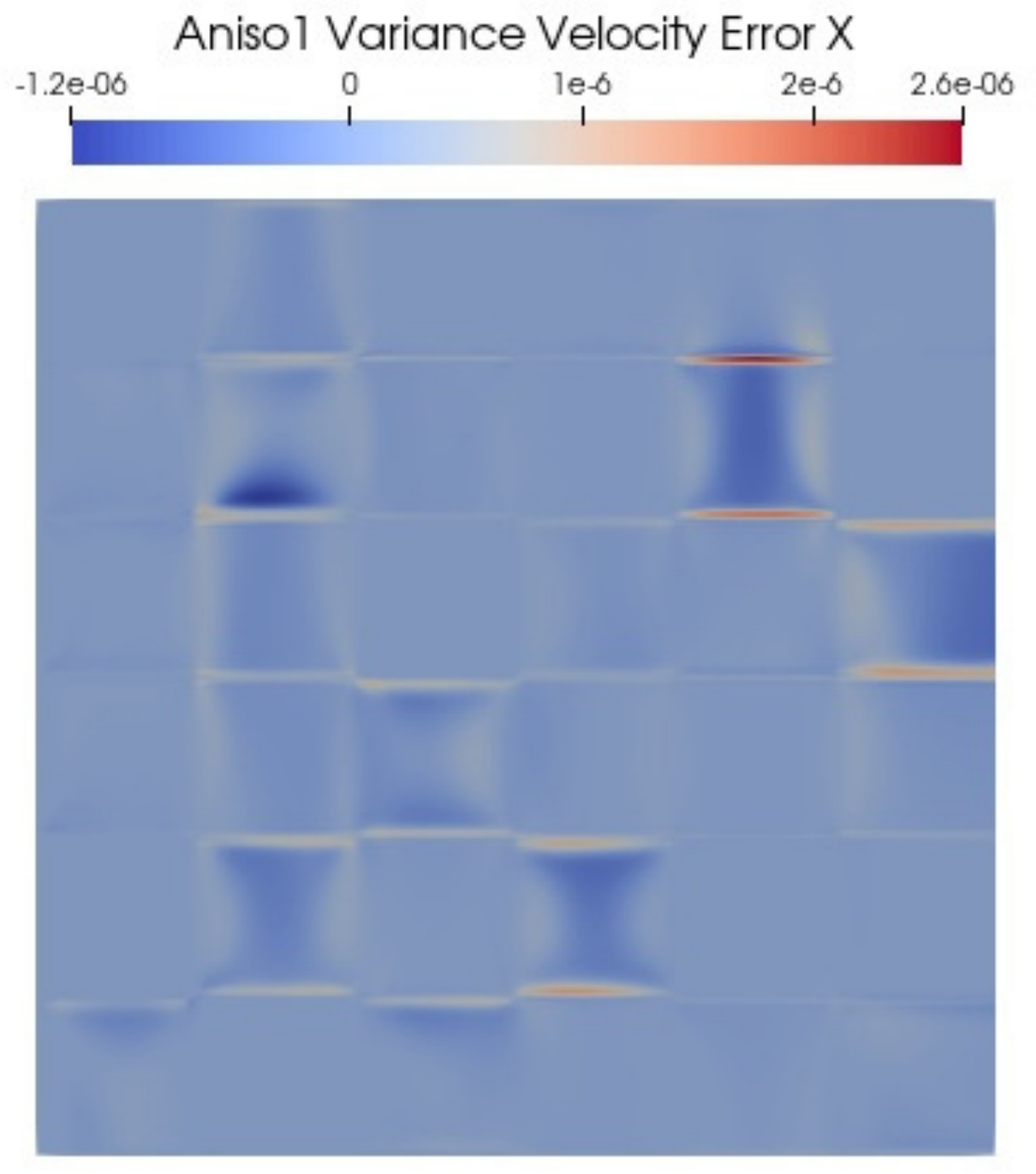} &
    \includegraphics[width=5cm]{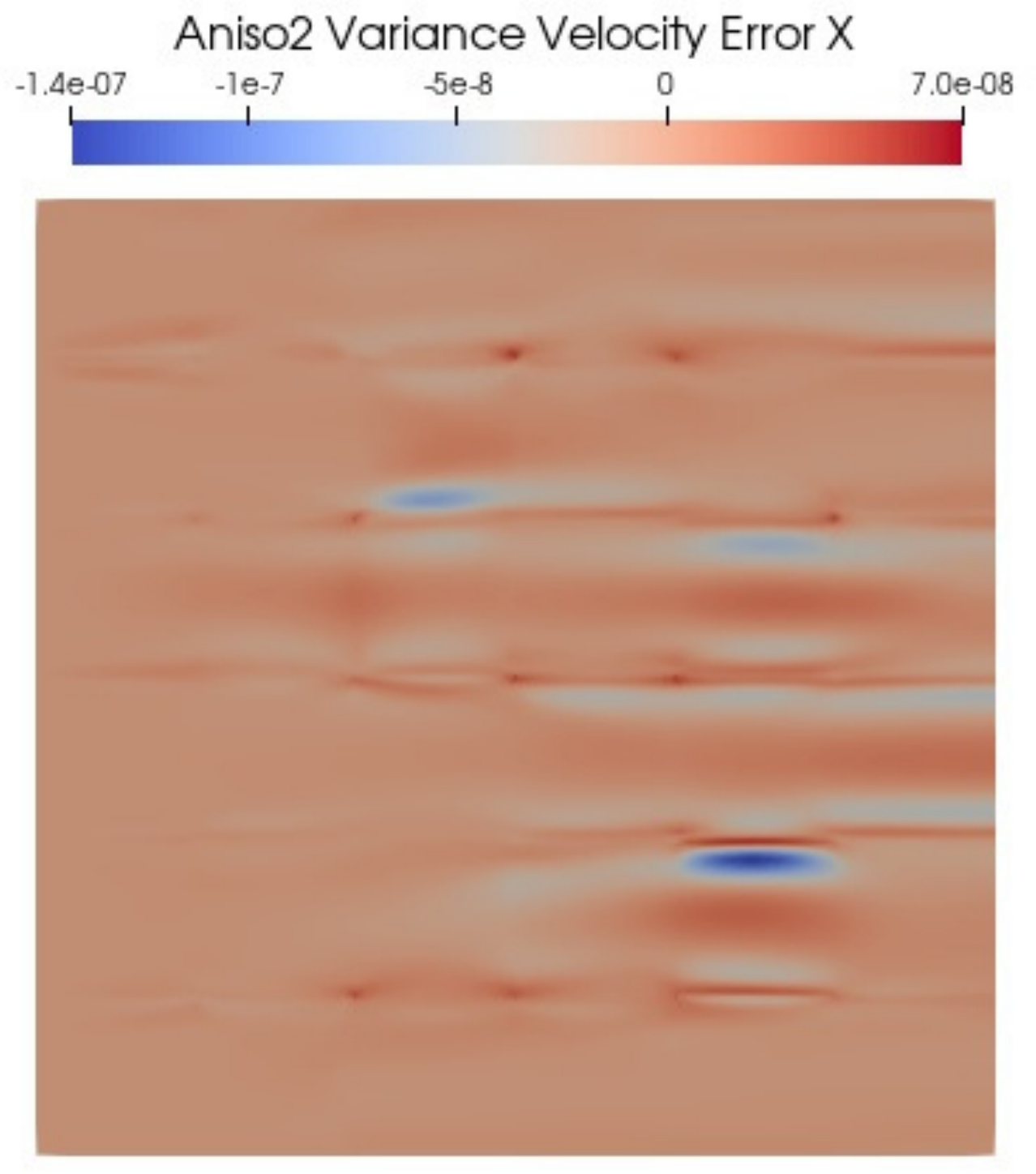}
  \end{tabular}
  \caption{Moment errors in the~$x$ velocity for each of the three problems}
  \label{fig:x-velocity-mean-and-variance-error}
\end{figure*}

\begin{figure*}[p!]
  \centering
  \begin{tabular}{ccc}
    \includegraphics[width=5cm]{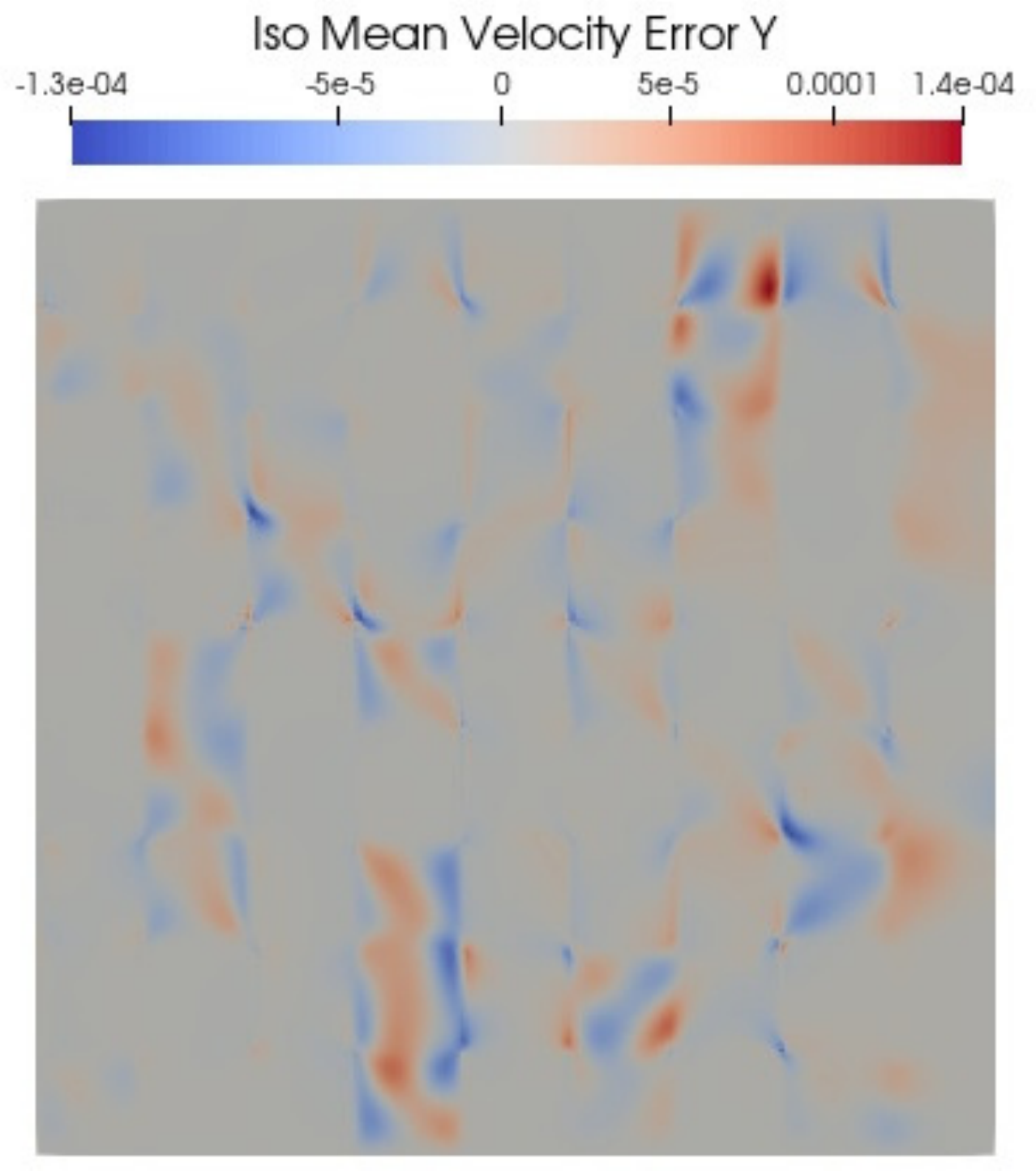} &
    \includegraphics[width=5cm]{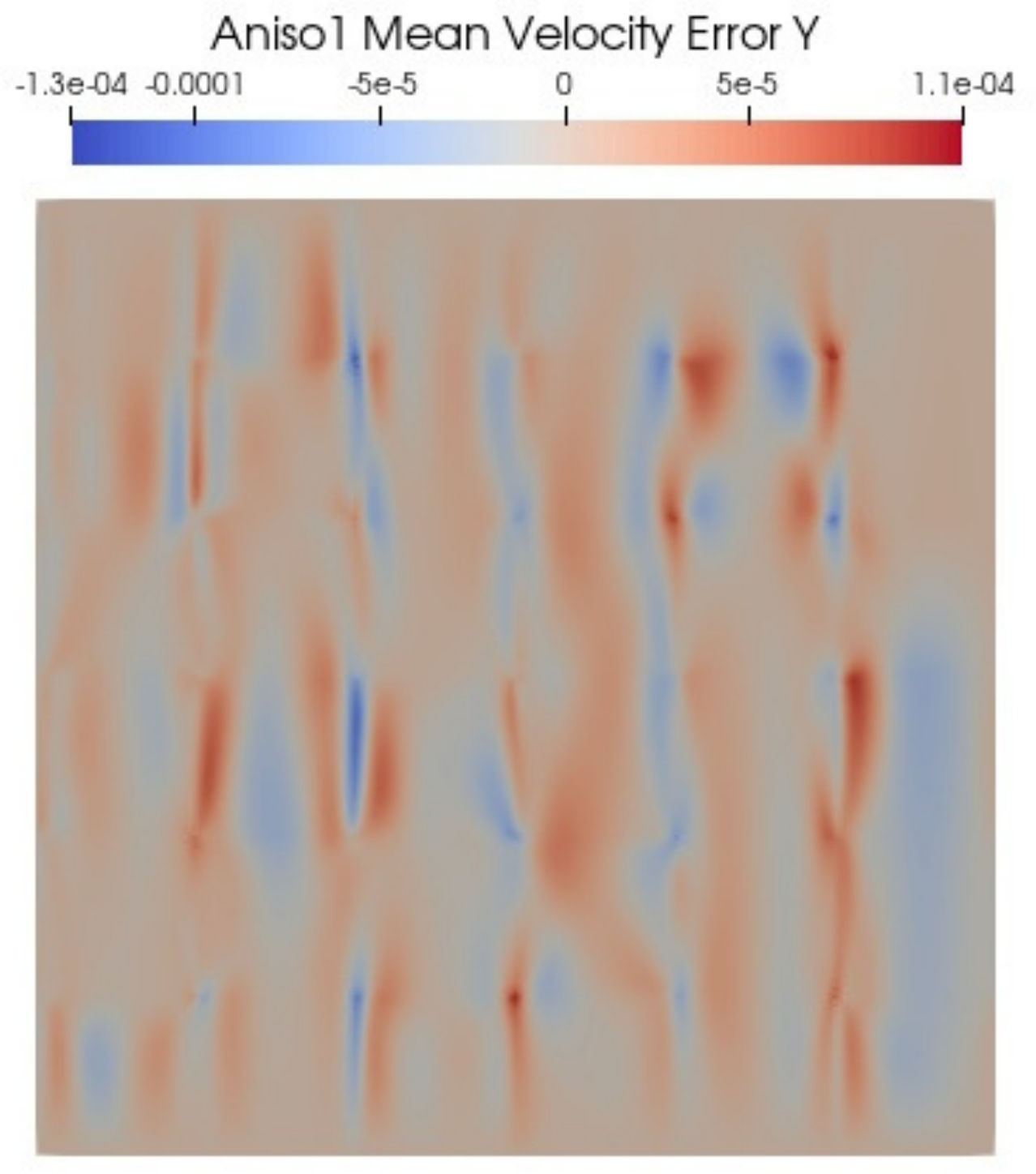} &
    \includegraphics[width=5cm]{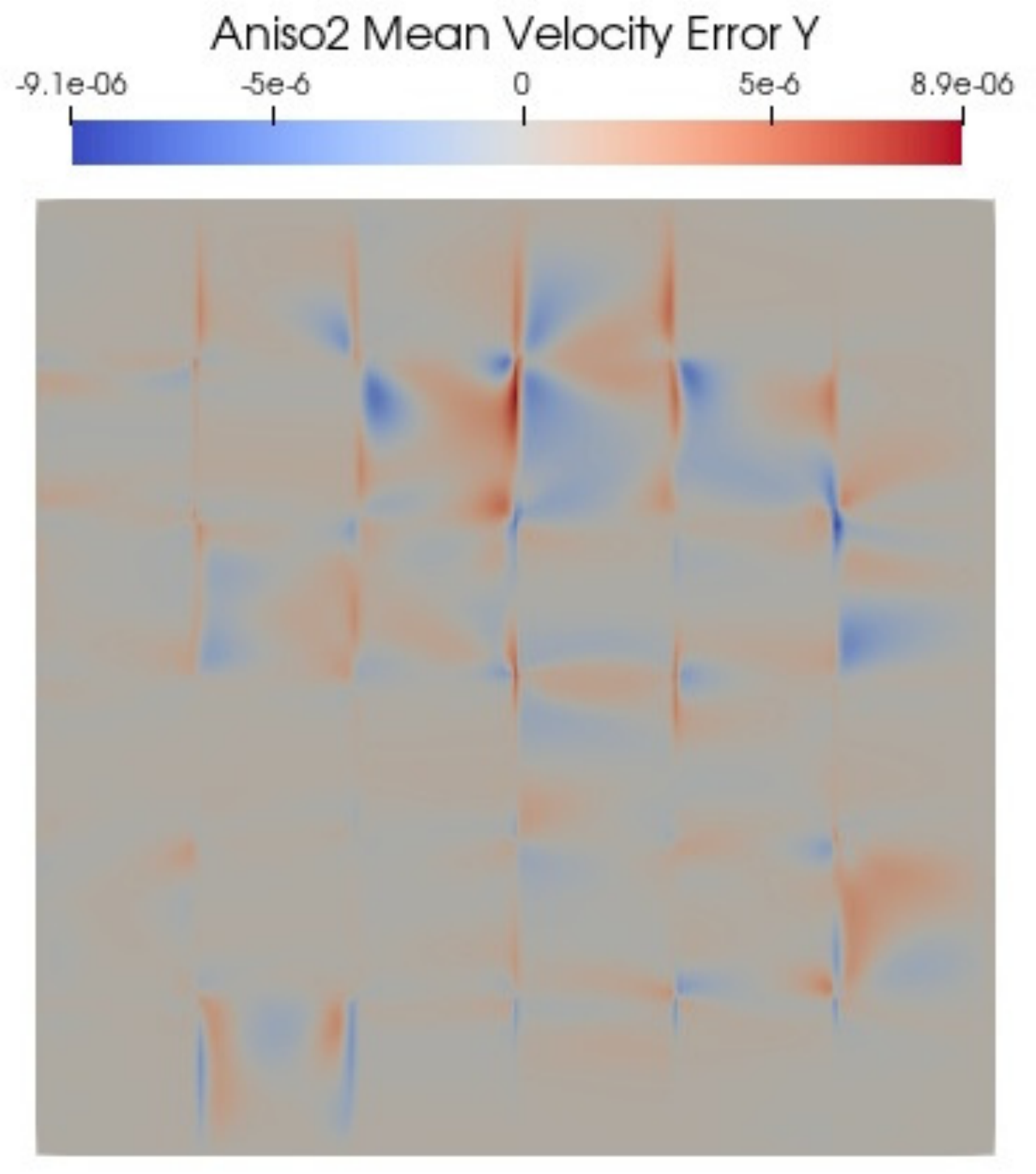} \\
    \includegraphics[width=5cm]{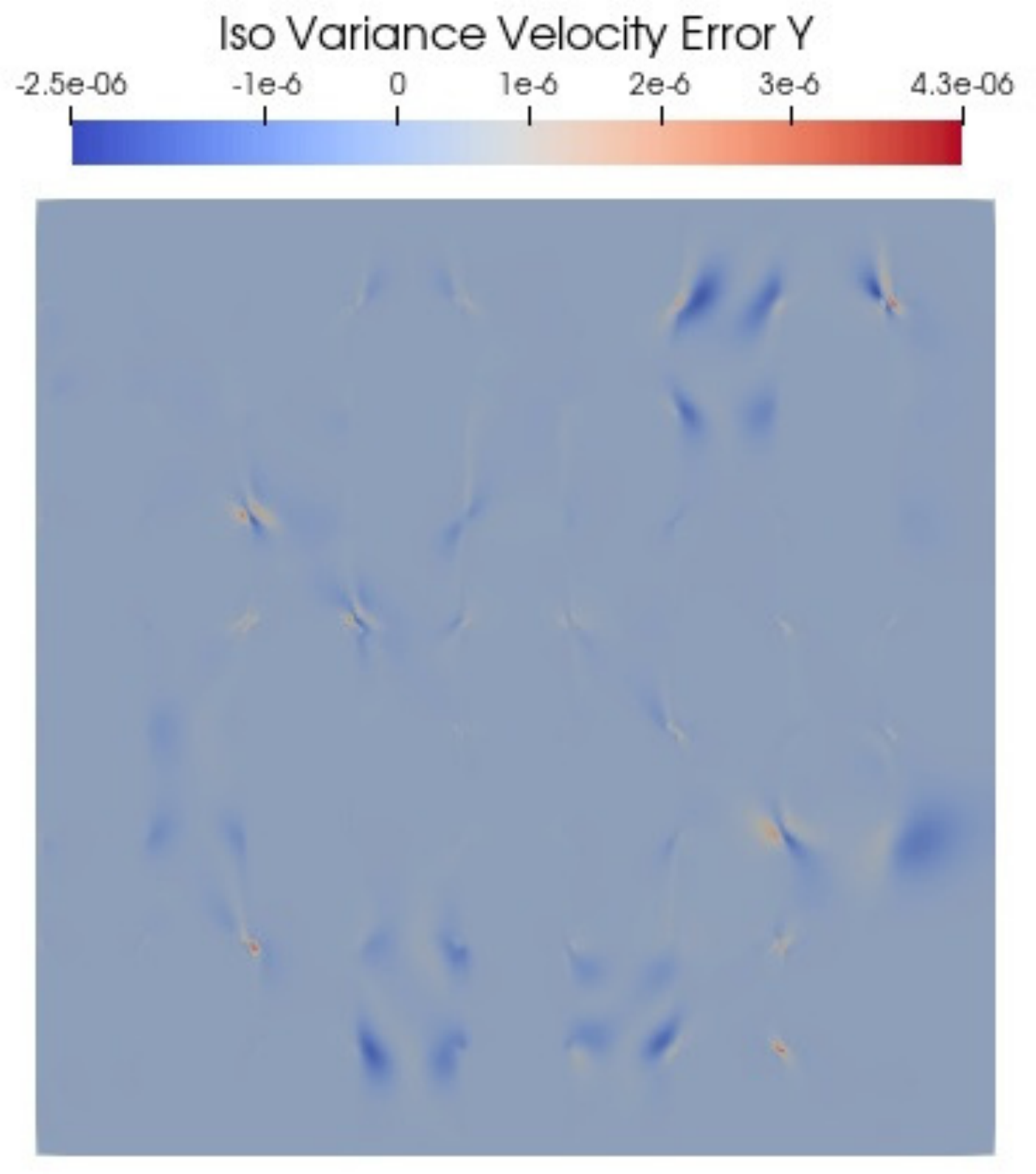} &
    \includegraphics[width=5cm]{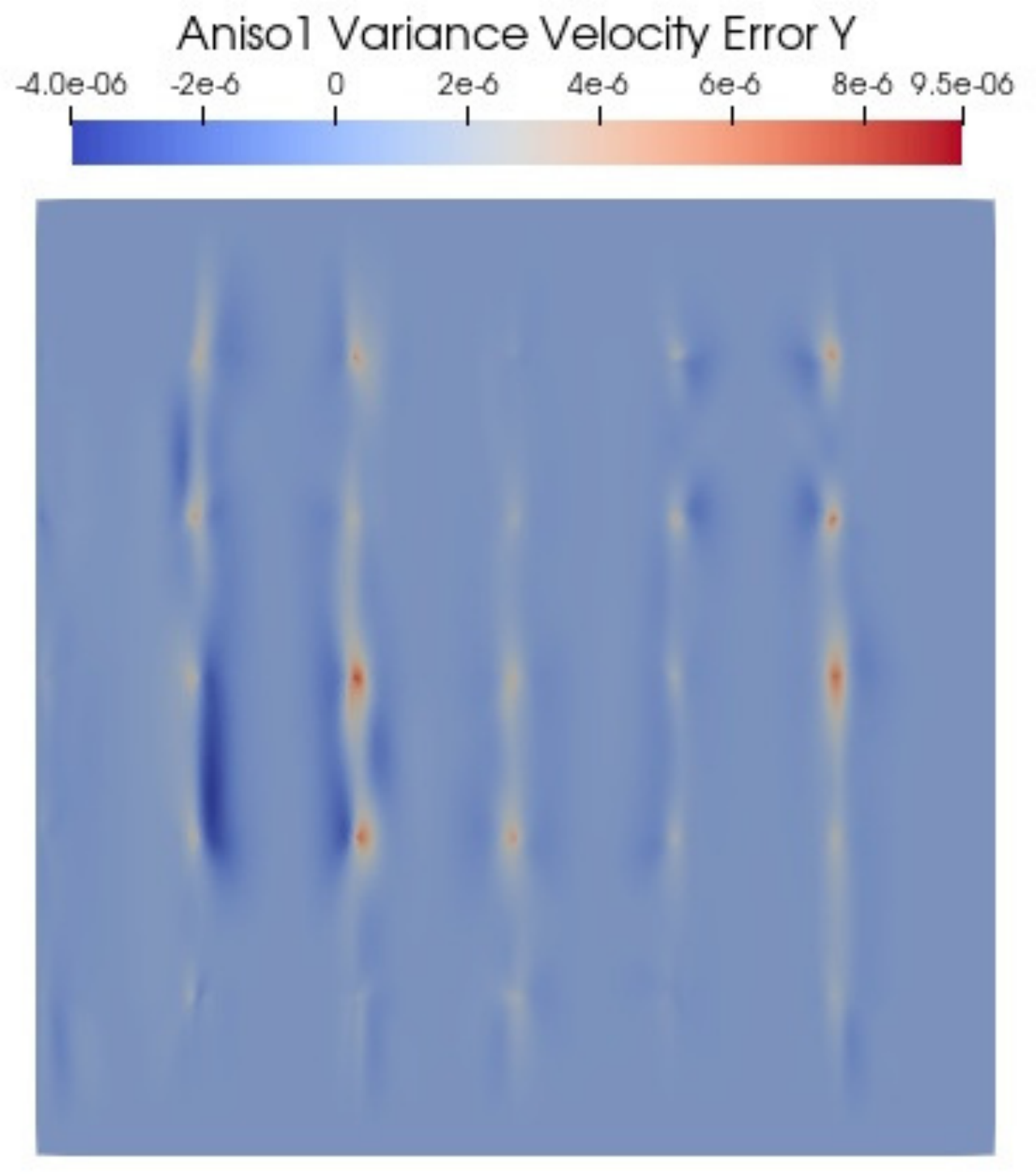} &
    \includegraphics[width=5cm]{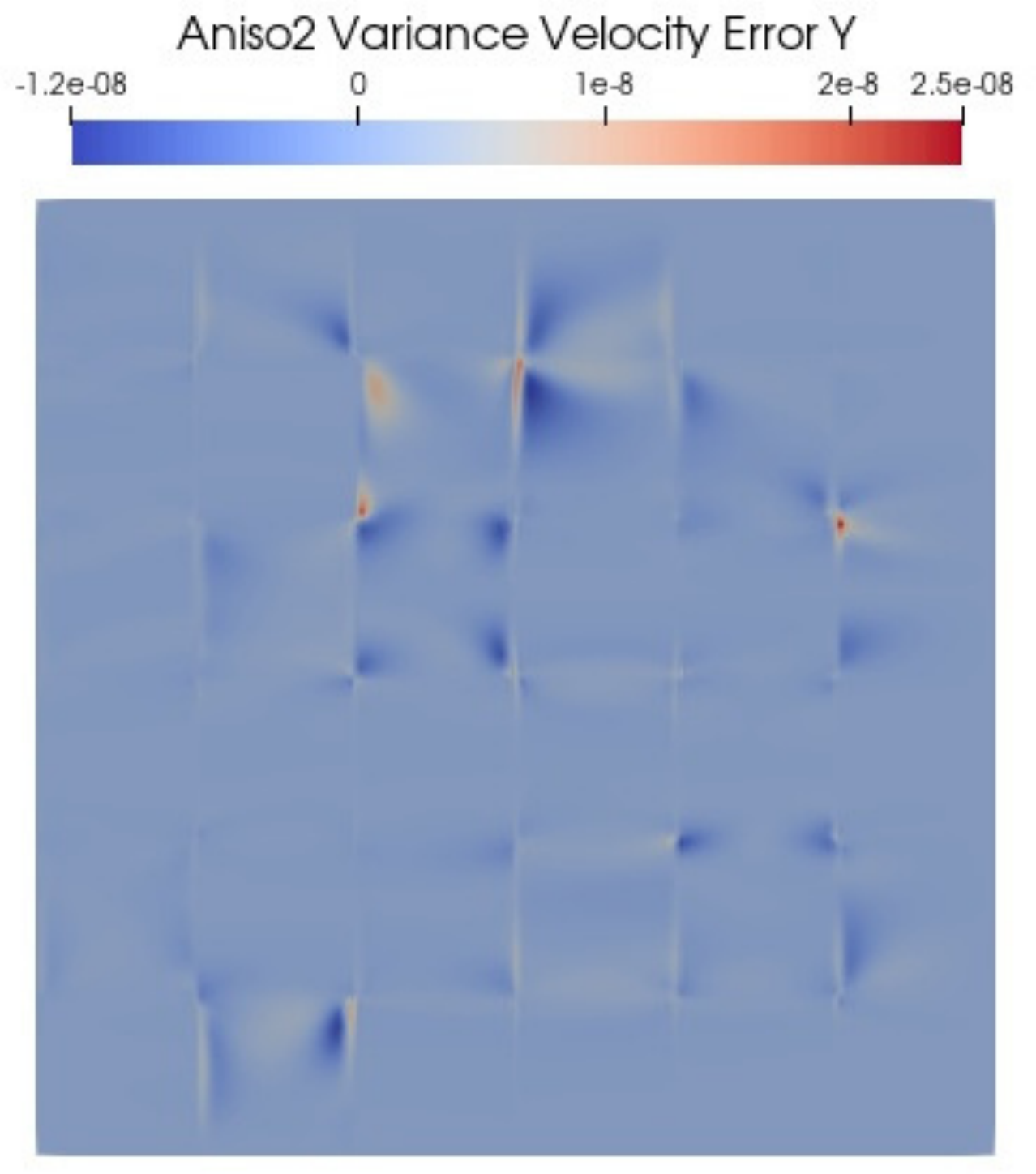}
  \end{tabular}
  \caption{Moment errors in the~$y$ velocity for each of the three problems}
  \label{fig:y-velocity-mean-and-variance-error}
\end{figure*}

The errors over the physical domain~$\domain$ are displayed in Figures~\ref{fig:pressure-mean-and-variance-error} (pressure),~\ref{fig:velocity-magnitude-mean-and-variance-error} (velocity magnitude),~\ref{fig:x-velocity-mean-and-variance-error} ($x$ velocity), and~\ref{fig:y-velocity-mean-and-variance-error} ($y$ velocity). All of these images were generated using ParaView~\cite{Paraview} from the computations with a reduced basis tolerance of~$0.01$ and ANOVA tolerance of~$10^{-6}$. Each figure consists of six subfigures arranged in a 2x3 grid. The columns correspond to the problems \emph{iso}, \emph{aniso1}, and \emph{aniso2}, in that order. The first row depicts the mean error and the second row depicts the variance error.

The pressure errors are depicted in Figure~\ref{fig:pressure-mean-and-variance-error}. The errors exhibit different behavior among the three problems. For the \emph{iso} problem, the errors appear to be evenly distributed throughout the domain although with clear increases at the corners of subdomains. However, for the \emph{aniso1} (resp., \emph{aniso2}), problem, we observe vertical (resp., horizontal) bands. Recall that the chosen permeabilities for the anisotropic problems favor either vertical or horizontal flow. The bands appear to reflect these favored permeabilities, as higher errors might be expected where values are greater. As with the \emph{iso} case, interfaces and corners between subdomains are emphasized. The heightened errors at these interfaces are likely due to higher order effects not captured in the ANOVA terms selected for the expansion.

The velocity errors are depicted in Figures~\ref{fig:velocity-magnitude-mean-and-variance-error},~\ref{fig:x-velocity-mean-and-variance-error}, and~\ref{fig:y-velocity-mean-and-variance-error}. As with the pressures errors, the errors in all cases are heightened along subdomain interfaces, as expected due to the exclusion of higher order ANOVA terms. Furthermore, the favored permeabilities for the anisotropic problems result in larger errors in the velocity directions which are favored, again likely due to the values being larger.

While the above analysis demonstrates that a smaller reduced basis tolerance greatly improves the accuracy of the moment computations, it comes with a cost of a larger number of high-fidelity solves and a larger reduced basis. Table~\ref{tab:rb-sizes} summarizes the number of high-fidelity solves for each problem and each choice of tolerances. For comparison, the number of collocation points is provided in parentheses. For all three problems, it appears that the number of high-fidelity solves almost doubles as the reduced basis tolerance is reduced by a factor of 10. However, the number of high-fidelity solves remains fairly constant even as the ANOVA tolerance increases. This suggests that a small reduced basis is able to fairly accurately approximate the space of high-fidelity solves.

\begin{table}[h!]
  \centering
  \begin{tabular}{lr|ccc|c}
    \cline{3-5}
    \, & \, & \multicolumn{3}{c|}{$\epsilon^{\mathrm{A}}$} & \, \\
    \, & \, & $10^{-4}$ & $10^{-5}$ & $10^{-6}$ & \, \\
    \hline
    \multicolumn{1}{ |l }{\multirow{3}{*}{$\epsilon^{\mathrm{RB}}$}} & 1 &	28 (1573) &	33 (27701) &	33 (43541) & \multicolumn{1}{c|}{\multirow{3}{*}{iso}} \\
    \multicolumn{1}{ |l }{} & 0.1 & 59 (4021)	& 61 (33605) &	62 (50885) & \multicolumn{1}{c|}{} \\
    \multicolumn{1}{ |l }{} & 0.01 & 110 (7285) &	122 (43541) &	123 (52165) & \multicolumn{1}{c|}{} \\
    \hline
    \multicolumn{1}{ |l }{\multirow{3}{*}{$\epsilon^{\mathrm{RB}}$}} & 1 & 	29 (1537) &	31 (19889) &	31 (37825) & \multicolumn{1}{c|}{\multirow{3}{*}{aniso1}} \\
    \multicolumn{1}{ |l }{} & 0.1 & 50 (2737) &	52 (27665) &	52 (40049) & \multicolumn{1}{c|}{} \\
    \multicolumn{1}{ |l }{} & 0.01 & 90 (3985) &	95 (36737) &	96 (41185) & \multicolumn{1}{c|}{} \\
    \hline
    \multicolumn{1}{ |l }{\multirow{3}{*}{$\epsilon^{\mathrm{RB}}$}} & 1 & 	34 (2737) &	37 (35665) &	37 (41185) & \multicolumn{1}{c|}{\multirow{3}{*}{aniso2}} \\
    \multicolumn{1}{ |l }{} & 0.1 & 73 (2737) &	77 (34609) &	77 (41185) & \multicolumn{1}{c|}{} \\
    \multicolumn{1}{ |l }{} & 0.01 & 	91 (3329) &	120 (34609) &	121 (41185) & \multicolumn{1}{c|}{} \\
    \hline
  \end{tabular}
  \caption{Comparison of number of high-fidelity solves used to form the reduced basis and number of ANOVA collocation points (in parentheses) for the three problems with varying tolerances}
  \label{tab:rb-sizes}
\end{table}

\begin{figure}[h!]
  \centering
  \begin{tabular}{c}
    \includegraphics[width=8.4cm]{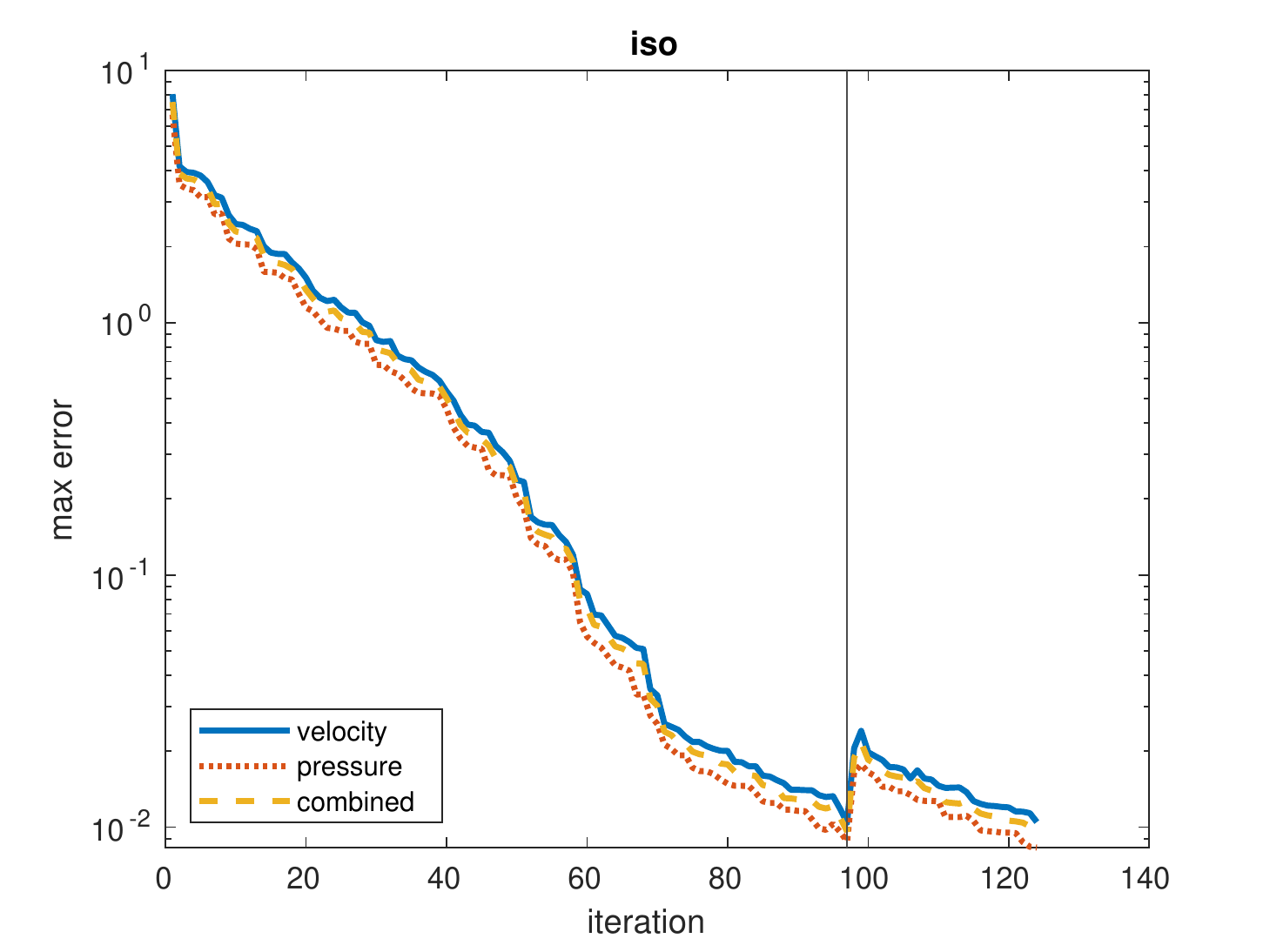} \\
    \includegraphics[width=8.45cm]{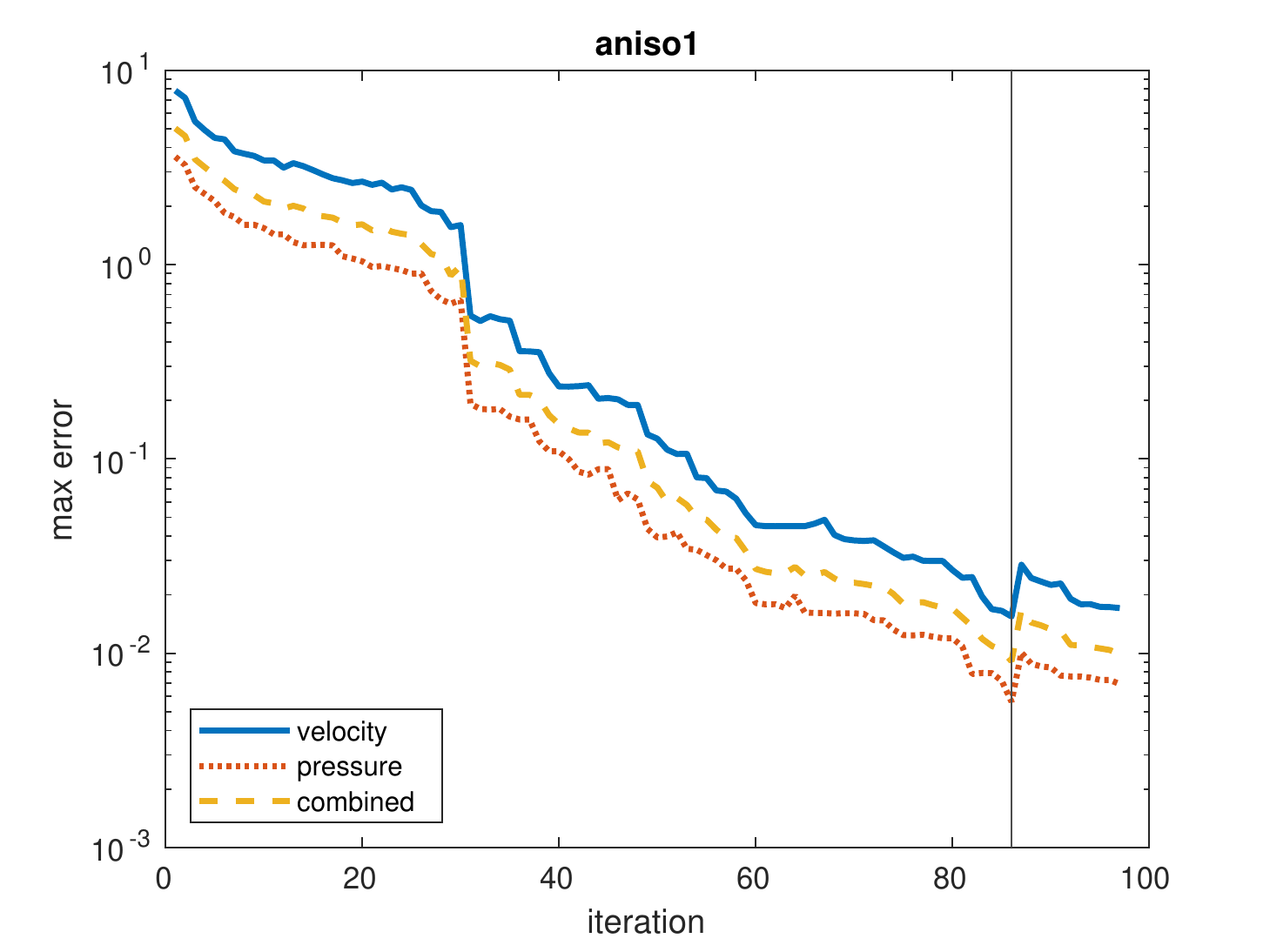} \\
    \includegraphics[width=8.4cm]{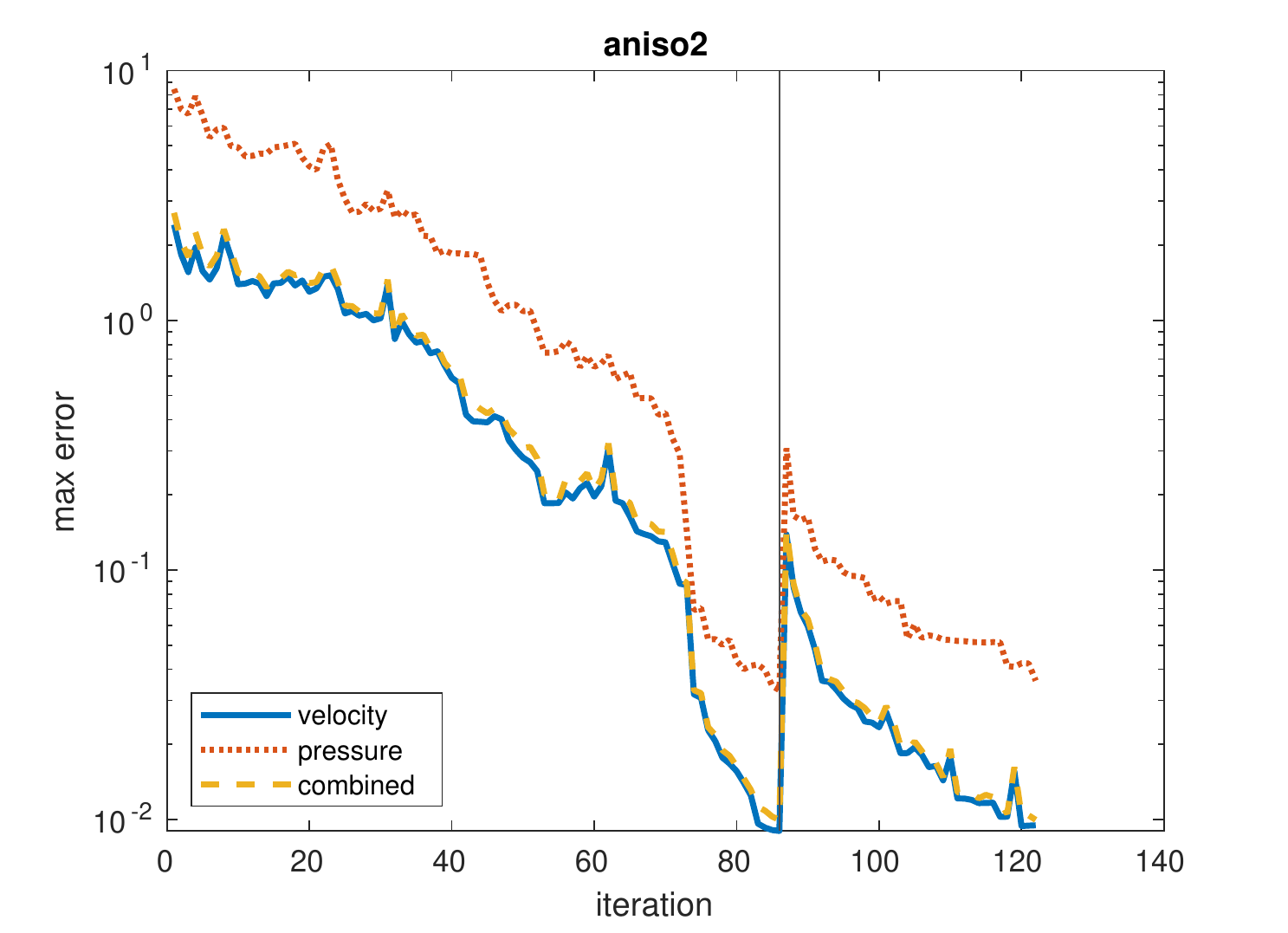}
  \end{tabular}
  \caption{Maximum values of error estimates for each training iteration using~$\epsilon^{\mathrm{A}} = 10^{-6}$ and~$\epsilon^{\mathrm{RB}} = 0.01$}
  \label{fig:max-aee-errors}
\end{figure}

The decrease in the a posteriori error estimates over each iteration of training is depicted in Figure~\ref{fig:max-aee-errors}. The velocity error estimate is depicted as the solid blue line, the pressure error estimate as the red dotted line, and the combined error estimate as the orange dashed line. Each error estimate depicted is normalized by the norm of the reduced basis solution to produce a relative error estimate. The figure displays, at each iteration, the largest relative error estimate of the given type. Since the largest error estimates for each type are not guaranteed to be taken from the same sample at any given iteration, and each of these error estimates are normalized by the reduced basis solution at that sample, there is no guaranteed ordering of the values (i.e., the combined error is not necessarily greater than the velocity or pressure error).

For the \emph{iso} problem, all three errors appear to be roughly equal. For the \emph{aniso1} problem, however, the velocity error is the largest with the combined error being in between velocity and pressure. Recall that the \emph{aniso1} problem favors vertical flow, but the boundary conditions impose horizontal flow. This contention between vertical and horizontal flow perhaps explains the increased difficulty of the reduced basis in approximating the velocity space. On the other hand, the pressure error is the greatest for the \emph{aniso2} problem with the velocity and combined errors being roughly equal. For this problem, horizontal flow is favored so that there is no contention between the flow specified by the boundary conditions and the permeabilities, allowing easier approximation of velocity.

The error estimates are displayed on a log scale, suggesting linear convergence. The black vertical line marks the iteration at which training terminated after the first ANOVA level. The increase in the error estimates at the start of the second ANOVA level is due to an increase in the number of collocation points used during training. In both the \emph{iso} and \emph{aniso1} problems, the increase in error between levels is small, suggesting that the reduced basis formed at the end of the first level generalized fairly well to the new collocation points. The \emph{aniso2} problem has a larger increase in error, but there was also a larger decrease in error towards the end of the level 1 training. This suggests that the level 1 reduced basis did not generalize well to the new collocation points. In all cases, however, the error appears to decrease during the second level at roughly the same rate as during the first level, despite the greater increase in the number of training samples.

The a posteriori error estimates only compute an upper bound on the actual error. It would be helpful to understand how tight these error estimates actually are. To study this, we selected a random subset of 1000 ANOVA collocation points from each problem using tolerances~$\epsilon^{\mathrm{RB}} = 0.01$ and~$\epsilon^{\mathrm{A}} = 10^{-6}$ and computed high-fidelity solutions at each selected point. Using the associated reduced basis for each problem, we also computed the reduced basis approximation. From this, we obtained exact errors for each collocation point and computed ratios of the a posteriori error estimate against the exact error. This ratio is called the \emph{effectivity} of the error estimator. These values are presented in Figure~\ref{fig:effectivities}. Velocity ratios are shown in blue upward facing triangles, pressure ratios in red downward facing triangles, and combined ratios in orange circles. Since the error estimates form upper bounds, we expect these ratios to be greater than 1. However, if the bound is tight, we expect the ratios to be close to 1. From the figure, it is apparent that all ratios are greater than 1, suggesting our error estimates are true upper bounds. However, the ratio is much greater than 1, suggesting that the upper bounds are not as tight as we would like them to be for an ideal setup. It is well known that the tightness of the error estimate is related to the condition number~\cite[Section 3.6.2]{Quarteroni-2016-RBM}. As demonstrated in Table~\ref{tab:condition-numbers}, the Stokes-Brinkman systems are highly ill-conditioned, resulting in loose upper bounds with the error estimates. Thus, the effectivities of Figure~{\ref{fig:effectivities} are expected. Nonetheless, despite the lack of sharpness in the upper bounds, our experiments provide numerical evidence that combining reduced basis methods with anchored ANOVA yields substantial benefits to estimating the moments of stochastic Stokes-Brinkman problems.

\begin{figure}[h!]
  \centering
  \begin{tabular}{c}
    \includegraphics[width=8.4cm]{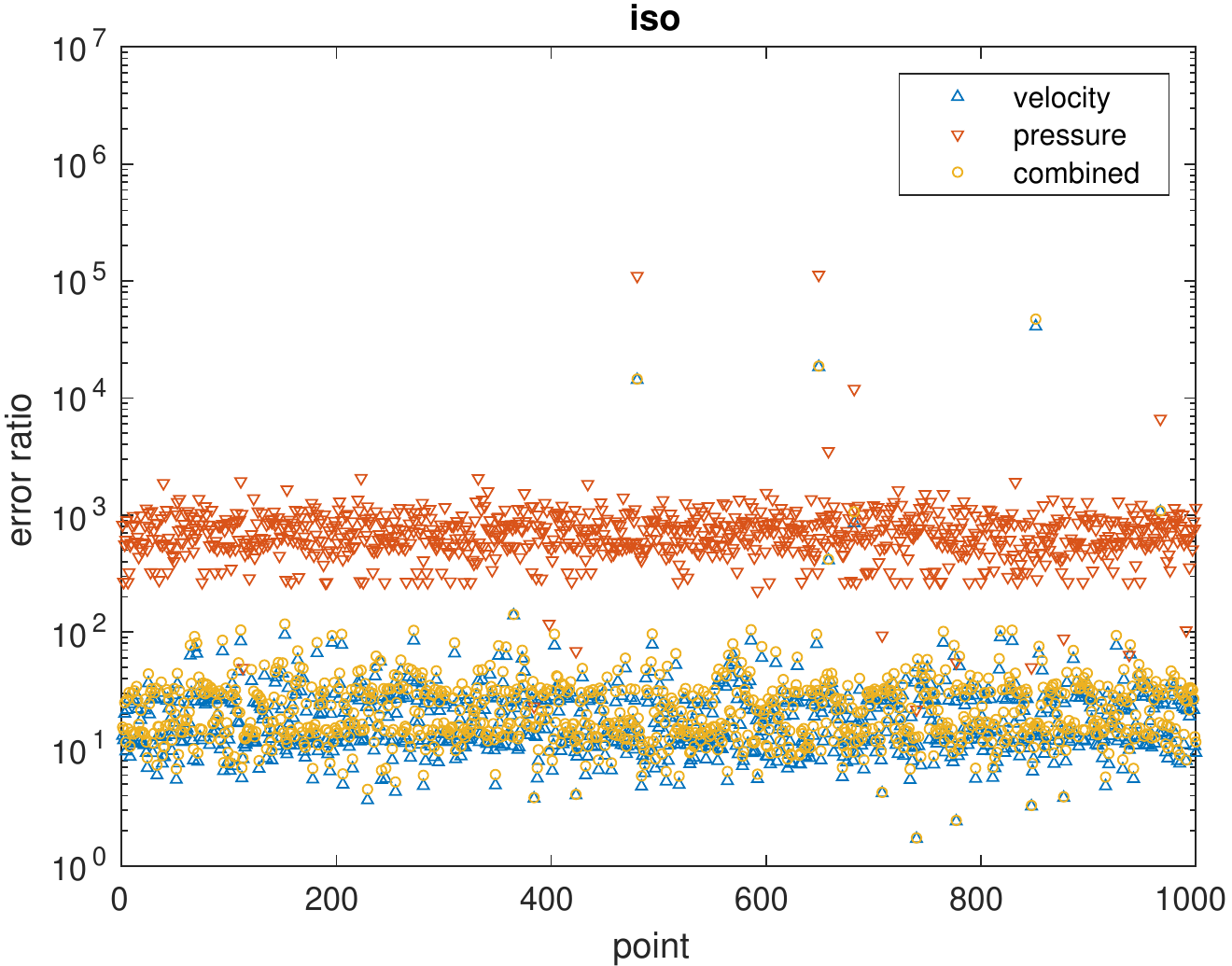} \\
    \includegraphics[width=8.4cm]{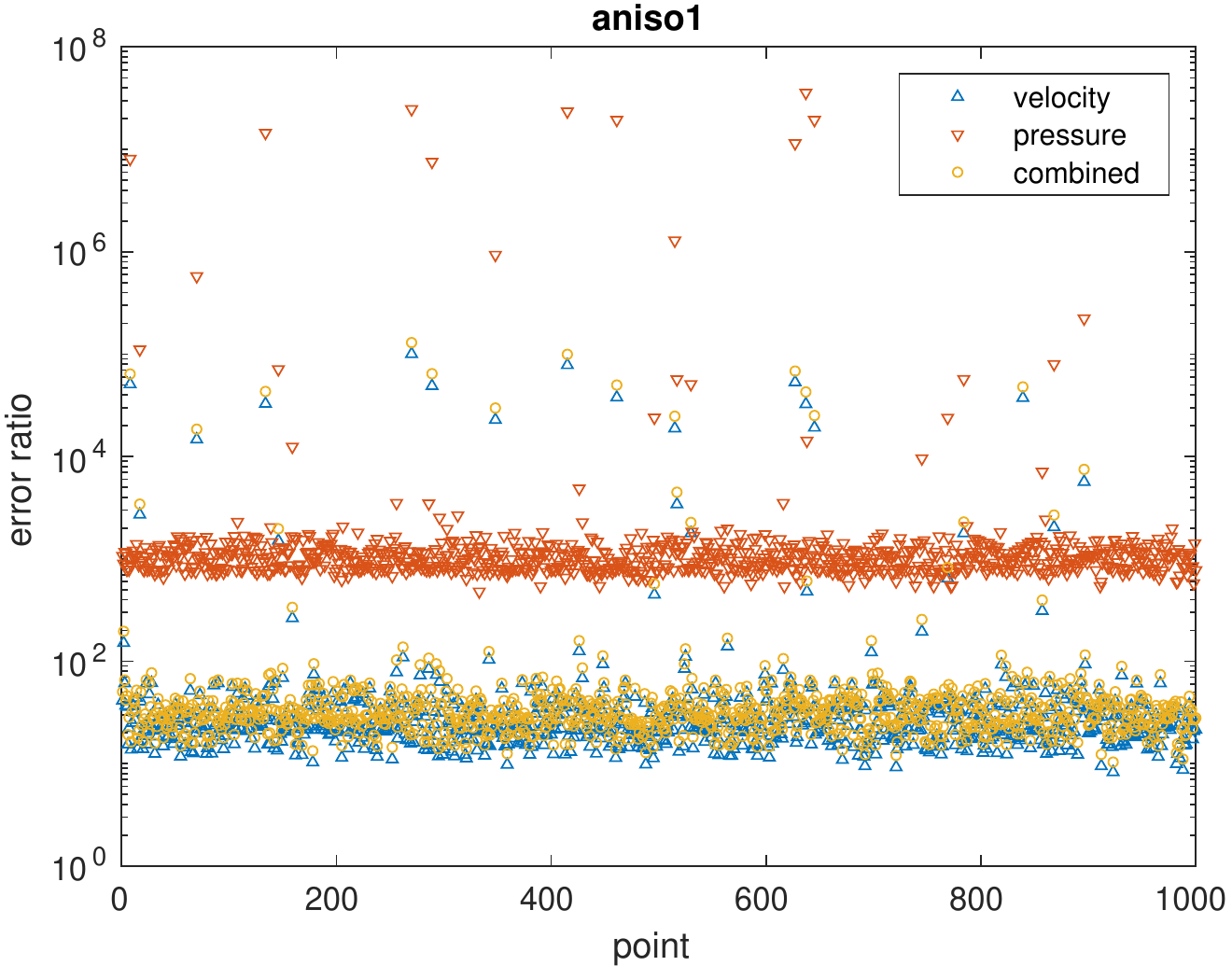} \\
    \includegraphics[width=8.4cm]{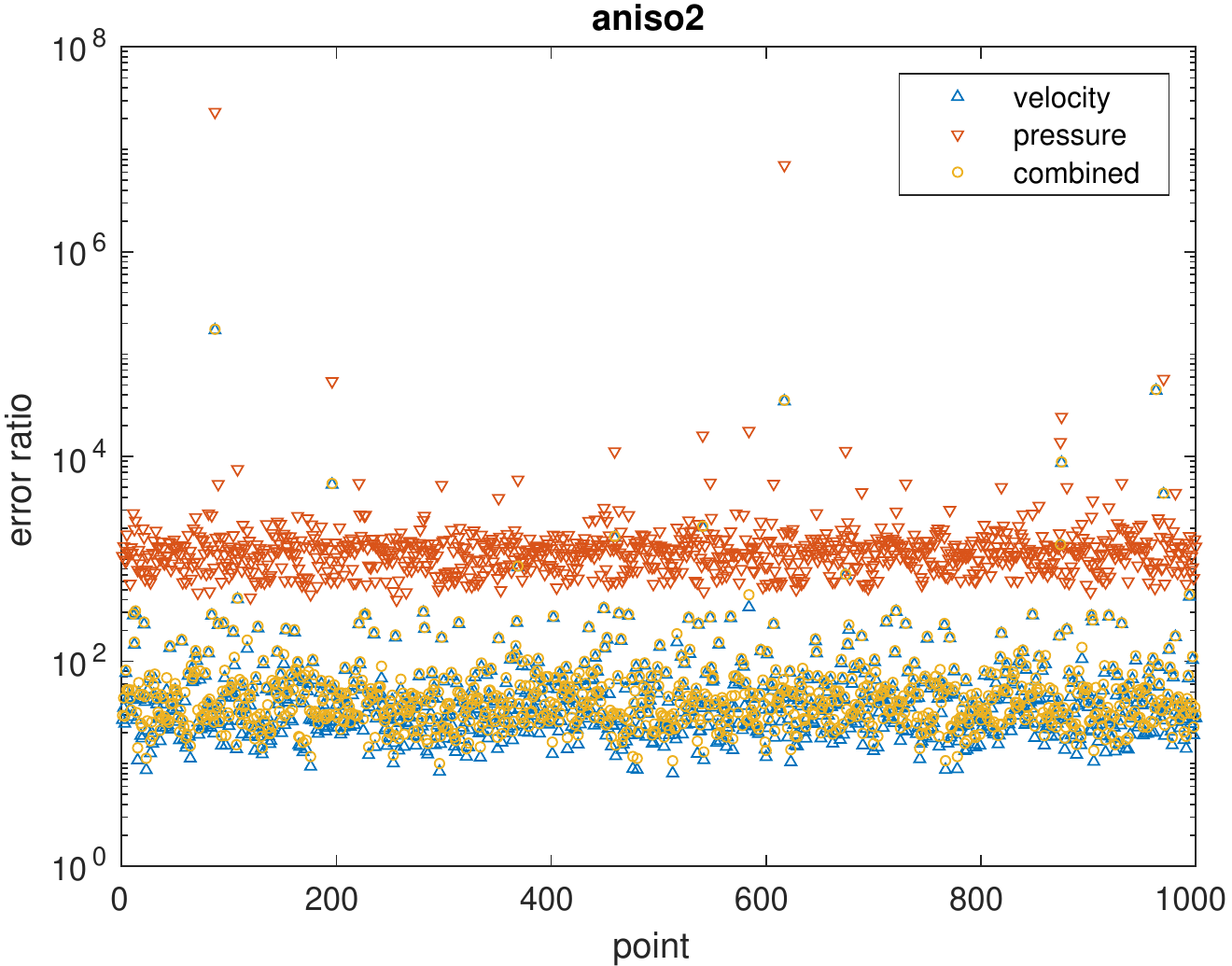}
  \end{tabular}
  \caption{Ratio of the error estimates~(\ref{eq:aee})--(\ref{eq:rb-pressure-bound}) against the exact errors from high-fidelity solves for each of the three problems at 1000 randomly sampled ANOVA collocation points}
  \label{fig:effectivities}
\end{figure}

\section{Conclusion}

We demonstrated that the use of the truncated ANOVA decomposition is effective in reducing the number of collocation points needed for accurate approximation of the statistical moments of several Stokes-Brinkman problems with stochastic permeability. Additionally, we showed that reduced basis methods yield significant savings by reducing the number of high-fidelity solves required to compute these moments. The reduced basis methods rely upon accurate and efficient a posteriori error estimates. We present such estimates based upon the Brezzi stability theory. While these error estimates allowed us to construct small reduced bases, they were not as sharp as desired due to the high-fidelity systems being highly ill-conditioned. Nonetheless, excellent reduced basis approximations were obtained.

\begin{acknowledgements}
  We would like to thank Prof. Howard C. Elman for sharing his notes regarding inf-sup stability for reduced basis methods for saddle-point problems.
\end{acknowledgements}

%
\section*{Conflict of interest}

The authors declare that they have no conflict of interest.

\bibliographystyle{spmpsci}      
\bibliography{Stokes-Brinkman-ANOVA}

\end{document}